\documentclass[11pt,a4paper]{article}
\usepackage{tikz}
\usetikzlibrary{patterns}

\usepackage{epsfig}
\usepackage{graphics}
\usepackage{amsmath,amssymb}
\usepackage{subfigure}
\usepackage{graphicx}
\usepackage{epstopdf}
\usepackage{hyperref}
\usepackage{pdfsync}
\usepackage{color}

\newtheorem{thm}{Theorem}
\newtheorem{lemma}{Lemma}
\newtheorem{proposition}{Proposition}
\newtheorem{definition}{Definition}
\newtheorem{remark}{Remark}

\numberwithin{equation}{section}
\numberwithin{thm}{section}
\numberwithin{lemma}{section}
\numberwithin{proposition}{section}
\numberwithin{definition}{section}
\numberwithin{remark}{section}
\numberwithin{figure}{section}
\numberwithin{table}{section}

\def\Alpha{\boldsymbol{\alpha}}
\def\Beta{\boldsymbol{\beta}}

\title{Minimax Goodness-of-Fit Testing in Ill-Posed Inverse Problems with Partially Unknown  Operators}
\author{{\em Cl\'ement ~Marteau}, \\
         Institut Math\'ematiques de Toulouse,
         INSA de Toulouse,\\
         Universit\'e de Toulouse,\\
         135, avenue de Rangueil, 31 077 Toulouse Cedex 4, France.\\
        {\em Email}:~\texttt{clement.marteau@math.univ-toulouse.fr}
          \\ \\
        and
          \\ \\
        {\em Theofanis ~Sapatinas},\\
        Department of Mathematics and Statistics,\\
        University of Cyprus,\\
        P.O. Box 20537,
        CY 1678 Nicosia,
        Cyprus.\\
        {\em Email}:~\texttt{fanis@ucy.ac.cy}}

\begin{document}
\maketitle

\begin{abstract}
We consider a Gaussian sequence model that contains ill-posed inverse problems as special cases. We assume that the associated operator is  partially unknown in the sense that its singular functions are known and the corresponding singular values are unknown but observed with Gaussian noise. For the considered model, we study the minimax goodness-of-fit testing problem. Working with certain ellipsoids in the space of squared-summable sequences of real numbers, with a ball of positive radius removed, we obtain lower and upper bounds for the minimax separation radius in the non-asymptotic framework, i.e., for fixed values of the involved noise levels. Examples of mildly and severely ill-posed inverse problems with ellipsoids of ordinary-smooth and super-smooth sequences are examined in detail and minimax rates of goodness-of-fit testing are obtained  for illustrative purposes. 


\end{abstract}

\newpage

\section{Introduction}
We consider the following Gaussian sequence model (GSM),
\begin{equation}
\label{1.0.0}
\begin{cases} Y_j=b_j \theta_j+\varepsilon \; \xi_j, & j \in \mathbb{N}, \\
X_j= b_j +\sigma \, \eta_j, & j \in \mathbb{N},
\end{cases}
\end{equation}
where $\mathbb{N} =\{1,2,\ldots \}$ is the set of natural numbers, $b=\{b_j\}_{j \in \mathbb{N}} > 0$ is an {\em unknown} sequence,  $\theta=\{\theta_j\}_{j\in\mathbb{N}} \in l^2(\mathbb{N})$ is the {\em unknown} signal of interest, $\xi=\{\xi_j\}_{j \in\mathbb{N}}$ and $\eta=\{\eta_j\}_{j \in\mathbb{N}}$ are sequences of independent standard Gaussian random variables (and independent of each other), and $\varepsilon, \sigma >0$ are {\em known} parameters (the noise levels). The observations are given by the sequence $(Y,X)=\{(Y_j,X_j)\}_{j \in \mathbb{N}}$ from the GSM (\ref{1.0.0}) and their joint law is denoted by $\mathbb{P}_{\theta,b}$. Here, $l^2(\mathbb{N})$ denotes the space of squared-summable sequence of real numbers, i.e., 
$
l^2(\mathbb{N}) =\{ \theta \in \mathbb{R}^\mathbb{N}:\; \|\theta\|^2:=\sum_{j \in \mathbb{N}}\theta_j^2 < +\infty \}. 
$\\

The GSM (\ref{1.0.0}) arises in the case of ill-posed inverse problems with {\em noisy} operators. Indeed, consider the Gaussian white noise model (GWNM)
\begin{equation}\label{1.0.2}
dY_{\varepsilon}(t)= Af(t)dt + \varepsilon \, dW(t), \quad t \in V,\\
\end{equation}
where $A$ is a linear bounded operator acting on a Hilbert
space ${\cal H}_1$  with values on another Hilbert space ${\cal H}_2$, $f(\cdot) \in {\cal H}_1$ is the unknown response function that one wants to detect or estimate, $W(\cdot)$ is a standard Wiener process on
$V \subseteq \mathbb{R}$, and $\varepsilon>0$ is a known parameter (the noise level). For the sake of simplicity, we only consider the case when $A$
is injective (meaning that $A$ has a trivial nullspace) and assume that $V=[0, 1]$, ${\cal H}_1 = L^2([0,1])$, $U \subseteq \mathbb{R}$ and ${\cal H}_2=L^2(U)$. In most cases of interest, $A$ is a compact operator (see, e.g., Chapter 2 of \cite{EHN_1996}). In particular, it admits a singular value decomposition (SVD) $(b_j, \psi_j, \varphi_j)_{j\in \mathbb{N}}$, in the sense that
\begin{equation}
A \varphi_j=b_j\psi_j, \quad A^{\star}\psi_j=b_j \varphi_j, \quad j \in \mathbb{N},
\label{eq:svdA}
\end{equation}
where $A^{\star}$ denotes the adjoint operator of $A$ -- here $(b^2_j)_{j \in \mathbb{N}}$ and $(\varphi_j)_{j \in \mathbb{N}}$ are, respectively, the eigenvalues and the eigenfunctions of $A^\star A$. Thus, the (first equation in) GSM (\ref{1.0.0}) arises where for all $j\in \mathbb{N}$ 
$$Y_j = \int_{0}^{1} \psi_j(t) dY_{\varepsilon}(t), \quad \theta_j =\int_{0}^{1} \varphi_j(t)f(t)dt,  \quad \xi_j=\int_{0}^{1} \psi_j(t)dW(t), \quad j\in \mathbb{N}, $$  
and $b^2_j>0$ (since $A$ is injective). In this case, the GWNM (\ref{1.0.2}) corresponds to a so-called ill-posed inverse problem since the inversion of $A^*A$ is not bounded. Possible examples of such decompositions arise with, e.g., convolution or Radon-transform operators, see, e.g., \cite{EHN_1996}. The effect of the ill-posedness of the model is clearly seen in the decay of the singular values $b_j$ as $j \to +\infty$. As $j \to +\infty$, $b_j \theta_j$ gets weaker and is then more difficult to perform inference on the sequence $\theta=\{ \theta_j\}_{j\in\mathbb{N}}$.\\

In the early literature, the compact operator $A$ (and, hence, its sequence  $b=\{b_j\}_{j \in \mathbb{N}}$ of singular values) was supposed to be {\em fully known}. (Note that, in this case, the second equation in the GSM (\ref{1.0.0}) does not appear.) We refer, e.g., to \cite{Cavalier_book}, \cite{riskhull}, \cite{cav1},  \cite{cav2}, \cite{cav4} (minimax estimation) and to \cite{LLM_2012}, \cite{ISS_2012} (minimax signal detection/minimax goodness-of-fit testing). Therein, minimax rates/oracle inequalities (estimation) and minimax separation radius/minimax separation rates (signal detection or goodness-of-fit testing) were established, amongst other investigations, for ill-posed inverse problems with smoothness conditions on the sequence of interest. \\

The case of an unknown compact operator $A$ that is observed with Gaussian noise has also been recently treated in the estimation literature, especially the situation where $A$ is {\em partially unknown}, see, e.g., \cite{cavalier2}, \cite{DHPV_2012}, \cite{JS_2013}. In these contributions, it is assumed for the corresponding SVD (\ref{eq:svdA}) that 
\begin{itemize}
\item the sequence of singular functions $(\psi,\varphi)=(\psi_j,\varphi_j)_{j\in \mathbb{N}}$ is {\em known}, 
\item the sequence of singular values $b=\{b_j\}_{j \in \mathbb{N}}$ is {\em unknown}\; but observed with some Gaussian noise.
\end{itemize}
In other words, the following sequence model is considered
\begin{equation*}
\label{eq:rand_eig}
X_j= b_j +\sigma \, \eta_j, \quad j \in \mathbb{N},
\end{equation*}
where $\eta=\{\eta_j\}_{j \in\mathbb{N}}$ is a sequence of independent standard Gaussian random variables (and independent of the standard Gaussian sequence $\xi=\{\xi_j\}_{j \in\mathbb{N}}$), and $\sigma >0$ is a known parameter (the noise level). Therefore, the second equation in the GSM (\ref{1.0.0}) is also readily available.\\

To practically motivate the GSM (\ref{1.0.0}), consider the following deconvolution model (see also \cite{cavalier2} for a complete discussion on this subject)
\begin{equation}
\label{eq:model-dec}
dY_{\varepsilon}(t) = g \star f (t) + \varepsilon\, dW(t ), \quad t \in [0,1], 
\end{equation}
where 
$$
g \star f(t)=\int_{0}^1 g(t-x)f(x)dx, \quad t \in [0,1],
$$ 
is the convolution between $g(\cdot)$ and $f(\cdot)$, $g(\cdot)$ is an unknown 1-periodic (convolution) kernel in $L^2([0,1])$, $f(\cdot)$ is an unknown 1-periodic signal in $L^2([0,1])$, $dY_{\varepsilon}(\cdot)$ is observed, $W(\cdot)$ is a standard Wiener process,  and $\varepsilon >0$ is the noise level. Let $\phi_j(\cdot)$, $j \in \mathbb{N}$, be the usual real trigonometric basis on $V$. The model (\ref{eq:model-dec}) is equivalent to the (first equation in the) GSM (\ref{1.0.0}) by a projection on the trigonometric basis $\phi_j(\cdot)$, $j \in \mathbb{N}$. In the case where the kernel $g(\cdot)$ is unknown (i.e., the sequence $(b_k)_{k\in\mathbb{N}}=(\langle g,\phi_k\rangle)_{k\in \mathbb{N}}$ is unknown), suppose that we can pass the trigonometric basis $\phi_j(\cdot)$, $j \in \mathbb{N}$, through the convolution kernel, i.e., to send {\em each} $\phi_j(\cdot)$, $j \in \mathbb{N}$,
as an input function $f(\cdot)$ and observe the corresponding $dY_{\varepsilon,j}(\cdot)$, $j \in \mathbb{N}$. In other words,we are able to obtain training data for the estimation of the unknown convolution kernel $g(\cdot)$ in this setting. In particular, we obtain exactly the two sequences of observations $Y_j$ and $X_j$, $j \in \mathbb{N}$, in the GSM (\ref{1.0.0}). In this case, the corresponding noise levels coincide, i.e., $\varepsilon=\sigma$.\\

To the best of our knowledge, there is no research work on minimax goodness-of-fit testing in ill-posed inverse problems with partially unknown operators. Our aim is to fill this gap. In particular, considering model (\ref{1.0.0}) and working with certain ellipsoids in the space of squared-summable sequences of real numbers, with a ball of positive radius removed, we obtain lower and upper bounds for the minimax separation radius in the non-asymptotic framework, i.e., for fixed values of $\varepsilon$ and $\sigma$. Examples of mildly and severely ill-posed inverse problems with ellipsoids of ordinary-smooth and super-smooth sequences are examined in detail and minimax rates of goodness-of-fit testing are obtained  for illustrative purposes.\\

The paper is organized as follows. Section \ref{s:stat_model} presents the considered statistical setting and a brief overview of the main results.  Section \ref{sec:sc-o-t} is devoted to the construction of the suggested testing procedure. A general upper bound on the maximal second kind error is then displayed and special benchmark examples are presented for illustrative purposes. The corresponding lower bounds are proposed in Section \ref{lastref}. Some concluding remarks and open questions are discussed in Section \ref{s:remarks}. Finally, all proofs and technical arguments are gathered in Section \ref{s:appendix}.   \\

Throughout the paper we set the following notations. For all $x,y\in \mathbb{R}$, $\delta_x(y)=1$ if $x=y$ and $\delta_x(y)=0$ if $x\not =y$. Also, $x \wedge y:= \min\{x,y\}$ and $x \vee y:= \max\{x,y\}$. Given two sequences $(c_j)_{j \in \mathbb{N}}$ and $(d_j)_{j \in \mathbb{N}}$ of real numbers, $c_j \sim d_j$ means that there exists $0<\kappa_0 \leq \kappa_1 <\infty$ such that $\kappa_0 \leq c_j/d_j \leq \kappa_1$ for all $j \in \mathbb{N}$. Let $\nu$ be either $\varepsilon$ or $\sigma$ or $(\varepsilon,\sigma)$, and let $\mathcal{V}$ be either $\mathbb{R}^+:=(0,+\infty)$ or $\mathbb{R}^+ \times \mathbb{R}^+$. Given two collections $(c_\nu)_{\nu \in \mathcal{V}}$ and $(d_\nu)_{\nu \in \mathcal{V}}$ of positive real numbers, $c_\nu \gtrsim d_\nu$ means that there exists $0 < \kappa_0 <+\infty$ such that $c_\nu \geq \kappa_0 \,d_\nu$ for all $\nu \in \mathcal{V}$. Similarly, $c_\nu \lesssim d_\nu$ means that there exists $0<\kappa_1<+\infty$ such that $c_\nu \leq \kappa_1 \,d_\nu$ for all $\nu \in \mathcal{V}$. 

\section{Minimax Goodness-of-Fit Testing}
\label{s:stat_model}

\subsection{The Statistical Setting}
Given observations $(Y,X)=\{(Y_j,X_j)\}_{j \in \mathbb{N}}$ from the GSM (\ref{1.0.0}), the aim is to compare the underlying (unknown) signal $\theta \in l^2(\mathbb{N})$ to a (known) benchmark signal $\theta_0$, i.e., to test 
\begin{equation}
H_0: \theta=\theta_0 \;\; \mathrm{versus} \;\; \ H_1: \theta-\theta_0 \in \mathcal{F},
\label{testing_pb0}
\end{equation}
for some given $\theta_0$ and a given subspace $\mathcal{F}$. The statistical setting (\ref{testing_pb0}) is known as {\em goodness-of-fit testing} when $\theta_0 \neq 0$ or {\em signal detection} when  $\theta_0 =0$.\\

\begin{remark}
\label{equiv_test}
{\rm Given observations from the GWNM (\ref{1.0.2}), the testing problem (\ref{testing_pb0}) is equivalent to
\begin{equation*}
H_0: f=f_0 \;\; \mathrm{versus} \;\; H_1: f-f_0 \in \tilde{\mathcal{F}},
\label{testing_pbf}
\end{equation*}
for a given benchmark function $f_0$ and a given subspace $\tilde{\mathcal{F}}$. In most cases, $\tilde{\mathcal{F}}$ contains functions $f \in L^2([0,1])$ that admit a Fourier series expansion with Fourier coefficients $\theta$ belonging to ${\cal F}$ (see, e.g., \cite{IS_2003}, Section 3.2.).} 
\end{remark}

The choice of the set $\mathcal{F}$ is important. Indeed, it should be rich enough in order to contain the true $\theta$. At the same time, if it is too rich, it will not be possible to control the performances of a given test due to the complexity of the problem. The common approach for such problems is to impose both a {\em regularity} condition (which characterizes the smoothness of the underlying signal) and an {\em energy} condition (which measures the amount of the underlying signal).  \\ 

Concerning the regularity condition, we will work with certain ellipsoids in $l^2(\mathbb{N})$. In particular, we assume that $\theta\in \mathcal{E}_{a}(R)$, the set $\mathcal{E}_{a}(R)$ being defined as
\begin{equation}
\mathcal{E}_{a}(R) = \left\lbrace \theta\in l^2(\mathbb{N}), \ \sum_{j \in \mathbb{N}} a_j^2 \theta_j^2 \leq R \right\rbrace,
\label{eq:ellips_a}
\end{equation}
where $a=(a_j)_{j\in \mathbb{N}}$ denotes a non-decreasing sequence of positive real numbers with $a_j \rightarrow +\infty$ as $j \rightarrow +\infty$,
and $R>0$ is a constant. The set $ \mathcal{E}_{a}(R)$ can be seen as a condition on the decay of $\theta$. The cases where $a$ increases very fast correspond to $\theta$ with a small amount of non-zero coefficients. In such a case, the corresponding signal can be considered as being `smooth'. Without loss of generality, in what follows, we set $R=1$, and write $\mathcal{E}_{a}$ instead of $\mathcal{E}_{a}(1)$. \\

Regarding the energy condition, it will be measured in the $l^2(\mathbb{N})$-norm. In particular, given $r_{\varepsilon,\sigma}>0$ (called the radius), which is allowed to depend on the noise levels $\varepsilon,\sigma >0$, we will consider $\theta\in \mathcal{E}_{a}$ such that $\| \theta\| > r_{\varepsilon,\sigma}$.  Given a smoothness sequence $a$ and a radius $r_{\varepsilon,\sigma}>0$, the set $\mathcal{F}$ can thus be defined as
\begin{equation}
\label{def:f-set}
\mathcal{F} := \Theta_a(r_{\varepsilon,\sigma}) = \left\lbrace \theta \in \mathcal{E}_a, \ \|\theta \| \geq r_{\varepsilon,\sigma} \right\rbrace.
\end{equation}
In other words, the set ${\cal F}$ is an ellipsoid in $l^2(\mathbb{N})$ with  a ball of radius $r_{\varepsilon,\sigma} >0$ removed. In many cases of interest, the set $\mathcal{F}$ provides constraints on the Fourier coefficients of $f \in L^2([0,1])$ in the model (\ref{1.0.2}) (see, e.g., \cite{IS_2003}, Section 3.2).   \\

We consider below the hypothesis testing setting (\ref{testing_pb0}) with $\theta_0 \neq 0$ (i.e., goodness-of-fit testing). Formally, given observations from the GSM (\ref{1.0.0}), for any {\em given} $\theta_0 \neq 0$, we will be dealing with the following goodness-of-fit testing problem
\begin{equation}
H_0: \theta=\theta_0 \quad \mathrm{versus} \quad H_1: \theta_0 \in \mathcal{E}_a,\; \theta-\theta_0\in \Theta_a(r_{\varepsilon,\sigma}),
\label{testing_pb1}
\end{equation}
where $\Theta_a(r_{\varepsilon,\sigma})$ is defined in (\ref{def:f-set}). The sequence $a$ being fixed, the main issue for the problem (\ref{testing_pb1}) is then to characterize the values of $r_{\varepsilon,\sigma} >0$ for which both hypotheses $H_0$ (called the null hypothesis) and $H_1$ (called the alternative hypothesis) are `separable'  (in a sense which will be made precise later on).\\

\begin{remark}
\label{rm:sd}
{\rm
We would like to stress that in the {\em standard}\; GSM (i.e., (\ref{1.0.0}) with $\sigma =0$),  signal detection (i.e., $\theta_0=0$) and goodness-of-fit testing (i.e., $\theta_0 \neq 0$) problems are equivalent as soon as the involved operator is injective. Indeed, without loss of generality, we can still replace the observed sequence $(Y_j)_{j\in \mathbb{N}}$ by $(\tilde Y_j)_{j\in \mathbb{N}}:=(Y_j -b_j\theta_{0,j})_{j\in \mathbb{N}}$. This is no more the case in the GSM (\ref{1.0.0}) since the sequence $(b_j)_{j\in\mathbb{N}}$ is unknown. Signal detection and goodness-of-fit problems should therefore be treated in a different manner. In this work, we only address the goodness-of-fit testing problem (\ref{testing_pb1}).  \\
}
\end{remark}

In the following, a (non-randomized) test $\Psi:=\Psi(Y,X)$ will be defined as a measurable function of the observation $(Y,X)=(Y_j,X_j)_{j\in\mathbb{N}}$ from GSM (\ref{1.0.0}) having values in the set $\lbrace 0,1 \rbrace$. By convention, $H_0$ is rejected if $\Psi=1$ and  $H_0$ is not rejected if $\Psi=0$. Then, given a test $\Psi$, we can investigate 
\begin{itemize}
\item the first kind error probability defined as 
\begin{equation}\Alpha_{\varepsilon,\sigma}(\Psi):= \mathbb{P}_{\theta_0,b}( \Psi =1),
\label{eq:type1-F}
\end{equation} 
which measures the probability to reject $H_0$ when $H_0$ is true (i.e., $\theta=\theta_0$); it is often constrained as being bounded by a prescribed level $\alpha \in ]0,1[$, and
\item the maximal second kind error probability defined as 
\begin{equation}
\Beta_{\varepsilon,\sigma}(\Theta_a(r_{\varepsilon,\sigma}),\Psi) := \sup_{\theta_0 \in \mathcal{E}_a,\,\theta-\theta_0\in \Theta_a(r_{\varepsilon,\sigma})} \mathbb{P}_{\theta,b}(\Psi=0),
\label{eq:type2}
\end{equation}
which measures the worst possible probability not to reject $H_0$ when $H_0$ is not true (i.e., when $\theta_0 \in \mathcal{E}_a$ and $\theta-\theta_0\in \Theta_a(r_{\varepsilon,\sigma})$); one would like to ensure that it is bounded by a prescribed level $\beta \in ]0,1[$.
\end{itemize}
\vspace{0.2cm}

For simplicity in our exposition, we will restrict ourselves to $\alpha$-level tests, 
i.e., tests $\Psi_\alpha$ satisfying $\Alpha_{\varepsilon,\sigma}(\Psi_\alpha)\leq \alpha$,
for any fixed value $\alpha \in ]0,1[$.\\

Let $\alpha,\beta\in ]0,1[$ be given, and let $\Psi_\alpha$ be an $\alpha$-level test. 

\begin{definition}
The separation radius of the $\alpha$-level test $\Psi_\alpha$ over the class $\mathcal{E}_a$ is defined as
$$ r_{\varepsilon,\sigma}(\mathcal{E}_a,\Psi_\alpha,\beta) := \inf \left\lbrace r_{\varepsilon,\sigma}>0: \ \Beta_{\varepsilon,\sigma}(\Theta_a(r_{\varepsilon,\sigma}),\Psi_\alpha)  \leq \beta\right\rbrace,$$
where the maximal second kind error probability $\Beta_{\varepsilon,\sigma}(\Theta_a(r_{\varepsilon,\sigma}),\Psi_\alpha)$ is defined in (\ref{eq:type2}). 
\label{eq:minimax_separation_radius}
\end{definition}

In some sense, the separation radius $r_{\varepsilon,\sigma}(\mathcal{E}_a,\Psi_\alpha,\beta)$ corresponds to the smallest possible value of the available signal $\| \theta - \theta_0 \|$ for which $H_0$ and $H_1$ can be `separated' by the $\alpha$-level test $\Psi_\alpha$ with maximal second kind error probability, bounded by a prescribed level $\beta \in ]0,1[$.

\begin{definition}
\label{def:minimax_radius}
The minimax separation radius $\tilde{r}_{\varepsilon,\sigma}:=\tilde{r}_{{\varepsilon,\sigma}}(\mathcal{E}_a, \alpha, \beta)>0$ over the class $\mathcal{E}_a$ is defined as
\begin{equation}
\tilde r_{{\varepsilon,\sigma}}:= \inf_{\tilde \Psi_\alpha: \,\Alpha_{\varepsilon,\sigma}(\tilde \Psi_\alpha)\leq \alpha} r_{\varepsilon,\sigma}(\mathcal{E}_a, \tilde\Psi_\alpha,\beta).
\label{eq:minimax_radius}
\end{equation}
\end{definition}

The minimax separation radius $\tilde r_{{\varepsilon,\sigma}}$ corresponds to the smallest radius $r_{{\varepsilon,\sigma}} >0$ such that there exists some $\alpha$-level test $\tilde \Psi_{\alpha}$ for which the maximal second kind error probability 
$\Beta_{\varepsilon,\sigma}(\Theta_a(r_{\varepsilon,\sigma}),\tilde \Psi_\alpha)$ is not greater than $\beta$.

\subsection{Summary of the Results}
\label{s:summary}

Our aim is to establish `optimal' separation conditions for the goodness-of-fit testing problem (\ref{testing_pb1}). This task requires, in particular, precise (non-asymptotic) controls of the first kind error probability $\Alpha_{\varepsilon,\sigma}(\Psi_\alpha)$ and the maximal second kind error probability $\Beta_{\varepsilon,\sigma}(\Theta_a(r_{\varepsilon,\sigma}),\Psi_\alpha)$ (of a specific test $\Psi_\alpha$ that will be made precise in Section \ref{sec:sc-o-t}) by prescribed levels $\alpha, \beta \in ]0,1[$, respectively. Such controls allow us to derive both upper and lower bounds on the minimax separation radius $\tilde r_{\varepsilon,\sigma}$, as summarized in the following theorem.\\

\begin{thm}
\label{thm:sum}
Let $\alpha,\beta \in ]0,1[$ be fixed, such that $\alpha \leq \beta$.  Consider the goodness-of-fit testing problem (\ref{testing_pb1}). Then, there exist explicit positive constants$^1$ $\tilde C(\alpha,\beta)>0$, $C_{\alpha,\beta}>0$, $c_{\alpha,\beta}>0$  and $\sigma_0 \in ]0,1[$ such that, for all\; $0<\sigma \leq \sigma_0$ and for each $\varepsilon >0$,
$$
(i) \quad \tilde r_{\varepsilon,\sigma}^2 \leq \inf_{D \in \mathbb{N}}\left[\tilde C(\alpha,\beta) 
\varepsilon^2 \sqrt{\sum_{j=1}^{D \wedge M_1}  b_j^{-4}} + \big(7+4\sqrt{\ln(2/\alpha)}\big) \left[\sigma^2 \ln^{3/2}(1/\sigma) \vee a^{-2}_{D \wedge M_0}\right] \right],
$$
and, for all $\varepsilon, \sigma>0$,
$$
(ii) \quad \tilde r^2_{\varepsilon,\sigma} \geq  \left\{ \frac{C^2_{\alpha,\beta}}{16}\, \sigma^2\, 
\max_{1 \leq D \leq M_2} [ b_D^{-2} a_D^{-2} ] \right\} \vee \left\{ \sup_{D \in \mathbb{N}} \left[c_{\alpha,\beta}\,
\varepsilon^2 \sqrt{\sum_{j=1}^D b_j^{-4}} \wedge a_D^{-2} \right] \right\},
$$
where the bandwidths $M_0, M_1$ and $M_2$ depend\footnote{For the sake of brevity, these quantities are made precise in the subsequent sections} on both $(b_j)_{j\in \mathbb{N}}$ and $\sigma$. 
\end{thm}

Theorem \ref{thm:sum} provides a precise description on the behavior of the minimax separation radius $\tilde r_{\epsilon,\sigma}$ in terms of the sequences $(a_j)_{j\in \mathbb{N}}$ and $(b_j)_{j\in \mathbb{N}}$ and of the noise levels $\epsilon$ and $\sigma$. It is worth pointing out that this control is \textit{non-asymptotic}. There is indeed a technical constraint on the value of $\sigma$ ($0< \sigma \leq \sigma_0$, $\sigma_0 \in ]0,1[$), but we do not assume its convergence towards $0$, i.e., we work with {\em fixed} values of the noise levels $\varepsilon$ and $\sigma$. \\

Then, we apply the above result on specific problems. Namely, we consider various behaviors for both sequences $(a_j)_{j\in \mathbb{N}}$ and $(b_j)_{j\in \mathbb{N}}$,  and discuss the properties of the associated minimax separation radii $\tilde r_{\varepsilon,\sigma}$. Concerning the eigenvalues $(b_j^2)_{j\in \mathbb{N}}$ of the operator $A^*A$, we will alternatively consider situations where 
$$ b_j   \sim j^{-t} \quad \mathrm{or} \quad b_j   \sim \exp(-jt), \quad \forall j\in \mathbb{N}, \;\; \text{for some} \;\; t >0.$$
The first case corresponds to the so-called \textit{mildly ill-posed} problems while the second one corresponds to \textit{severely ill-posed} problems.  Concerning the ellipsoids $\mathcal{E}_a$, i.e., the sequence $(a_j)_{j\in \mathbb{N}}$, two different kinds of smoothness will be investigated, namely, 
$$ a_j   \sim j^{s} \quad \text{or} \quad a_j   \sim \exp(js), \quad \forall j\in \mathbb{N},  \;\; \text{for some} \;\; s>0,$$
the so-called {\em ordinary-smooth}  and {\em super-smooth} cases, respectively.
In the above scenarios, we apply Theorem \ref{thm:sum} and describe the associated upper and lower bounds on the minimax separation radius $\tilde r_{\varepsilon,\sigma}$. They are, respectively, displayed in Table \ref{tab:sep-u} and Table \ref{tab:sep-l}.\\
\\

\begin{table}[h]
\centering
\begin{tabular}{|c|c|c|}
\hline
{\em Goodness-of-Fit }&{\em ordinary-smooth} &{\em super-smooth} \\
{\em Testing Problem} & $a_j \sim j^{s}$ & $a_j \sim \exp\{js\}$ \\
\hline {\em mildly ill-posed} &$\varepsilon^{4s/(2s+2t+1/2)} \vee [\sigma \ln^{3/4}(1/\sigma)]^{2[(s/t)\wedge 1]}$&$
\varepsilon^2(\ln (1/\varepsilon))^{2t+1/2} \vee \sigma^2 \ln^{3/2}(1/\sigma)$\\
 $b_j\sim j^{-t}$ & & \\
\hline
{\em severely ill-posed} & $(\ln (1/\varepsilon))^{-2s} \vee [\ln (1/\sigma \ln^{-1/2}(1/\sigma))]^{-2s}$&$
\varepsilon^{2s/(s+t)} \vee [\sigma \ln^{1/2}(1/\sigma)]^{2[(s/t)\wedge 1]}$\\
$b_j \sim \exp\{-jt\}$ & & \\
\hline
\end{tabular}
\caption{\textrm{Minimax goodness-of-fit testing with {\em unknown} singular values: {\em upper bounds} on the minimax separation radius $\tilde r^2_{\varepsilon,\sigma}$ for $0 <\varepsilon \leq \varepsilon_0$, $\varepsilon_0 \in ]0,1[$, and $0 <\sigma \leq \sigma_0$, $\sigma_0 \in ]0,1[$,  for all $t,s>0$}.}
\label{tab:sep-u}
\end{table}

\begin{table}[h]
\centering
\begin{tabular}{|c|c|c|}
\hline
{\em Goodness-of-Fit }&{\em ordinary smooth} &{\em super smooth} \\
{\em Testing Problem} & $a_j \sim j^{s}$ & $a_j \sim \exp\{js\}$ \\
\hline {\em mildly ill-posed} &$\varepsilon^{4s/(2s+2t+1/2)} \vee \sigma^{2[(s/t)\wedge 1]}$&$
\varepsilon^2(\ln \varepsilon^{-1})^{2t+1/2} \vee \sigma^2 $\\
 $b_j\sim j^{-t}$ & & \\
\hline
{\em severely ill-posed} & $(\ln (1/\varepsilon))^{-2s} \vee (\ln (1/\sigma))^{-2s}$&$
\varepsilon^{2s/(s+t)} \vee \sigma^{2[(s/t)\wedge 1]}$\\
$b_j \sim \exp\{-jt\}$ & & \\
\hline
\end{tabular}
\caption{\textrm{Minimax goodness-of-fit testing with {\em unknown} singular values: {\em lower bounds} on the minimax separation radius $\tilde r^2_{\varepsilon,\sigma}$ for $0 <\varepsilon \leq \varepsilon_0$, $\varepsilon_0 \in ]0,1[$, and $0 <\sigma \leq \sigma_0$, $\sigma_0 \in ]0,1[$, for all $t ,s>0$}.}
\label{tab:sep-l}
\end{table}

Looking at these tables, both lower and upper bounds coincide in every considered case, up to a logarithm term that depends on the noise level $\sigma$. Hence, Theorem \ref{thm:sum} provides a sharp control on the minimax separation radius $\tilde r_{\varepsilon,\sigma}$ in various settings.  The interesting property of such minimax separation radii is that they have the same structure whatever the considered situation: a maximum between two terms depending, respectively, on the noise levels $\epsilon$ and $\sigma$. It is also worth pointing out that the first term depending on $\epsilon$ corresponds to the minimax separation radius in the case where the operator is known (i.e., $\sigma=0$), as displayed in Table \ref{tab:sep}.   \vspace{0.5cm}

\begin{table}[h]
\centering
\begin{tabular}{|c|c|c|}
\hline
{\em Goodness-of-Fit }&ordinary-smooth &super-smooth \\
{\em Testing Problem} & $a_j \sim j^{s}$ & $a_j \sim \exp\{js\}$ \\
\hline 
mildly ill-posed &$\varepsilon^{4s/(2s+2t+1/2)}$&$
\varepsilon^2(\ln \varepsilon^{-1})^{2t+1/2}$\\
 $b_j\sim j^{-t}$ & & \\
\hline
severely ill-posed & $(\ln \varepsilon^{-1})^{-2s}$&$
\varepsilon^{2s/(s+t)}$\\ 
$b_j \sim \exp\{-jt\}$ & & \\
\hline
\end{tabular}
\caption{\textit{Minimax goodness-of-fit testing with {\em known} singular values: the separation rates $\tilde r^2_{\varepsilon}$ for $0 <\varepsilon \leq \varepsilon_0$, $\varepsilon_0 \in ]0,1[$,  for all $t,s>0$}.}
\label{tab:sep}
\end{table}

The results displayed in Theorem \ref{thm:sum} and Tables \ref{tab:sep-u}, \ref{tab:sep-l} can also be understood as follows.  Two problems are at hand: detection of the underlying signal (with a minimax separation radius that only depends on $\epsilon$) and detection of the `frequencies' $j$ for which the terms $b_j$ can be replaced by observations $X_j$ without loss of precision (with a minimax separation radius that depends only on $\sigma$). The final minimax separation radius is then the maximum of these two terms, i.e., the signal detection hardness is related to the most difficult underlying problem.  We stress that such phenomenon has already been discussed in the minimax estimation framework, see e.g., \cite{DHPV_2012}, \cite{JS_2013}.

%

\section{Upper Bound on the Minimax Separation Radius}
\label{sec:sc-o-t} 

In this section, we first propose an $\alpha$-level testing procedure. Then, we investigate its maximal second kind error probability and establish a non-asymptotic upper bound on the minimax separation radius (which corresponds to item $(i)$ of Theorem \ref{thm:sum}). Finally, in Section \ref{subsec:AUB}, we provide a control of the upper bounds for minimax separation radii for the specific cases displayed in Table \ref{tab:sep-u}.

\subsection{The Spectral Cut-Off Test} 
For a given $\theta_0 \neq 0$, the aim of the goodness-of-fit testing problem (\ref{testing_pb1}) is to determine whether or not $\theta=\theta_0$. In particular, for any given $j \in \mathbb{N}$, one would like to infer the corresponding value $\theta_j$ from the observation $(Y,X)=(Y_j,X_j)_{j\in\mathbb{N}}$ from GSM (\ref{1.0.0}). Typically, for any given $j \in \mathbb{N}$, one may use the `naive' estimate ${\hat \theta}_j$ of $\theta_j$, defined by
$$
{\hat \theta}_j:=\frac{Y_j}{X_j}=\frac{b_j}{X_j} \theta_j + \varepsilon \frac{1}{X_j} \xi_j, \quad j \in \mathbb{N}.
$$
In order to ensure a `good' approximation of $\theta_j$ by ${\hat \theta}_j$ (in a sense which will be made precise later on), a precise control of the ratio $b_j/X_j$ is required. 
To this end, we want to avoid coefficients for which $X_j \lesssim \sigma$, namely for which the observation $X_j$ is of the order of the corresponding noise level $\sigma$, that does not have `discriminatory' power . Therefore, we will restrict ourselves to coefficients $X_j$ with indices $1 \leq j \leq M$, where the bandwidth $M$ is defined by
\begin{equation}
M := \inf\{j \in \mathbb{N}:\, |X_j| \leq \sigma h_j\}-1,
\label{eqn:M}
\end{equation}
where, for all $j\in\mathbb{N}$,
\begin{equation}
h_j = 16 \sqrt{ \ln \left( \frac{\kappa j^2}{\alpha} \right)} + \sqrt{2 \ln \left( \frac{10}{\alpha} \right)},
\label{eq:h}
\end{equation}
for some $\kappa>\exp(1)$.

\begin{remark}{\rm
The value of $\kappa$ is, in some sense, related to the value of the first kind error probability of the suggested testing procedure. We will see below that the value $\kappa = 5(3\pi^2 +12)/6$ is convenient to our purpose. We stress that $\kappa$ is not \textit{a regularization parameter}:  an `optimal' value of $\kappa$ only allows to get `optimal' constants in the final results but will not change the order of the corresponding minimax separation rates. Finding optimal constants is outside the scope of this work.}
 \end{remark}

\begin{figure}
\begin{center}
\begin{tikzpicture}[scale=1.5]
\draw [thick,->] (-0.5,0) -- (5,0);
\draw (5.25,0) node {$j$};
\draw [thick, ->] (0,-0.5) -- (0,4);
\draw (-0.15,-0.15) node {$0$};
\draw (-0.15,4.2) node {$b_j$};
\draw [black,domain=0.3:5,thin] plot (\x, {(\x)^(-1)});
\draw [thin,dotted] (0.8,0) -- (0.8,1.25);
\draw (0.9,-0.2) node {$M_0$};
\draw[fill] (0.8,0) circle (0.03cm);
\draw [thin,dotted] (0,1.25) -- (0.8,1.25);
\draw[fill] (0,1.25) circle (0.03cm);
\draw (-0.9,1.25) node {$\sigma h_{0,M_0}$};
\draw [thin,dotted] (2.3,0) -- (2.3,0.43); 
\draw (2.4,-0.2) node {$M_1$};
\draw[fill] (2.3,0) circle (0.03cm);
\draw [thin,dotted] (0,0.43) -- (2.3,0.43);
\draw[fill] (0,0.43) circle (0.03cm);
\draw (-0.8,0.43) node {$\sigma h_{1,M_1}$};
%
\draw [blue,domain=0.3:5,thin] plot [samples=100] (\x, {(\x)^(-1)+ (rnd-0.5)/1.5});
\draw [black,domain=0:5,thin,dashed] plot (\x, {0.7+(ln(1+\x/3))});
\end{tikzpicture}
\end{center}
\caption{\textit{An illustration of the spatial positions of the bandwidths $M_0$ and $M_1$, defined in (\ref{eqn:M0-M1-M}). The decreasing solid curve corresponds to the values of the sequence $b=(b_j)_{j\in \mathbb{N}}$ with respect to the index $j \in \mathbb{N}$, while the oscillating curve demonstrates one realization of the random sequence $X=(X_j)_{j\in \mathbb{N}}$ according to the GSM (\ref{1.0.0}). 
The increasing dashed curve draws the behavior of the sequence $\sigma h_j$.
For the corresponding random `bandwidth' $M$ defined in  (\ref{eqn:M0-M1-M}), Lemma \ref{lem:tech_lem_1} shows that  $M \in [M_0,M_1[$ with high probability.}}
\label{fig:Fig-1}
\end{figure}
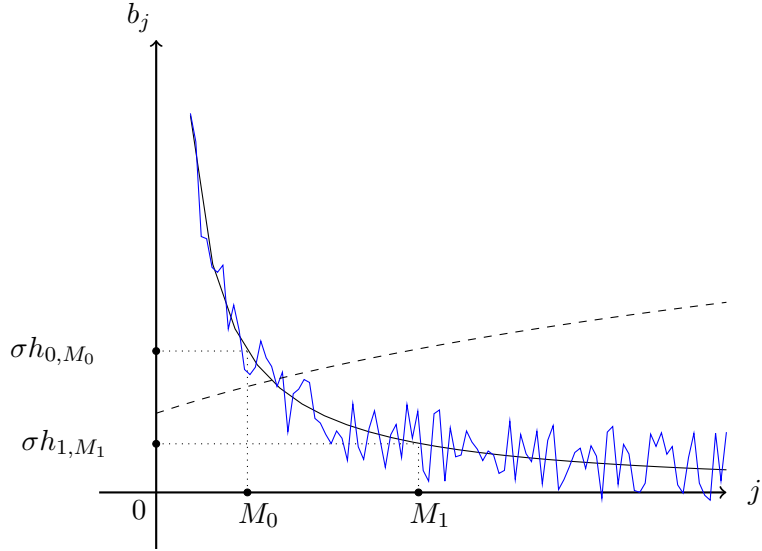

 The bandwidth $M$ is a random variable but can be controlled in the sense that
$M \in [M_0,M_1[$ with high probability (see Lemma \ref{lem:tech_lem_1} for precise computations and Figure \ref{fig:Fig-1} for a graphical illustration), where the bandwidths $M_0$ and $M_1$ are defined by
\begin{equation}
\label{eqn:M0-M1-M}
\begin{cases}
M_0 :=\inf\{j \in \mathbb{N}:\, b_j \leq \sigma h_{0,j}\}-1, &\\
M_1 :=\inf\{j \in \mathbb{N}:\, b_j \leq \sigma h_{1,j}\}, &
\end{cases}
\end{equation}
and the sequences $h_0=(h_{0,j})_{j \in \mathbb{N}}$, $h_1=(h_{1,j})_{j \in \mathbb{N}}$ satisfy 
\begin{eqnarray}
h_{0,j} &=& 18 \sqrt{ \ln \left( \frac{\kappa j^2}{\alpha} \right)} + \sqrt{2 \ln \left( \frac{10}{\alpha} \right)}, \label{eq:h_0} \\
h_{1,j} &=& 16 \sqrt{ \ln \left( \frac{\kappa j^2}{\alpha} \right)}, \label{eq:h_1} 
\end{eqnarray}
for all $j\in\mathbb{N}$. The sequences $h=(h_j)_{j \in \mathbb{N}}$, $h_0=(h_{0,j})_{j \in \mathbb{N}}$ and $h_1=(h_{1,j})_{j \in \mathbb{N}}$ in the definition of $M_0$, $M_1$ and $M$ allow a `uniform' control of the standard Gaussian sequence $\eta=(\eta_j)_{j \in \mathbb{N}}$ (associated with $X=(X_j)_{j \in \mathbb{N}}$), for all  $1 \leq j \leq M_1$ (see Lemmas \ref{lem:tech_lem_1}, \ref{lem:tech_lem_2} and \ref{lem:tech_lem_3} in Section \ref{s:appendix}).\\

We are now in the position to construct a (spectral cut-off) testing procedure. According to the methodology proposed earlier in the literature (see e.g. \cite{Baraud}, \cite{IS_2003} or \cite{MS_2014}),  our test will be based on an estimation of $\| \theta-\theta_0\|^2$. For any fixed $D \in \mathbb{N}$, consider the test statistic
\begin{equation}
\label{eq:sco-tests-tatist}
T_{D,M}:=\sum_{j=1}^{D \wedge M}\left( \frac{Y_j}{X_j}-\theta_{j,0}\right)^2.
\end{equation}
Given a prescribed level $\alpha\in ]0,1[$ for the kind error probability, the associated spectral cut-off test is then defined as
\begin{equation}
\label{eq:sco-test}
\Psi_{D,M}:=\mathbf{1}\{T_{D,M} > t_{1-\alpha,D}(X)\},
\end{equation}
where 
\begin{equation}
\label{eqn:comp-a-quantile}
t_{1-\alpha,D}(X):=\varepsilon^2 \sum_{j=1}^{D \wedge M}X_j^{-2} + C(\alpha) 
\varepsilon^2 \sqrt{\sum_{j=1}^{D \wedge M}X_j^{-4}} +(1+\sqrt{x_{\alpha/2}}) \left[\sigma^2 \ln^{3/2}(1/\sigma) \vee a_{D \wedge M}^{-2}\right],
\end{equation}
and 
\begin{equation}
C(\alpha)=3 \sqrt{x_{\alpha/2}} + 2 x_{\alpha/2}>0, \quad  x_\gamma := \ln(1/\gamma) \ \forall \gamma \in ]0,1[.
\label{eq:ca}
\end{equation}
In other words, if the 'estimator' $T_{D,M}$ of $\| \theta-\theta_0\|^2$ is greater than the fixed threshold $t_{1-\alpha,D}(X)$, $\theta$ and $\theta_0$ are very unlikely to be close to each other, and we will reject $H_0$. \\

\begin{remark}
{\rm Under $H_0$, $Y_j=b_j\theta_{j,0}+\varepsilon \xi_j$, $j \in \mathbb{N}$, and, hence,
$$
T_{D,M}=\sum_{j=1}^{D \wedge M}\left[\left(\frac{b_j}{X_j} -1\right)\theta_{j,0}+\varepsilon X_j^{-1}\xi_j\right]^2.
$$
Therefore, the law of $T_{D,M}$ is not available and, thus, its corresponding $(1-\alpha)$-quantile is not computable in practice, since the sequence $b=(b_j)_{j \in \mathbb{N}}$ is unknown. However, Proposition \ref{prop:size-sco-test} below ensures that the threshold $t_{1-\alpha,D}(X)$ defined in (\ref{eqn:comp-a-quantile}) provides a computable upper bound on this quantile. 
}
\end{remark}

First, we focus on the first kind error probability. The following proposition states that the spectral cut-off test $\Psi_{D,M}$ defined in (\ref{eq:sco-test})-(\ref{eqn:comp-a-quantile}), is an $\alpha$-level test.

\begin{proposition}
\label{prop:size-sco-test}
Let $\alpha \in ]0,1[$ be fixed. Consider the goodness-of-fit testing problem (\ref{testing_pb1}). Then, setting $\kappa =5(3\pi^2+12)/6$, there exists $\sigma_0 \in ]0,1[$ such that, for all \,$0<\sigma \leq \sigma_0$ and for each $\varepsilon >0$, the spectral cut-off test $\Psi_{D,M}$, defined in (\ref{eq:sco-test})-(\ref{eqn:comp-a-quantile}), is an $\alpha$-level test,  i.e.,
\begin{equation}
\label{eq:sco-test-size}
\Alpha_{\varepsilon,\sigma}(\Psi_{D,M}) \leq \alpha.
\end{equation}
\end{proposition}

\noindent 
The proof is postponed to Section \ref{s:proof1}.

\begin{remark}{\rm In order to shed light on the term $\sigma_0$,  we provide bellow a heuristic argument. Note that, under $H_0$, thanks to a (rough) Taylor expansion,
$$
T_{D,M}\simeq \sum_{j=1}^{D \wedge M}\left[\varepsilon b_j^{-1}\xi_j + \sigma b_j^{-1}\theta_{j,0} \eta_j \right]^2.
$$
Compared to the `noise-free' case (i.e., $\sigma =0$), we have in some sense to deal with the additional term $\sigma b_j^{-1} \theta_{j,0}\eta_j$. Two scenarios are at hand
\begin{itemize}
\item If $\sup_j b_j^{-1} a_j^{-1} \leq C_0$, the expected amount of additional signal is
$$ \sigma^2 \sum_{j=1}^{D\wedge M} b_j^{-2}\theta_j^2 \leq \sigma^2 C_0 \|\theta\|^2,$$ 
which is of the order of the classical parametric rate $\sigma^2$. However, since $C_0$ is unknown, we use a rough standard deviation control on this additional term, which requires a logarithmic term (i.e., $\ln^{3/2}(1/\sigma)$) in the right hand side of (\ref{eqn:comp-a-quantile}). We stress that this logarithmic term can be removed if the knowledge of $C_0$ is assumed.
\item On the other hand, we can prove that $\sigma b_j^{-1} \eta_j$ (see Lemma \ref{lem:tech_lem_4}) is bounded with controlled probability, according to the construction of the bandwidth $M$ given in (\ref{eqn:M}). In such case, the additional term can be controlled by the `bias' $a^{-2}_{D\wedge M}$. \end{itemize}
Due to the additional logarithmic term mentioned above, the first kind error probability can be controlled as soon as $\sigma$ is small enough (i.e., $0<\sigma \leq \sigma_0$ for some $\sigma_0\in ]0,1[$). Unsurprisingly, it is impossible to retrieve any kind of information on the observations if the noise level $\sigma$ is too large.
}

\end{remark}

\subsection{A Non-Asymptotic Upper Bound}

We now turn our attention to the the maximal second error probability. The following proposition provides, for each noise level $\varepsilon >0$ and for noise level $\sigma$ small enough, an upper bound for the separation radius $r_{\varepsilon,\sigma}(\mathcal{E}_a,\Psi_{D,M},\beta)$ of the spectral cut-off test $\Psi_{D,M}$ defined in (\ref{eq:sco-tests-tatist})-(\ref{eqn:comp-a-quantile}). \\

\begin{proposition}
\label{thm:power-sco-test}
Let $\alpha,\beta \in ]0,1[$ be fixed, such that $\alpha \leq \beta$. Consider the goodness-of-fit testing problem (\ref{testing_pb1}). Let $\Psi_{D,M}$ be the spectral cut-off test, defined in (\ref{eq:sco-test})-(\ref{eqn:comp-a-quantile}). Then, there exists $\sigma_0 \in ]0,1[$ such that, for all\; $0<\sigma \leq \sigma_0$ and for each $\varepsilon >0$,
\begin{equation} 
\label{eq:re-ub-1}
r_{\varepsilon,\sigma}^2(\mathcal{E}_a,\Psi_{D,M},\beta) \leq \tilde{C}(\alpha,\beta) 
\varepsilon^2 \sqrt{\sum_{j=1}^{D \wedge M_1}  b_j^{-4}} + (7+4\sqrt{x_{\alpha/2}})\left[\sigma^2 \ln^{3/2}(1/\sigma) \vee a^{-2}_{D \wedge M_0}\right],
\end{equation}
where 
\begin{equation}
\tilde C(\alpha,\beta) = 16 (C(\alpha) + 3 \sqrt{x_{\beta/2}}).
\label{eq:constant_cyprus}
\end{equation}
\end{proposition}
The proof is postponed to Section \ref{s:proof2}.

\begin{remark}{\rm
According to Proposition \ref{thm:power-sco-test}, given a radius $r_{\varepsilon,\sigma} >0$, then
$$
r^2_{\varepsilon,\sigma} \geq \tilde{C}(\alpha,\beta) 
\varepsilon^2 \sqrt{\sum_{j=1}^{D \wedge M_1}  b_j^{-4}} + (7+4\sqrt{x_{\alpha/2}}) \left[\sigma^2 \ln^{3/2}(1/\sigma) \vee a^{-2}_{D \wedge M_0}\right] \; \Rightarrow \;
\Beta_{\varepsilon,\sigma}(\Theta_a(r_{\varepsilon,\sigma}),\Psi_{D,M}) \leq \beta, 
$$
and, hence,
$$
\tilde{r}^2_{\varepsilon,\sigma} \leq \inf_{D \in \mathbb{N}} \left[ \tilde{C}(\alpha,\beta) \varepsilon^2 \sqrt{\sum_{j=1}^{D \wedge M_1}  b_j^{-4}} + (7+4\sqrt{x_{\alpha/2}}) \left[\sigma^2 \ln^{3/2}(1/\sigma) \vee a^{-2}_{D \wedge M_0}\right] \right].
$$}
\end{remark}

\noindent

Note that the upper bound on the separation radius 
$r_{\varepsilon,\sigma}^2(\mathcal{E}_a,\Psi_{D,M},\beta)$ given in (\ref{eq:re-ub-1}) depends on two antagonistic terms, namely, $\varepsilon^2 \sqrt{\sum_{j=1}^{D \wedge M_1}  b_j^{-4}}$ and $\sigma^2 \ln^{3/2}(1/\sigma) \vee a^{-2}_{D \wedge M_0}$. Ideally, one would like to make this upper bound as small as possible, i.e., to obtain the weakest possible condition on $\|\theta-\theta_0\|$ such that, for any fixed $\beta \in ]0,1[$, 
$\Beta_{\varepsilon,\sigma}(\Theta_a(r_{\varepsilon,\sigma}),\Psi_{D,M}) \leq \beta$. Therefore, one would like to select $D:=D^\star$ such that
$$
D^\star := \arg\min_{D \in \mathbb{N}} \left\{\tilde{C}(\alpha,\beta) 
\varepsilon^2 \sqrt{\sum_{j=1}^{D \wedge M_1}  b_j^{-4}} + (7+4\sqrt{x_{\alpha/2}}) \left[\sigma^2 \ln^{3/2}(1/\sigma) \vee a^{-2}_{D \wedge M_0}\right]\right\},
$$
where $\tilde C(\alpha,\beta)$ is defined in (\ref{eq:constant_cyprus}). However, this `optimal' bandwith $D^\star$ is not available in practice since the sequence $b=(b_j)_{j \in \mathbb{N}}$ is not assumed to be known. To this end, we use instead the bandwidth $D:=D^\dagger$ defined as 
\begin{equation}
\label{eq:D-Dagger}
D^\dagger := \arg\min_{D \in \mathbb{N}} \left\{\tilde C(\alpha,\beta) 
\varepsilon^2 \sqrt{\sum_{j=1}^{D \wedge M}  X_j^{-4}} + (7+4\sqrt{x_{\alpha/2}}) \left[\sigma^2 \ln^{3/2}(1/\sigma) \vee a^{-2}_{D \wedge M}\right]\right\},
\end{equation}
The following theorem illustrates the performances of the corresponding spectral cut-off test $\Psi_{D^\dagger,M}$, defined in (\ref{eq:sco-test}), with $D:=D^\dagger$, defined in (\ref{eq:D-Dagger}).

\begin{thm}
\label{cor:power-sco-test}
Let $\alpha,\beta \in ]0,1[$ be fixed, such that $\alpha \leq \beta$. Consider the goodness-of-fit testing problem (\ref{testing_pb1}). Let $\Psi_{D^\dagger,M}$ be the spectral cut-off test, defined in (\ref{eq:sco-test}) with $D:=D^\dagger$, defined in (\ref{eq:D-Dagger}). Then, there exists $\sigma_0 \in ]0,1[$ such that, for all\; $0<\sigma \leq \sigma_0$ and for each $\varepsilon >0$,
\begin{equation}
\label{eq:alpha-level-cor}
\Alpha_{\varepsilon,\sigma}(\Psi_{D^\dagger,M}) \leq \alpha
\end{equation}
and
\begin{equation} 
\label{eq:re-ub-1}
r_{\varepsilon,\sigma}^2(\mathcal{E}_a,\Psi_{D^\dagger,M},\beta) \leq \inf_{D \in \mathbb{N}}\left[\tilde C(\alpha,\beta) 
\varepsilon^2 \sqrt{\sum_{j=1}^{D \wedge M_1}  b_j^{-4}} + (7+4\sqrt{x_{\alpha/2}}) \left[\sigma^2 \ln^{3/2}(1/\sigma) \vee a^{-2}_{D \wedge M_0}\right] \right],
\end{equation}
where the constant $\tilde C(\alpha,\beta)$ has been introduced in (\ref{eq:constant_cyprus}). 
\end{thm}
The proof of Theorem \ref{cor:power-sco-test} is postponed to Section \ref{s:proof3}.

\begin{remark} {\rm
According to Theorem \ref{cor:power-sco-test}, given a radius $r_{\varepsilon,\sigma} >0$, then
\begin{eqnarray*}
r^2_{\varepsilon,\sigma} &\geq &\inf_{D \in \mathbb{N}}\left[ \tilde C(\alpha,\beta) 
\varepsilon^2 \sqrt{\sum_{j=1}^{D \wedge M_1}  b_j^{-4}} + (7+4\sqrt{x_{\alpha/2}}) \left[\sigma^2 \ln^{3/2}(1/\sigma) \vee a^{-2}_{D \wedge M_0}\right] \right] \\
& & \hspace{1cm}\Rightarrow 
\Beta_{\varepsilon,\sigma}(\Theta_a(r_{\varepsilon,\sigma}),\Psi_{D^\dagger,M}) \leq \beta, 
\end{eqnarray*}
and, hence,
\begin{equation}
\tilde r^2_{\varepsilon,\sigma} \leq \inf_{D \in \mathbb{N}}\left[\tilde C(\alpha,\beta) 
\varepsilon^2 \sqrt{\sum_{j=1}^{D \wedge M_1}  b_j^{-4}} + (7+4\sqrt{x_{\alpha/2}})\ \left[\sigma^2 \ln^{3/2}(1/\sigma) \vee a^{-2}_{D \wedge M_0}\right] \right].
\label{eq:min-sep-rate-cor}
\end{equation}}
This upper bound corresponds to item $(i)$ of Theorem \ref{thm:sum}.
\end{remark}

\subsection{Upper Bounds: Specific Cases}
\label{subsec:AUB}

Our aim in this section is to determine an explicit value (in terms of the noise levels $\varepsilon$ and $\sigma$) for the upper bounds on the minimax separation radius $\tilde r_{\varepsilon, \sigma}$ obtained in Theorem \ref{cor:power-sco-test} above. To this end, we will consider well-known specific cases regarding the behavior of both sequences $(a_j)_{j\in \mathbb{N}}$ and $(b_j)_{j\in \mathbb{N}}$. According to the existing literature, we will essentially deal  with mildly and severely ill-posed problems with ellipsoids of ordinary-smooth and super-smooth functions (see also Section \ref{s:summary} for formal definitions).

\begin{thm}
\label{thm:asupper}
Consider the goodness-of-fit testing problem (\ref{testing_pb1}) when observations are given by (\ref{1.0.0}), and the signal of interest has smoothness governed by (\ref{eq:ellips_a}). Then, 
\begin{itemize}
\item[(i)] If $b_j \sim j^{-t}$, $t>0$, and $a_j \sim j^s$, $s>0$, for all $j \in \mathbb{N}$, then,
there exists $\varepsilon_0, \sigma_0 \in ]0,1[$ such that, for all $0 < \varepsilon \leq \varepsilon_0$ and $0 <\sigma \leq \sigma_0$, the minimax separation radius $\tilde r_{\varepsilon, \sigma}$ satisfies
\begin{equation*}
\tilde r^2_{\varepsilon, \sigma} \lesssim \varepsilon^{\frac{4s}{2s+2t+1/2}}
\vee \left[\sigma \ln^{3/4}(1/\sigma) \right]^{2 \left(\frac{s}{t} \wedge 1\right)}.
\end{equation*}

\item[(ii)] If $b_j \sim j^{-t}$, $t>0$, and $a_j \sim\exp\{js\}$, $s>0$, for all $j \in \mathbb{N}$, then,
there exists $\varepsilon_0, \sigma_0 \in ]0,1[$ such that, for all $0 < \varepsilon \leq \varepsilon_0$ and $0 <\sigma \leq \sigma_0$, the minimax separation radius $\tilde r_{\varepsilon, \sigma}$ satisfies
$$
\tilde r^2_{\varepsilon, \sigma} \lesssim \varepsilon^2 \left[\ln\left(1/\varepsilon\right)\right]^{\left(2t+\frac{1}{2}\right)} \vee \sigma^{2} \ln^{\frac{3}{2}}\left(1/\sigma\right).
$$

\item[(iii)] If $b_j \sim \exp\{-jt\}$, $t>0$, and $a_j \sim j^s$, $s>0$, for all $j \in \mathbb{N}$, then,
there exists $\varepsilon_0, \sigma_0 \in ]0,1[$ such that, for all $0 < \varepsilon \leq \varepsilon_0$ and $0 <\sigma \leq \sigma_0$, the minimax separation radius $\tilde r_{\varepsilon, \sigma}$ satisfies
$$
\tilde r^2_{\varepsilon, \sigma} \lesssim \left[\ln\left(1/\varepsilon\right)\right]^{-2s} \vee \left[\ln\left(\frac{1}{\sigma \ln^{1/2}(1/\sigma)}\right)\right]^{-2s}.
$$

\item[(iv)] If $b_j \sim \exp\{-jt\}$, $t>0$, and $a_j \sim \exp\{js\}$, $s>0$, for all $j \in \mathbb{N}$, then,
there exists $\varepsilon_0, \sigma_0 \in ]0,1[$ such that, for all $0 < \varepsilon \leq \varepsilon_0$ and $0 <\sigma \leq \sigma_0$, the minimax separation radius $\tilde r_{\varepsilon, \sigma}$ satisfies
$$
\tilde r^2_{\varepsilon, \sigma} \lesssim \varepsilon^{\frac{2s}{s+t}}
\vee [\sigma\ln^{1/2}(1/\sigma)]^{2 \left(\frac{s}{t} \wedge 1\right)}.
$$
\end{itemize}
\end{thm}

The proof is postponed to Section \ref{s:p_rates}. The main task is to compute the asymptotic trade-off between both antagonistic terms 
$\varepsilon^2 \sqrt{\sum_{j=1}^{D \wedge M_1}  b_j^{-4}} $ and $\left[\sigma^2 \ln^{3/2}(1/\sigma) \vee a^{-2}_{D \wedge M_0}\right] $ in the upper bounds on the minimax separation radius $\tilde r_{\varepsilon, \sigma}$ displayed in (\ref{eq:min-sep-rate-cor}).

\section{Lower Bounds on the Minimax Separation Radius}
\label{lastref}
We establish a non-asymptotic lower bound on the minimax separation radius (which corresponds to item $(ii)$ of Theorem \ref{thm:sum}). In order to do this, we consider two special cases of the GSM (\ref{1.0.0}), namely the situations where
\begin{itemize}
\item[(a)] $\varepsilon=0$:  the signal is observed without noise but the eigenvalues of the operator at hand are still noisy, and 
\item[(b)] $\sigma=0$: the `classical' model (see, e.g., \cite{ISS_2012} or \cite{LLM_2012}) where the eigenvalues of the operator at hand are known.  
\end{itemize}
Both models (a) and (b) correspond to some `extreme' situations but provide, in some sense, a benchmark for the problem at hand.  We first establish a lower bound for the case (a) in Section \ref{s:lb1} and recall the lower bound for the case (b) (that has already been discussed in, e.g., \cite{Baraud}, \cite{LLM_2012} or \cite{MS_2014}) in Section \ref{s:lb2}. Then, we establish in Section \ref{s:lb3} that the minimax separation radius associated to goodness-of-fit testing problem (\ref{testing_pb1}) is always greater than the maximum of the minimax separation radii associated to the cases (a) and (b). Finally, in Section \ref{subsec:ALB}, we provide a control of the lower bounds for minimax separation radii for the specific cases displayed in Table \ref{tab:sep-l}.

\subsection{Lower Bounds for a GSM with $\varepsilon=0$}
\label{s:lb1}
We consider the GSM (\ref{1.0.0}) with $b=\bar b$ and $\varepsilon=0$, i.e.,
\begin{equation}
\label{eq:lower-bound-e=0}
\begin{cases} Y_j=\bar b_j \theta_j, & j \in \mathbb{N}, \\
X_j= \bar b_j +\sigma \, \eta_j, & j \in \mathbb{N}.
\end{cases}
\end{equation}
For a given sequence $b=(b)_{j \in \mathbb{N}}$, define
$$
\mathcal{B}(b)=\{\nu \in l^2(\mathbb{N}):\; C_0 |b_j| \leq |\nu_j| \leq C_1 | b_j|, \;j \in \mathbb{N},\; 0<C_0 \leq 1 \leq C_1 <+\infty\}.
$$

Given observations from the GSM (\ref{eq:lower-bound-e=0}), for any given $\theta_0 \neq 0$ and $\bar b\in \mathcal{B}(b)$, we consider the following goodness-of-fit testing problem
\begin{equation}
H_0: \theta=\theta_0 \quad \mathrm{versus} \quad H_1: \theta_0 \in \mathcal{E}_a,\; \theta-\theta_0\in \Theta_a(r_\sigma), \bar b\in \mathcal{B}(b),
\label{testing_pb1-e=0}
\end{equation}
where $\Theta_a(r_\sigma)=\left\lbrace \mu \in \mathcal{E}_a, \ \|\mu \| \geq r_\sigma \right\rbrace$.\\

Our aim below is to provide a lower bound 
on the minimax separation radius $\tilde r_{0,\sigma}$, defined as
$$
\tilde r_{0,\sigma}:= \inf_{\tilde \Psi_\alpha: \,\Alpha_{0,\sigma}(\tilde \Psi_\alpha)\leq \alpha} r_{0,\sigma}(\mathcal{E}_a, \tilde\Psi_\alpha,\beta),
$$
where $r_{0,\sigma}(\mathcal{E}_a,\Psi_\alpha,\beta) $ is the separation radius of any given $\alpha$-level test $\Psi_\alpha$, defined as
$$ 
r_{0,\sigma}(\mathcal{E}_a,\Psi_\alpha,\beta) := \inf \left\lbrace r_{\sigma}>0: \ \Beta_{0,\sigma,b}(\Theta_a(r_{\sigma}),\mathcal{B}(b),\Psi_\alpha)  \leq \beta\right\rbrace,$$
and  $\beta_{0,\sigma,b}(\Theta(r_\sigma),\mathcal{B}(b),\Psi_\alpha)$ is the associated maximal second kind error probability, defined as
$$
\beta_{0,\sigma,b}(\Theta(r_\sigma),\mathcal{B}(b),\Psi_\alpha):=
\sup_{\substack{\theta_0 \in \mathcal{E}_a,\; \theta-\theta_0\in \Theta_a(r_\sigma) \\ \bar b\in \mathcal{B}(b)}} \mathbb{P}_{\theta, \bar b}(\Psi_\alpha=0).
$$  

The following proposition states a lower bound for the minimax separation radius $\tilde r_{0,\sigma}$ of the goodness-of-fit testing problem (\ref{testing_pb1-e=0}).\\

\begin{proposition}
\label{prop:1-e=0}
Assume that $(Y,X)=(Y_j,X_j)_{j\in\mathbb{N}}$ are observations from the GSM \ref{eq:lower-bound-e=0}) and consider the goodness-of-fit testing problem (\ref{testing_pb1-e=0}). Let $\alpha \in ]0,1[$ and $\beta \in\, ]0,1-\alpha[$ be given. Then, for every $\sigma >0$, the minimax separation radius $\tilde r_{0,\sigma}$ is lower bounded by 
\begin{equation}
\label{eq:prop-e=0}
\tilde r_{0,\sigma} \geq \frac{C_{\alpha,\beta}}{4}\, \sigma\, \max_{1 \leq D \leq M_2} \left[b_{D}^{-1} a_{D}^{-1}\right],
\end{equation}
where
\begin{equation}
\label{eq:M2-def}
M_2:=\sup \left\{ D \in \mathbb{N}:\; C_{\alpha,\beta}\, \sigma b_D^{-1} \leq 2 \quad \text{and} \quad G_D(C_0,C1) \geq \frac{1}{\sqrt{1+4(1-\alpha-\beta)^2}}\right\}
\end{equation}
with
\begin{equation}
\label{eq:C-alpha-beta}
C_{\alpha,\beta}=\ln (1+4(1-\alpha-\beta)^2)>0 \;\; \text{and} \;\;
G_D(C_0,C_1)=
\frac{1}{\sigma \sqrt{2\pi}} \int_{C_0b_{D}}^{C_1b_{D}}\exp\left\{-\frac{1}{2\sigma^2}(t-b_{D})^2\right\}dt,
\end{equation}
for some constants $0<C_0 \leq 1 \leq C_1 <+\infty$.
\end{proposition}
The proof is postponed to Section \ref{s:prlowerb}.

\begin{remark}
{\rm } Note that 
$$
G_D(C_0,C1) = \Phi \left((C_1-1)\frac{b_D}{\sigma} \right)-\Phi \left((C_0-1)\frac{b_D}{\sigma} \right),
$$
where $\Phi(\cdot)$ is the cumulative distribution function of the standard Gaussian distribution. Hence, 
$$
G_D(C_0,C1) \geq \frac{1}{\sqrt{1+4(1-\alpha-\beta)^2}} \Leftrightarrow b_D \geq \sigma K,
$$
where $K:=K(C_0,C_1,\alpha,\beta)>0$. Then $M_2$ in (\ref{eq:M2-def}) can be re-expressed as
\begin{equation}
\label{eq:new-M2}
M_2:= \sup \left\{ D \in \mathbb{N}: b_D \geq \sigma [K \vee C_{\alpha,\beta}/2]\right\}.
\end{equation} 
This expression $M_2$ in (\ref{eq:new-M2}) can be compared to the respective expressions of $M_0$ and $M_1$ defined in (\ref{eqn:M0-M1-M}). In particular, we point-out that there is no logarithmic term involved in $M_2$. 
\end{remark}

\subsection{Lower Bounds for the GSM when $\sigma=0$}
\label{s:lb2}
We consider the GSM (\ref{1.0.0}) with $\sigma=0$, i.e.,
\begin{equation}
\label{eq:lower-bound-s=0}
\begin{cases} Y_j=b_j \theta_j+ \varepsilon \xi_j, & j \in \mathbb{N}, \\
X_j=  b_j, & j \in \mathbb{N}.
\end{cases}
\end{equation}
Note that, in this case, the above model can be re-expressed as
\begin{equation}
\label{eq:model-s=0}
Y_j=b_j\theta_j+\varepsilon \xi_j, \quad j \in \mathbb{N},
\end{equation}
where $b=(b_j)_{j \in \mathbb{N}}$ is a {\em known} positive sequence.\\

The following proposition states a lower bound for the minimax separation radius $\tilde r_{\varepsilon,0}$, defined in (\ref{def:minimax_radius}) with $\sigma=0$, of the following goodness-of-fit testing problem 
\begin{equation}
H_0: \theta=\theta_0 \quad \mathrm{versus} \quad H_1: \theta_0 \in \mathcal{E}_a,\; \theta-\theta_0\in \Theta_a(r_{\varepsilon,0}),
\label{testing_pb-2}
\end{equation}
where $\Theta_a(r_{\varepsilon,0})$ is defined in (\ref{def:f-set}) with $\sigma=0$.\\

\begin{proposition}
\label{prop:1-s=0}
Assume that $Y=(Y_j)_{j\in\mathbb{N}}$ are observations from the GSM \ref{eq:model-s=0}) and consider the goodness-of-fit testing problem (\ref{testing_pb-2}). Let $\alpha \in ]0,1[$ and $\beta \in\, ]0,1-\alpha[$ be given. Then, for every $\varepsilon >0$, the minimax separation radius $\tilde r_{\varepsilon,0}$ is lower bounded by 
\begin{equation}
\label{eq:prop-s=0}
\tilde r^2_{\varepsilon,0} \geq \sup_{D \in \mathbb{N}} \left[c_{\alpha,\beta}\,
\varepsilon^2 \sqrt{\sum_{j=1}^D b_j^{-4}} \wedge a_D^{-2} \right],
\end{equation}
where
\begin{equation}
\label{eq:c-alpha-beta}
c_{\alpha,\beta}=(2 \ln (1+4(1-\alpha-\beta)^2))^{1/4}>0.
\end{equation}
\end{proposition}

\noindent 
The proof of the Proposition \ref{prop:1-s=0} with detailed arguments are related discussion can be found in e.g., \cite{Baraud}, \cite{LLM_2012} and \cite{MS_2014}.

\subsection{A Combined Lower Bound}
\label{s:lb3}

The following result provides a lower bound on the minimax separation radius $\tilde r_{\varepsilon,\sigma}$ for the goodness-of-fit testing problem (\ref{testing_pb1}). This lower bound corresponds to item $(ii)$ of Theorem \ref{thm:sum}.

\begin{thm}
\label{thm-comb-lb}
Consider the GSMs (\ref{1.0.0}), (\ref{eq:lower-bound-e=0}) and (\ref{eq:lower-bound-s=0}). Denote by $\tilde r_{\varepsilon,\sigma}$, $\tilde r_{0,\sigma}$ and $\tilde r_{\varepsilon,0}$ the corresponding minimax separation radii. Then, for every $\varepsilon >0$ and $\sigma >0$,
\begin{equation}
\label{eq:com-lower-bound-1}
\tilde r_{\varepsilon,\sigma} \geq \tilde r_{0,\sigma} \vee \tilde r_{\varepsilon,0}.
\end{equation}
In particular,
\begin{equation}
\label{eq:com-lower-bound-2}
\tilde r^2_{\varepsilon,\sigma} \geq  \left\{ \frac{C^2_{\alpha,\beta}}{16}\, \sigma^2\, 
\max_{1 \leq D \leq M_2} [ b_D^{-2} a_D^{-2} ] \right\} \vee \left\{ \sup_{D \in \mathbb{N}} \left[c_{\alpha,\beta}\,
\varepsilon^2 \sqrt{\sum_{j=1}^D b_j^{-4}} \wedge a_D^{-2} \right] \right\},
\end{equation}
where
$C_{\alpha,\beta}$ is given in (\ref{eq:C-alpha-beta}), $M_2$ is given in (\ref{eq:M2-def}) and 
$c_{\alpha,\beta}$  is given in (\ref{eq:c-alpha-beta}).
\end{thm}

\noindent 
The proof of Theorem \ref{thm-comb-lb} is postponed to Section \ref{s:prlower_m}.

\begin{remark}{\rm
At a first sight, the upper and lower bounds respectively displayed in (i) and (ii) of Theorem \ref{thm:sum} do not exactly match up. However, a closer look at the involved formulas indicates that both quantities contain terms that have similar behaviors. This is, in some sense, confirmed in Section  \ref{subsec:ALB} below where specific sequences $(a_j)_{j\in \mathbb{N}}$ and $(b_j)_{j\in \mathbb{N}}$ are treated. 
}
\end{remark}

\subsection{Lower Bounds: Specific Cases}
\label{subsec:ALB}

Our aim in this section is to determine an explicit value (in terms of the noise levels $\varepsilon$ and $\sigma$) for the lower bounds on the minimax separation radius $\tilde r_{\varepsilon, \sigma}$ obtained in Theorem \ref{thm-comb-lb} above for the specific sequences $(a_j)_{j\in \mathbb{N}}$ and $(b_j)_{j\in \mathbb{N}}$ considered in Section \ref{subsec:AUB}.

\begin{thm}
\label{thm:alower}
Consider the goodness-of-fit testing problem (\ref{testing_pb1}) when observations are given by (\ref{1.0.0}), and the signal of interest has smoothness governed by (\ref{eq:ellips_a}). Then, 
\begin{itemize}
\item[(i)] If $b_j \sim j^{-t}$, $t>0$, and $a_j \sim j^s$, $s>0$, for all $j \in \mathbb{N}$, then,
there exists $\varepsilon_0, \sigma_0 \in ]0,1[$ such that, for all $0 < \varepsilon \leq \varepsilon_0$ and $0 <\sigma \leq \sigma_0$, the minimax separation radius $\tilde r_{\varepsilon, \sigma}$ satisfies
\begin{equation*}
\tilde r^2_{\varepsilon, \sigma} \gtrsim \varepsilon^{\frac{4s}{2s+2t+1/2}}
\vee \sigma^{2 \left(\frac{s}{t} \wedge 1\right)}.
\end{equation*}

\item[(ii)] If $b_j \sim j^{-t}$, $t>0$, and $a_j \sim\exp\{js\}$, $s>0$, for all $j \in \mathbb{N}$, then,
there exists $\varepsilon_0, \sigma_0 \in ]0,1[$ such that, for all $0 < \varepsilon \leq \varepsilon_0$ and $0 <\sigma \leq \sigma_0$, the minimax separation radius $\tilde r_{\varepsilon, \sigma}$ satisfies
$$
\tilde r^2_{\varepsilon, \sigma} \gtrsim \varepsilon^2 \left[\ln\left(1/\varepsilon\right)\right]^{\left(2t+\frac{1}{2}\right)} \vee \sigma^{2}.
$$

\item[(iii)] If $b_j \sim \exp\{-jt\}$, $t>0$, and $a_j \sim j^s$, $s>0$, for all $j \in \mathbb{N}$, then,
there exists $\varepsilon_0, \sigma_0 \in ]0,1[$ such that, for all $0 < \varepsilon \leq \varepsilon_0$ and $0 <\sigma \leq \sigma_0$, the minimax separation radius $\tilde r_{\varepsilon, \sigma}$ satisfies
$$
\tilde r^2_{\varepsilon, \sigma} \gtrsim \left[\ln\left(1/\varepsilon\right)\right]^{-2s} \vee \left[\ln\left(1/\sigma \right)\right]^{-2s}.
$$

\item[(iv)] If $b_j \sim \exp\{-jt\}$, $t>0$, and $a_j \sim \exp\{js\}$, $s>0$, for all $j \in \mathbb{N}$, then,
there exists $\varepsilon_0, \sigma_0 \in ]0,1[$ such that, for all $0 < \varepsilon \leq \varepsilon_0$ and $0 <\sigma \leq \sigma_0$, the minimax separation radius $\tilde r_{\varepsilon, \sigma}$ satisfies
$$
\tilde r^2_{\varepsilon, \sigma} \gtrsim \varepsilon^{\frac{2s}{s+t}}
\vee \sigma^{2 \left(\frac{s}{t} \wedge 1\right)}.
$$
\end{itemize}
\end{thm}

The proof is postponed to Section \ref{s:p_rateslow}. As in the case of the upper bound, the main task is to compute the trade-off between both different antagonistic terms involved  in the lower bound on the minimax separation radius $\tilde r_{\varepsilon,\sigma}$ displayed in (\ref{eq:com-lower-bound-2}).

\section{Concluding Remarks}
\label{s:remarks}

The main conclusion of this work is that goodness-of-fit testing in an inverse problem setting is `feasible', even in the specific situation where some uncertainty is observed on the operator at hand in the model (\ref{1.0.0}). We have established `optimal' separation conditions for the goodness-of-fit testing problem  (\ref{testing_pb1}) via a sharp control of the associated minimax separation radius. \\

We stress that several outcomes and open questions are still of interest. We can mention, among others, 
\begin{itemize}
\item \textbf{Adaptivity:} As proved in Theorem \ref{cor:power-sco-test}, the test $\Psi_{D^\dagger,M}$ introduced in (\ref{eq:sco-tests-tatist})-(\ref{eq:ca}) with $D^\dagger$ defined in (\ref{eq:D-Dagger}) is powerful in the sense that its separation radius is equal (up to constant) to the minimax one. However, this test strongly depends on the sequence $a=(a_j)_{j\in\mathbb{N}}$ that characterizes the smoothness of the signal of interest. In practice, this sequence is unknown and \textit{adaptive} procedures are necessary (see, e.g., \cite{IS_2003} or \cite{ISS_2012}). 
\item \textbf{Signal detection:} We have already mentioned in Remark \ref{rm:sd} that \textit{signal detection} is different from \textit{goodness-of-fit testing}  (\ref{testing_pb1}) when the GSM (\ref{1.0.0}) is at hand. In this work, we were concerned with the case where $\theta_0 \neq 0$ (goodness-of-fit testing). However, some attention should also be paid in the future to the case where $\theta_0 = 0$ (signal detection). In particular, testing methodologies and related minimax separation radii are  quite different from those presented above.  
\item \textbf{Errors-in-variables model:} Density model with measurement errors have been at the core of several statistical studies in the past decades (see, e.g., \cite{meister} for an overview). Formally, given a sample of independent and identical distributed random variables  $(Y_i)_{i=1,2,\ldots,n}$ satisfying
$$ Y_i = X_i + \epsilon _i \quad i=1,2,\ldots, n, $$
the aim is to produce some inference on the {\em unknown} density of the $X_i$ denoted by $f$, the $\epsilon_i$ corresponding to some error, with {\em known} density $\eta$. This appears to be an inverse (deconvolution) problem since the $Y_i$ are associated to the convolved density $f*\eta$. In a goodness-of-fit testing task, this model has been discussed in \cite{Butucea} and minimax separation rates (in the asymptotic minimax testing framework) have been established in various settings. In the spirit of our contribution, it could be interesting to propose methods taking into account some possible uncertainty on the density $\eta$ at hand. 
\end{itemize}
 All these topics require special attention that is beyond the scope of this paper. Nevertheless, they provide an avenue for future research.


\section{Appendix: Proofs}
\label{s:appendix}

\subsection{Useful Lemmas}
\label{s:ulemma}
The constant $C>0$ and $0 < \tau <1$ below will vary from place to place.\\

\noindent
The following lemma is inspired by Lemma 6.1 of \cite{cavalier2}. 

\begin{lemma}
\label{lem:tech_lem_1}
Let $M, M_0, M_1$ be defined as in (\ref{eqn:M})-(\ref{eq:h_1}) where $\alpha\in ]0,1[$ and $\kappa \geq \exp(1)$ are fixed values.  Define the event 
\begin{equation}
\label{eq:M_set}
\mathcal{M}=\{M_0 \leq M < M_1\}.
\end{equation}
Then, for any $\sigma \in ]0,1[$,
\begin{equation}
\label{eq:M_set_prob}
\mathbb{P}(\mathcal{M}^c) \leq \frac{\alpha}{10} + \frac{\alpha\pi^2}{6\kappa}.
\end{equation}
\end{lemma}

\noindent 
\textsc{Proof of Lemma \ref{lem:tech_lem_1}.} It is easily seen that
\begin{eqnarray*}
\mathbb{P}(M \geq M_1) = \mathbb{P}\left(\bigcap_{j=1}^{M_1} \{|X_j|>\sigma h_j\}\right) &\leq& \mathbb{P}\left(|X_{M_{1}}| > \sigma h_{M_1} \right), \\
&\leq& \mathbb{P}\left(|b_{M_{1}}| + \sigma |\eta_{M_{1}}|> \sigma h_{M_1} \right), \\
&\leq& \mathbb{P}\left(|\eta_{M_{1}}|> h_{M_1}-h_{1,M_1}\right), \\
& = &  \mathbb{P}\left(|\eta_{M_{1}}|> \sqrt{2\ln\left( \frac{10}{\alpha}\right)} \right),
\end{eqnarray*}
where the sequences $(h_j)_{j\in \mathbb{N}}$ and $(h_{1,j})_{j\in \mathbb{N}}$  are defined in (\ref{eq:h}) and (\ref{eq:h_1}) respectively.  Using the bound
\begin{equation}
\frac{1}{2\pi} \int_x^{+\infty} e^{-\frac{x^2}{2}} dx \leq \frac{1}{x} \frac{e^{-\frac{x^2}{2}}}{\sqrt{2\pi}} \quad \forall x>0,
\label{eq:inequality_g}
\end{equation}
we get
\begin{equation}
\mathbb{P}(M \geq M_1) \leq  \frac{2}{\sqrt{2\pi}} \frac{\alpha}{10} \frac{1}{\sqrt{2\ln(10/\alpha)}} \leq \frac{\alpha}{10},
\label{eq:control_M1}
\end{equation}
since $\sqrt{2\ln(10/\alpha)} > 1$ for all $\alpha\in ]0,1[$. In the same spirit,
\begin{eqnarray*}
\mathbb{P}(M < M_0) = \mathbb{P}\left(\bigcup_{j=1}^{M_0} \{|X_j|\leq \sigma h_j\}\right) &\leq& \sum_{j=1}^{M_0}\mathbb{P}\left(|X_j| \leq \sigma h_j \right), \\
&\leq& \sum_{j=1}^{M_0}\mathbb{P}\left(|b_j| - \sigma |\eta_j|\leq \sigma h_j \right), \\
&\leq& \sum_{j=1}^{M_0} \mathbb{P}\left(\sigma |\eta_j|\geq |b_j| - \sigma h_j    )\right), \\
&\leq& \sum_{j=1}^{M_0} \mathbb{P}\left(|\eta_j|\geq  h_{0,j}- h_j\right). 
\end{eqnarray*}
According to the respective definition of $(h_j)_{j\in \mathbb{N}}$, $(h_{0,j})_{j\in \mathbb{N}}$ (see (\ref{eq:h}) and (\ref{eq:h_0})), and using again inequality (\ref{eq:inequality_g}), we obtain
\begin{eqnarray}
\mathbb{P}(M \geq M_1) 
&\leq& \frac{2}{\sqrt{2\pi}} \sum_{j=1}^{M_0} \frac{1}{h_{0,j}-h_j} \exp\left\{-\frac{1}{2}(h_{0,j}-h_j)^2\right\}, \nonumber\\
&\leq& \frac{2}{\sqrt{2\pi}} \sum_{j=1}^{M_0} \frac{1}{2\sqrt{\ln\left(\frac{\kappa j^2}{\alpha} \right)}} \frac{\alpha}{\kappa j^2}, \nonumber\\
&\leq& \frac{\alpha}{\kappa} \sum_{j \in \mathbb{N}} \frac{1}{j^2} = \frac{\alpha \pi^2}{6\kappa},
\label{eq:control_M0}
\end{eqnarray}
on noting that $\sum_{j \in \mathbb{N}} \frac{1}{j^2} =\pi^2/6$. Since 
$$
\mathbb{P}(\mathcal{M}^c) \leq \mathbb{P}(M < M_0) + \mathbb{P}(M \geq M_1),
$$
the lemma follows, thanks to (\ref{eq:control_M1}) and (\ref{eq:control_M0}). \begin{flushright} $\Box$ \end{flushright}

\begin{lemma}
\label{lem:tech_lem_2}
Let $M$ be defined as in (\ref{eqn:M}) and (\ref{eq:h_0}) where $\alpha\in ]0,1[$ and $\kappa \geq \exp(1)$ are fixed values. Define the event 
\begin{equation}
\label{eq:B_set}
\mathcal{B}=\bigcap_{j=1}^{M}\left\{\sigma |\eta_j| \leq \frac{b_j}{2}\right\}.
\end{equation}
Then, for any $\sigma \in ]0,1[$ 
\begin{equation}
\label{eq:M_set_prob}
\mathbb{P}(\mathcal{B}^c) \leq \frac{\alpha}{10}+ \frac{\alpha\pi^2}{3\kappa}.
\end{equation}
\end{lemma}

\noindent 
\textsc{Proof of Lemma \ref{lem:tech_lem_2}} Using the definitions of $\mathcal{M}$ and $M_1$, simple calculations give
\begin{eqnarray*}
\mathbb{P}(B^c)&=&\mathbb{P}(B^c \cap \mathcal{M})+\mathbb{P}(B^c \cap \mathcal{M}^c), \\
&\leq& \mathbb{P}\left(\bigcup_{j=1}^{M_1-1} \left\{\sigma |\eta_j| > \frac{b_j}{2} \right\} \right) +
\mathbb{P}(\mathcal{M}^c),\\
&\leq&\sum_{j=1}^{M_1-1} \mathbb{P}\left(|\eta_j| >\frac{1}{2}h_{1,j}\right) +
\mathbb{P}(\mathcal{M}^c).
\end{eqnarray*}
Using (\ref{eq:h_1}), Lemma \ref{lem:tech_lem_1} and (\ref{eq:inequality_g}), we obtain
\begin{eqnarray}
\mathbb{P}(B^c) &\leq& \frac{2}{\sqrt{2\pi}} \sum_{j=1}^{M_1} \frac{1}{\sqrt{8^2\ln\left(\frac{\kappa j^2}{\alpha} \right)}} \frac{\alpha}{\kappa j^2} + \frac{\alpha}{10} + \frac{\alpha\pi^2}{6\kappa} \nonumber \\
&\leq& \frac{\alpha}{10}+ \frac{\alpha\pi^2}{3\kappa}.
\end{eqnarray}
Hence, the lemma holds true. \begin{flushright} $\Box$ \end{flushright}

\begin{lemma}
\label{lem:tech_lem_3}
Let $\theta  \in \mathcal{E}_a$ be given. Let $M$ be defined as in (\ref{eqn:M}) and (\ref{eq:h_0}) where $\alpha\in ]0,1[$ and $\kappa \geq \exp(1)$ are fixed values. Then, for any $\sigma \in ]0,1[$ and for any $D \in \mathbb{N}$,
$$
\mathbb{P} \left(\sum_{j=1}^{D \wedge M}\left(\frac{b_j}{X_j} -1\right)^2\theta_{j}^2 \geq \sigma^2 \ln^{3/2}(1/\sigma) \vee a_{D \wedge M}^{-2}\right) \leq \frac{\alpha}{5} + \frac{\alpha}{\kappa}\left(\frac{\pi^2}{2} +2\right)+  C \exp\{-\ln^{1+\tau}(1/\sigma)\},
$$
for some $C>0$ and $0 <\tau <1$.
\end{lemma}

\noindent 
\textsc{Proof of Lemma \ref{lem:tech_lem_3}.} Using Lemma \ref{lem:tech_lem_1}, Lemma \ref{lem:tech_lem_2} and a Taylor expansion as in Lemma 6.6 of \cite{cavalier2}, we get, for all $j \leq M$,
$$
\frac{b_j}{X_j}=\frac{1}{1+\sigma b_j^{-1} \eta_j}=1-\sigma b_j^{-1} \eta_j +\sigma^2 \zeta_j^{-2}\eta_j^2,
$$
where $\zeta_{j}^{-1} \leq 8 b_j^{-1}$ on the even $\mathcal{B}$ defined in (\ref{eq:B_set}). Hence
\begin{eqnarray}
\lefteqn{ \mathbb{P} \left(\sum_{j=1}^{D \wedge M}\left(\frac{b_j}{X_j} -1\right)^2\theta_{j}^2 \geq \sigma^2 \ln^{3/2}(1/\sigma) \vee a_{D \wedge M}^{-2}\right)}\nonumber \\ 
&=& 
\mathbb{P} \left(\sum_{j=1}^{D \wedge M} (-\sigma b_j^{-1} \eta_j +\sigma^2 \zeta_j^{-2}\eta_j^2)^2 \theta_{j}^2 \geq \sigma^2 \ln^{3/2}(1/\sigma) \vee a_{D \wedge M}^{-2}\right), \nonumber \\
&\leq &
\mathbb{P} \left( 2 \sigma^2 \sum_{j=1}^{D \wedge M} b_j^{-2} \theta_{j}^2 \eta_j^2 + 2 \sigma^4  \sum_{j=1}^{D \wedge M} \zeta_j^{-4}\theta_{j}^2 \eta_j^4  \geq  \sigma^2 \ln^{3/2}(1/\sigma) \vee a_{D \wedge M}^{-2}\right). \nonumber
\end{eqnarray}
Therefore
\begin{eqnarray}
\lefteqn{ \mathbb{P} \left(\sum_{j=1}^{D \wedge M}\left(\frac{b_j}{X_j} -1\right)^2\theta_{j}^2 \geq \sigma^2 \ln^{3/2}(1/\sigma) \vee a_{D \wedge M}^{-2}\right)}  \nonumber \\ 
& \leq & \mathbb{P} \left( \left\lbrace 2 \sigma^2 \sum_{j=1}^{D \wedge M} b_j^{-2} \theta_{j}^2 \eta_j^2 \geq  \frac{1}{2} \left[\sigma^2 \ln^{3/2}(1/\sigma) \vee a_{D \wedge M}^{-2}\right] \right\rbrace \cap (\mathcal{B}\cap\mathcal{M})\right) \nonumber \\
& & \hspace{0cm}  +\mathbb{P} \left( \left\lbrace 2 \sigma^4  \sum_{j=1}^{D \wedge M} \zeta_j^{-4}\theta_{j}^2 \eta_j^4 \geq  \frac{1}{2} \left[\sigma^2 \ln^{3/2}(1/\sigma) \vee a_{D \wedge M}^{-2}\right] \right\rbrace \cap (\mathcal{B}\cap\mathcal{M})\right) + \mathbb{P}\left( (\mathcal{B}\cap\mathcal{M})^c\right), \nonumber \\
& := & T_1 + T_2 +   \mathbb{P}\left( (\mathcal{B}\cap\mathcal{M})^c\right).
\label{cyprus1}
\end{eqnarray}

We concentrate bellow our attention on the term $T_1$ defined as
$$
\mathbb{P} \left( \left\lbrace 2 \sigma^2 \sum_{j=1}^{D \wedge M} b_j^{-2} \theta_{j}^2 \eta_j^2 \geq  \frac{1}{2} \left[\sigma^2 \ln^{3/2}(1/\sigma) \vee a_{D \wedge M}^{-2}\right]\right\rbrace \cap (\mathcal{B}\cap\mathcal{M})\right).
$$ 
We consider the two following possible scenarios: (i) $a_j^{-1}b_j^{-1} \leq C_0$ as $j \rightarrow +\infty$, for some $C_0>0$, and (ii) $a_j^{-1}b_j^{-1} \rightarrow +\infty$ as $j \rightarrow +\infty$.\\

Consider first scenario (i). Then, using again (\ref{eq:inequality_g})
\begin{eqnarray}
T_1 & \leq & \mathbb{P} \left( \left\{ 2 \sigma^2 \max_{1\leq j \leq D \wedge M}(b_j^{-2} a_j^{-2}\eta_j^2) \geq \frac{1}{2} \sigma^2 \ln^{3/2}(1/\sigma) \vee a_{D \wedge M}^{-2}\right\} \cap (\mathcal{B}\cap\mathcal{M})\right) \nonumber \\
& \leq & 
 \mathbb{P} \left( \left\{ 2 C_0^2 \sigma^2 \max_{1\leq j \leq D \wedge M}(\eta_j^2) \geq \frac{1}{2} \sigma^2 \ln^{3/2}(1/\sigma) \right\} \cap (\mathcal{B}\cap\mathcal{M})\right)
\nonumber \\
& \leq & 
 \sum_{j=1}^{M_1-1} \mathbb{P} \left(\eta_j \geq \frac{1}{2C_0}  \ln^{3/4}(1/\sigma)  \right), \nonumber\\
 & \leq & \frac{2 M_1 }{\sqrt{2\pi}} \frac{2C_0}{ \ln^{3/4}(1/\sigma) } \exp \left(- \frac{C\ln^{3/2}(1/\sigma)}{8C_0^2} \right) \leq  C \exp\{-\ln^{1+\tau}(1/\sigma)\}.
\end{eqnarray}
for some constants $C,\tau \in \mathbb{R}^+$. A similar bound occurs for the term $T_2$ for this scenario. 

Consider now the second scenario (ii). Then
\begin{eqnarray}
T_1 &\leq & 
\mathbb{P} \left( \left\{ 2 b_{D \wedge M}^{-2} a_{D \wedge M}^{-2} \sigma^2 \max_{1\leq j \leq D \wedge M}(\eta_j^2) \geq \frac{1}{2} a_{D \wedge M}^{-2}\right\} \cap (\mathcal{B}\cap\mathcal{M})\right)
\nonumber \\
& \leq &
\mathbb{P} \left( \left\{ \sigma^2 \max_{1\leq j \leq D \wedge M}(\eta_j^2) \geq  \frac{1}{4} b_{D \wedge M}^{2} \right\} \cap (\mathcal{B}\cap\mathcal{M})\right)
\nonumber \\ 
& \leq &
\sum_{j=1}^{M_1-1} \mathbb{P} \left(\sigma^2 \eta_j^2 \geq  \frac{1}{4} b_{M_1-1}^{2}   \right), \nonumber 
\end{eqnarray}
since the sequence $(b_j)_{j\in \mathbb{N}}$ is non-increasing. Using (\ref{eq:inequality_g}), we get
\begin{eqnarray}
T_1 &\leq& \sum_{j=1}^{M_1-1} \mathbb{P} \left( \eta_j \geq  \frac{1}{2} h_{1,M_1-1}\right), \nonumber \\ 
&\leq& \frac{2M_1}{\sqrt{2\pi}} \frac{2}{h_{1,M_1-1}} \exp\left( - \frac{h_{1,M_1-1}^{2}}{8} \right), \nonumber\\ 
& \leq &  M_1 \exp\left( -   \ln\left( \frac{\kappa M_1^2}{\alpha} \right) \right) \nonumber \\
& \leq & M_1 \times \frac{\alpha}{\kappa M_1^2} \leq \frac{\alpha}{\kappa}.
\label{cyprus2}
\end{eqnarray}
By similar computations, we get 
\begin{eqnarray}
T_2 &:=& \mathbb{P} \left( \left\lbrace 2 \sigma^4  \sum_{j=1}^{D \wedge M} \zeta_j^{-4}\theta_{j}^2 \eta_j^4 \geq  \frac{1}{2} \sigma^2 \ln^{3/2}(1/\sigma) \vee a_{D \wedge M}^{-2}\right\rbrace \cap (\mathcal{B}\cap\mathcal{M})\right), \nonumber  \\
& \leq & \mathbb{P} \left( \left\lbrace 2 \times 8^4 \sigma^4  \sum_{j=1}^{D \wedge M} b_j^{-4}\theta_{j}^2 \eta_j^4 \geq  \frac{1}{2}  a_{D \wedge M}^{-2} \right\rbrace \cap (\mathcal{B}\cap\mathcal{M})\right), \nonumber \\
& \leq & \mathbb{P} \left( \left\lbrace 2 \times 8^4 \sigma^4  \max_{j=1..D \wedge M}  \eta_j^4 \geq  \frac{1}{2}  b_{D \wedge M}^4 \right\rbrace \cap (\mathcal{B}\cap\mathcal{M})\right), \nonumber  \\
& \leq & \sum_{j=1}^{M_1-1} \mathbb{P} \left( |\eta_j | \geq \frac{1}{8 \sqrt{2}} h_{1,M_1-1} \right), \nonumber \\
& \leq & \frac{2}{\sqrt{2\pi}} \frac{8\sqrt{2} M_1}{h_{1,M_1}} \exp\left( - \frac{1}{4 \times 8^2}  4\times 8^2 \ln \left( \frac{\kappa M_1^2}{\alpha} \right) \right) \leq \frac{\alpha}{\kappa}.
\label{cyprus3}
\end{eqnarray}
Hence, the lemma follows from Lemmas \ref{lem:tech_lem_1}, \ref{lem:tech_lem_2} and (\ref{cyprus1})-(\ref{cyprus3}). \begin{flushright} $\Box$ \end{flushright}

\begin{lemma}
\label{lem:tech_lem_4}
Let 
$$
Z_j = \nu_j + v_j \omega_j, \quad j \in \mathbb{N},
$$
where $\omega=(\omega_j)_{j \in \mathbb{N}}$ is a sequence of independent standard Gaussian random variables. For all $D\in \mathbb{N}$, define
$$
T=\sum_{j=1}^D Z_j^2 \quad \text{and} \quad \Sigma=\sum_{j=1}^D v_j^4 + 2 \sum_{j=1}^D v_j^2 \nu_j^2.
$$
Then, for all $x>0$,
\begin{eqnarray}
\label{eq:LLM-1}
\mathbb{P} \left(T -\mathbb{E}(T) > 2 \sqrt{\Sigma x} + 2x \sup_{1\leq j\leq D}(v_j^2) \right) &\leq& \exp(-x) \\
\label{eq:LLM-2}
\mathbb{P} \left(T -\mathbb{E}(T) < -2 \sqrt{\Sigma x}  \right) &\leq& \exp(-x).
\end{eqnarray}
\end{lemma}

\noindent 
\textsc{Proof of Lemma \ref{lem:tech_lem_4}} The proof is given in Lemma 2 of \cite{LLM_2012}.  \begin{flushright} $\Box$ \end{flushright}

\subsection{Non-Asymptotic Upper bounds}
\label{s:upper_bounds}

\subsubsection{Proof of Proposition \ref{prop:size-sco-test}} 
\label{s:proof1}

By definition,
$$
\Alpha_{\varepsilon,\sigma}(\Psi_{D,M}):=\mathbb{P}_{\theta_0,b}(\Psi_{D,M}=1)
=\mathbb{P}_{\theta_0,b}(T_{D,M} >
t_{1-\alpha,D}(X)).$$
Conditionally to the sequence $X=(X_j)_{j \in \mathbb{N}}$, for each $1 \leq j \leq 
D \wedge M$, the random variable $X_j^{-1} Y_j-\theta_{j,0}$ is Gaussian with mean $\nu_j=(b_j/X_j -1)\theta_{j,0}$ and standard deviation $v_j=\varepsilon X_j^{-1}$. In particular, for all $D\in \mathbb{N}$
\begin{equation}
\mathbb{E}_{\theta_0,b}(T_{D,M} \mid X):=
\mathbb{E}_{\theta_0,b}\left[\sum_{j=1}^{D \wedge M}\left(\frac{Y_j}{X_j}-\theta_{j,0}\right)^2 \mid  X \right]=
\sum_{j=1}^{D \wedge M}\left(\frac{b_j}{X_j} -1\right)^2\theta^2_{j,0}+\varepsilon^2 \sum_{j=1}^{D \wedge M} X_j^{-2}.
\label{eq:more}
\end{equation}
For all $D\in \mathbb{N}$, define
$$
\Sigma_{D,M}:= \varepsilon^4 \sum_{j=1}^{D \wedge M} X_j^{-4} + 
 \varepsilon^2 \sum_{j=1}^{D \wedge M}X_j^{-2}\left(\frac{b_j}{X_j} -1\right)^2\theta^2_{j,0}.
$$
Applying Lemma \ref{lem:tech_lem_4} with $T=T_{D,M}$, $\Sigma=\Sigma_{D,M}$ and $x=x_{\alpha/2}:=\ln(2/\alpha)$, we get 
\begin{equation}
\mathbb{P}_{\theta_0,b} \left(T_{D,M} -\mathbb{E}_{\theta_0}(T_{D,M} \mid X) > 2 \sqrt{\Sigma_{D,M} x_{\alpha/2}} + 2\varepsilon^2 x_{\alpha/2} \max_{1\leq j\leq D \wedge M}(X_j^{-2}) \mid X \right) \leq \frac{\alpha}{2}.
\label{eq:interm}
\end{equation}
Using the inequalities $\sqrt{a+b} \leq \sqrt{a}+\sqrt{b}$ and $ab \leq a^2/2+b^2/2$ for $a,b>0$, it is easily seen that
\begin{eqnarray}
\sqrt{\Sigma_{D,M}} &\leq& \varepsilon^2  \sqrt{\sum_{j=1}^{D \wedge M} X_j^{-4}}
+ \sqrt{\varepsilon^2 \sum_{j=1}^{D \wedge M}X_j^{-2}\left(\frac{b_j}{X_j} -1\right)^2\theta^2_{j,0}} \nonumber \\
&\leq& \varepsilon^2  \sqrt{\sum_{j=1}^{D \wedge M} X_j^{-4}} +
\sqrt{\varepsilon^2 \max_{1 \leq j \leq D \wedge M} X_j^{-2} \sum_{j=1}^{D \wedge M}\left(\frac{b_j}{X_j} -1\right)^2\theta^2_{j,0}}\nonumber \\
&\leq& \varepsilon^2  \sqrt{\sum_{j=1}^{D \wedge M} X_j^{-4}} +
\frac{1}{2}\varepsilon^2 \max_{1 \leq j \leq D \wedge M} X_j^{-2} +
\frac{1}{2} \sum_{j=1}^{D \wedge M}\left(\frac{b_j}{X_j} -1\right)^2\theta^2_{j,0}.
\label{eq:sigma-bound}
\end{eqnarray}
According to (\ref{eq:more})-(\ref{eq:sigma-bound}), we obtain the following bound 
$$
\mathbb{P}_{\theta_0,b} \left(T_{D,M} > (1+\sqrt{x_{\alpha/2}}) \sum_{j=1}^{D \wedge M}\left(\frac{b_j}{X_j} -1\right)^2\theta^2_{j,0} + \varepsilon^2 \sum_{j=1}^{D \wedge M} X_j^{-2} + C(\alpha)  \varepsilon^2  \sqrt{\sum_{j=1}^{D \wedge M} X_j^{-4}}  \mid X \right) \leq \frac{\alpha}{2},
$$
where the constant $C(\alpha)$ is defined in (\ref{eq:ca}).
Since $\mathbb{E}[\mathbb{E}(V \mid W)]=\mathbb{E}(V)$ for any random variables $V$ and $W$, the previous inequality leads to
$$
\mathbb{P}_{\theta_0,b} \left(T_{D,M} > (1+\sqrt{x_{\alpha/2}}) \sum_{j=1}^{D \wedge M}\left(\frac{b_j}{X_j} -1\right)^2\theta^2_{j,0} + \varepsilon^2 \sum_{j=1}^{D \wedge M} X_j^{-2} + C(\alpha)  \varepsilon^2  \sqrt{\sum_{j=1}^{D \wedge M} X_j^{-4}} \right) \leq \frac{\alpha}{2}.
$$ 
Then, by defining 
$$
\mathcal{A} = \left\{  \sum_{j=1}^{D \wedge M}\left(\frac{b_j}{X_j} -1\right)^2\theta_{j,0}^2 < \sigma^2 \ln^{3/2}(1/\sigma) \vee a_{D \wedge M}^{-2}   \right\},
$$
and applying Lemma \ref{lem:tech_lem_3}, we immediately get
\begin{eqnarray*}
\Alpha_{\varepsilon,\sigma}(\Psi_{D,M}) &\leq& \mathbb{P}_{\theta_0,b}(\{T_{D,M} >
t_{1-\alpha,D}(X)\} \cap \mathcal{A}) + \mathbb{P}(\mathcal{A}^c) \nonumber \\
& \leq &\frac{\alpha}{2} + \frac{\alpha}{5} + \frac{\alpha}{6\kappa} (3\pi^2 +12)+ C \exp\{-\ln^{1+\tau}(1/\sigma),\\
& = & \frac{7\alpha}{10} + \frac{\alpha}{\kappa} \left(\frac{\pi^2}{2} +2\right))+ C \exp\{-\ln^{1+\tau}(1/\sigma),
\end{eqnarray*}
for some $C>0$ and $0<\tau <1$.  In particular, setting
$$ \kappa = 5 \left(\frac{\pi^2}{2} + 2\right),$$
there exists $\sigma_0 \in ]0,1[$ such that, for all $\sigma \leq \sigma_0$ and for each $\varepsilon >0$, 
$$
\Alpha_{\varepsilon,\sigma}(\Psi_{D,M}) \leq \alpha.
$$
This concludes the proof of the proposition. \begin{flushright} $\Box$ \end{flushright}
\medskip

\subsubsection{Proof of Proposition \ref{thm:power-sco-test}} 
\label{s:proof2}

Let $\theta, \theta_0 \in \mathcal{E}_a$ and $\theta-\theta_0\in \Theta_a(r_{\varepsilon,\sigma})$. Then
\begin{eqnarray}
\label{eqn:power-sco-test-F-1}
\mathbb{P}_{\theta,b}(\Psi_{D,M}=0)&=&\mathbb{P}_{\theta,b}(\{\Psi_{D,M}=0\}\cap(\mathcal
{B}\cap\mathcal{M})) + \mathbb{P}_{\theta,b}(\{\Psi_{D,M}=0\}\cap(\mathcal
{B}\cap\mathcal{M})^c)\nonumber\\
&:=&T_1+T_2.
\end{eqnarray}

\noindent
{\em Control of $T_2$:} Using Lemma \ref{lem:tech_lem_1}, Lemma \ref{lem:tech_lem_2} and elementary probabilistic arguments, we get
\begin{eqnarray}
\label{eqn:control-2}
T_2:= \mathbb{P}_{\theta,b}(\{\Psi_{D,M}=0\}\cap(\mathcal{B}\cap\mathcal{M})^c) &\leq&
\mathbb{P}((\mathcal{B}\cap\mathcal{M})^c) \nonumber\\
&\leq& \mathbb{P}(\mathcal{B}^c) + \mathbb{P}(\mathcal{M}^c) \nonumber\\
&\leq& \frac{\alpha}{5} +  \frac{\alpha}{\kappa} \left(\frac{\pi^2}{2} + 2 \right)  \leq  \frac{\beta}{5} +  \frac{\beta}{\kappa} \left(\frac{\pi^2}{2} + 2 \right),
\end{eqnarray}
since $\beta > \alpha$. \\

\noindent
{\em Control of $T_1$:} Define $t_{\beta/2,D}(\theta,X)$ to be the $\beta/2$-quantile of $T_{D,M}$, conditionally on $X$, i.e.,
$$
\mathbb{P}_{\theta,b}(T_{D,M} \leq t_{\beta/2,D}(\theta,X) \mid X) \leq \frac{\beta}{2}.
$$
Then, by elementary probabilistic arguments, we get
\begin{eqnarray}
\label{eqn:control-1}
T_1&:=& \mathbb{P}_{\theta,b}(\{\Psi_{D,M}=0\}\cap \{\mathcal
{B}\cap\mathcal{M}\}), \nonumber\\ 
&=& \mathbb{E} \big[ \mathbb{P}_{\theta,b}\left(\{\Psi_{D,M}=0\} \mid X\right) \mathbf{1}\{\mathcal{B}\cap\mathcal{M}\}\big], \nonumber\\
&=&\mathbb{E} \big[  \mathbb{P}_{\theta,b}\left(T_{D,M} \leq t_{1-\alpha,D}(X) \mid X\right) \mathbf{1}\{\mathcal{B}\cap\mathcal{M}\}\big], \nonumber\\
&\leq& \frac{\beta}{2}\; \mathbb{E} \big[ \mathbf{1}\{t_{1-\alpha,D}(X) \leq t_{\beta/2,D}(\theta,X)\}\mathbf{1}\{\mathcal{B}\cap\mathcal{M}\} \big] \nonumber \\
&& \hspace{-0.3cm}+\hspace{0.3cm}
\mathbb{E} \big[ \mathbf{1}\{t_{1-\alpha,D}(X) > t_{\beta/2,D}(\theta,X)\}\mathbf{1}\{\mathcal{B}\cap\mathcal{M}\} \big] ,\nonumber\\
&\leq& \frac{\beta}{2} + \mathbb{E} \big[ \mathbf{1}\{t_{1-\alpha,D}(X) > t_{\beta/2,D}(\theta,X)\}\mathbf{1}\{\mathcal{B}\cap\mathcal{M}\} \big], \nonumber\\
&\leq& \frac{\beta}{2} + \mathbb{P}_{\theta,b} \big( \{t_{1-\alpha,D}(X) > t_{\beta/2,D}(\theta,X)\} \cap 
\{\mathcal{B}\cap\mathcal{M}\}\big). 
\end{eqnarray}
Our next task is to provide a lower bound for $t_{\beta/2,D}(\theta,X)$.
Under $H_1$, conditionally to the sequence $X=(X_j)_{j \in \mathbb{N}}$, for each $1 \leq j \leq 
D \wedge M$, the random variable $X_j^{-1} Y_j-\theta_{j,0}$ is Gaussian with mean $\nu_j$ and standard deviation $v_j$ defined as
$$\nu_j=\left(\frac{b_j}{X_j} -1 \right)\theta_{j}+(\theta_j-\theta_{j,0}) \quad \mathrm{and} \quad v_j=\varepsilon X_j^{-1}.$$
In particular, 
\begin{eqnarray}
\mathbb{E}_{\theta,b}(T_{D,M} \mid X) &=&
\sum_{j=1}^{D \wedge M}\left[ \left(\frac{b_j}{X_j} -1\right)\theta_{j}+
(\theta_j-\theta_{j,0})\right]^2 +
\varepsilon^2 \sum_{j=1}^{D \wedge M} X_j^{-2} \nonumber \\
&=&
\sum_{j=1}^{D \wedge M}\nu_j^2 +
\varepsilon^2 \sum_{j=1}^{D \wedge M} X_j^{-2}.
\end{eqnarray}
Let
\begin{eqnarray}
\tilde \Sigma_{D,M} &:=& \varepsilon^4 \sum_{j=1}^{D \wedge M} X_j^{-4} + 
 \varepsilon^2 \sum_{j=1}^{D \wedge M}X_j^{-2}\left[ \left(\frac{b_j}{X_j} -1\right)\theta_{j}+
(\theta_j-\theta_{j,0})\right]^2 \nonumber \\
&=& 
\varepsilon^4 \sum_{j=1}^{D \wedge M} X_j^{-4} + 
 \varepsilon^2 \sum_{j=1}^{D \wedge M}X_j^{-2}\nu_j^2.
\end{eqnarray}
Using Lemma \ref{lem:tech_lem_4} with $T=T_{D,M}$, $\Sigma=\tilde \Sigma_{D,M}$ and $x=x_{\beta/2}:=\ln(2/\beta)$, we obtain
\begin{eqnarray}
&& \mathbb{P}_{\theta,b} \left( T_{D,M} < \sum_{j=1}^{D \wedge M}\nu_j^2 +
\varepsilon^2 \sum_{j=1}^{D \wedge M} X_j^{-2} - 2 \sqrt{\tilde \Sigma_{D,M} {x_{\beta/2}}}   \mid X \right) \leq \frac{\beta}{2} \nonumber \\
&\Rightarrow& t_{\beta/2,D}(\theta,X) \geq \sum_{j=1}^{D \wedge M}\nu_j^2 +
\varepsilon^2 \sum_{j=1}^{D \wedge M} X_j^{-2} - 2 \sqrt{\tilde \Sigma_{D,M} {x_{\beta/2}}}.
\label{eq:b-quat-lower}
\end{eqnarray}
Therefore, using (\ref{eqn:comp-a-quantile}) and (\ref{eq:b-quat-lower}), we get
\begin{eqnarray}
&& \mathbb{P}_{\theta,b} \big( \{t_{1-\alpha,D}(X) > t_{\beta/2,D}(\theta,X)\} \cap 
\{\mathcal{B}\cap\mathcal{M}\}\big) \nonumber\\ 
&\leq&
\mathbb{P}_{\theta,b} \left(\left\{ \sum_{j=1}^{D \wedge M} \nu_j^2 < \left(  C(\alpha) + 2\sqrt{x_{\beta/2}} \right)  
\varepsilon^2 \sqrt{\sum_{j=1}^{D \wedge M}  X_j^{-4}} \right. \right. \nonumber \\
&& \left. \left.
+ \; (1+\sqrt{x_{\alpha/2}}) \left[\sigma^2 \ln^{3/2}(1/\sigma) \vee a^{-2}_{D \wedge M}\right]  +2\sqrt{x_{\beta/2}} \sqrt{\varepsilon^2 \sum_{j=1}^{D \wedge M}  X_j^{-2}\nu_j^2} \right\} \cap
\{\mathcal{B}\cap\mathcal{M}\} \right)
\nonumber\\
&\leq&
\mathbb{P}_{\theta,b} \left(\left\{\frac{1}{2} \sum_{j=1}^{D \wedge M} \nu_j^2 < C(\alpha,\beta) 
\varepsilon^2 \sqrt{\sum_{j=1}^{D \wedge M}  X_j^{-4}} + (1+\sqrt{x_{\alpha/2}}) [\sigma^2 \ln^{3/2}(1/\sigma) \vee a^{-2}_{D \wedge M}]  \right\}\cap 
\{\mathcal{B}\cap\mathcal{M}\} \right), \nonumber
\end{eqnarray}
where 
\begin{equation}
C(\alpha,\beta):= C(\alpha)+ 3\sqrt{x_{\beta/2}},
\label{eq:cab}
\end{equation} 
and $C(\alpha)$ is defined in (\ref{eq:ca}). Note that, for any $a,b \in \mathbb{R}$, using the Young inequality $2ab \leq \gamma a^2 + \gamma^{-1} b^2$  for $\gamma =1/2$ we get $(a+b)^2 \geq a^2/2-b^2$. Applying the latter inequality with
$$
a=\theta_j-\theta_{j,0}, \quad 
b= \left(\frac{b_j}{X_j} -1\right)\theta_{j}, \quad j=1,\ldots,D \wedge M,
$$
and using Lemma \ref{lem:tech_lem_3}, we arrive at
\begin{eqnarray}
&&\mathbb{P}_{\theta,b} \big( \{t_{1-\alpha,D}(X) > t_{\beta/2,D}(\theta,X)\} \cap 
\{\mathcal{B}\cap\mathcal{M}\}\big) \nonumber\\
&\leq&
\mathbb{P}_{\theta,b} \left(\left\{\sum_{j=1}^{D \wedge M} (\theta_j-\theta_{j,0})^2 < 4C(\alpha,\beta) 
\varepsilon^2 \sqrt{\sum_{j=1}^{D \wedge M}  X_j^{-4}}  \right. \right. \nonumber \\
&& \left. \left. + \, 4(1+\sqrt{x_{\alpha/2}}) \left[\sigma^2 \ln^{3/2}(1/\sigma) \vee a^{-2}_{D \wedge M}\right]  + 2\sum_{j=1}^{D \wedge M} \left(\frac{b_j}{X_j} -1\right)^2\theta_{j}^2 \right\}\cap \{\mathcal{B}\cap\mathcal{M}\} \right) \nonumber\\
&\leq&
\mathbb{P}_{\theta,b} \left(\left\{\sum_{j=1}^{D \wedge M} (\theta_j-\theta_{j,0})^2 < 4C(\alpha,\beta) 
\varepsilon^2 \sqrt{\sum_{j=1}^{D \wedge M}  X_j^{-4}} + (6+4\sqrt{x_{\alpha/2}}) \left[\sigma^2 \ln^{3/2}(1/\sigma) \vee a^{-2}_{D \wedge M}\right]  \right\}\cap \{\mathcal{B}\cap\mathcal{M}\} \right), \nonumber\\
& & \hspace{1cm} + \frac{\alpha}{5} + \frac{\alpha}{\kappa} \left(\frac{\pi^2}{2} +2\right)+  C \exp\{-\ln^{1+\tau}(1/\sigma)\}.     \nonumber
\end{eqnarray}
Using the fact that $\theta \in \mathcal{E}_a$, we get 
\begin{eqnarray} 
&&\mathbb{P}_{\theta,b} \big( \{t_{1-\alpha,D}(X) > t_{\beta/2,D}(\theta,X)\} \cap \{\mathcal{B}\cap\mathcal{M}\}\big) \nonumber\\
&\leq&
\mathbb{P}_{\theta,b} \left(\left\{\|\theta-\theta_0\|^2 < 4C(\alpha,\beta) 
\varepsilon^2 \sqrt{\sum_{j=1}^{D \wedge M}  X_j^{-4}} \right. \right. \nonumber \\
&& \left. \left. + \, (6+4\sqrt{x_{\alpha/2}}) \left[\sigma^2 \ln^{3/2}(1/\sigma) \vee a^{-2}_{D \wedge M}\right]  + \sum_{j> D \wedge M}(\theta_j-\theta_{j,0})^2\right\}\cap \{\mathcal{B}\cap\mathcal{M}\} \right)  \nonumber \\
& & \hspace{2cm} + \frac{\alpha}{5} + \frac{\alpha}{\kappa} \left(\frac{\pi^2}{2} +2\right)+  C \exp\{-\ln^{1+\tau}(1/\sigma)\}  \nonumber\\
&\leq&
\mathbb{P}_{\theta,b} \left(\left\{\|\theta-\theta_0\|^2 < 4C(\alpha,\beta) 
\varepsilon^2 \sqrt{\sum_{j=1}^{D \wedge M}  X_j^{-4}} + (7+4\sqrt{x_{\alpha/2}}) \left[\sigma^2 \ln^{3/2}(1/\sigma) \vee a^{-2}_{D \wedge M}\right]  \right\}\cap \{\mathcal{B}\cap\mathcal{M}\} \right) \nonumber\\
& & \hspace{2cm} + \frac{\alpha}{5} + \frac{\alpha}{\kappa} \left(\frac{\pi^2}{2} +2\right)+  C \exp\{-\ln^{1+\tau}(1/\sigma)\}.
\label{eq:tab-control}
\end{eqnarray}
To conclude the proof, note that on the event $\{\mathcal{B}\cap\mathcal{M}\}$, we have
\begin{equation}
M_0 \leq M < M_1 \quad \text{and} \quad \frac{b_j}{X_j} \in \left[\frac{2}{3}, 2\right] \quad \forall \; j=1,\ldots,M.
\label{eq:M0-M1}
\end{equation}
Hence, using (\ref{eq:tab-control}) and (\ref{eq:M0-M1})
\begin{eqnarray}
&&\mathbb{P}_{\theta,b} \big( \{t_{1-\alpha,D}(X) > t_{\beta/2,D}(\theta,X)\} \cap 
\{\mathcal{B}\cap\mathcal{M}\}\big) \nonumber\\
&\leq&
\mathbb{P}_{\theta,b} \left( \|\theta-\theta_0\|^2 < 16C(\alpha,\beta) 
\varepsilon^2 \sqrt{\sum_{j=1}^{D \wedge M_1}  b_j^{-4}} + (7+4\sqrt{x_{\alpha/2}}) \left[\sigma^2 \ln^{3/2}(1/\sigma) \vee a^{-2}_{D \wedge M_0}\right]  \right) \nonumber \\
& & \hspace{2cm} + \frac{\alpha}{5} + \frac{\alpha}{\kappa} \left(\frac{\pi^2}{2} +2\right)+  C \exp\{-\ln^{1+\tau}(1/\sigma)\}, \nonumber \\
& = & \frac{\alpha}{5} + \frac{\alpha}{\kappa} \left( \frac{\pi^2}{2} +2\right)+  C \exp\{-\ln^{1+\tau}(1/\sigma)\}, \nonumber
\end{eqnarray}
as soon as
\begin{equation}
\|\theta-\theta_0\|^2 \geq \tilde{C}(\alpha,\beta) 
\varepsilon^2 \sqrt{\sum_{j=1}^{D \wedge M_1}  b_j^{-4}} + (7+4\sqrt{x_{\alpha/2}}) \left[\sigma^2 \ln^{3/2}(1/\sigma) \vee a^{-2}_{D \wedge M_0}\right],
\label{eq:final-norm-1}
\end{equation}
where $\tilde C(\alpha,\beta)=16 C(\alpha,\beta)$ is defined in (\ref{eq:constant_cyprus}). Therefore, for any fixed $\beta \in ]\alpha,1[$, (\ref{eq:final-norm-1}) implies that, there exists $\sigma_0 \in ]0,1[$ such that, for all $0<\sigma<\sigma_0$ and for each $\varepsilon >0$,
$$
\mathbb{P}_{\theta,b}(\Psi_{D,M}=0) \leq  \frac{7\beta}{10} + \frac{\beta}{\kappa} \left(\frac{\pi^2}{2} +2\right)+ C \exp\{-\ln^{1+\tau}(1/\sigma)\} \leq \beta,
$$
for some $C>0$ and $0<\tau <1$, which, in turn, implies that (\ref{eq:re-ub-1}) holds true. The last part of the theorem is a direct consequence of (\ref{eq:type2}) and (\ref{eq:re-ub-1}). This completes the proof of the proposition. \begin{flushright} $\Box$ \end{flushright}

\subsubsection{Proof of Theorem \ref{cor:power-sco-test}} 
\label{s:proof3}

The validity of (\ref{eq:alpha-level-cor}) can be immediately derived from Proposition \ref{prop:size-sco-test} taking into account that Lemma \ref{lem:tech_lem_3} is still valid with $D:=D^\dagger$ (that depends on the sequence $X=(X_j)_{j \in \mathbb{N}}$). For the proof of (\ref{eq:re-ub-1}), note first that (\ref{eqn:power-sco-test-F-1}), (\ref{eqn:control-2}) and (\ref{eqn:control-1}) still holds true with $D:=D^\dagger$. In the same spirit, is is easy to see that Lemma \ref{lem:tech_lem_3} is still valid when the bandwidth $D$ is measurable with respect to the sequence $(X_k)_{k\in \mathbb{N}}$. Hence, the same inequality than (\ref{eq:tab-control}) can be obtained with $D:=D^\dagger$, namely
\begin{eqnarray*}
&&\mathbb{P}_{\theta,b} \big( \{t_{1-\alpha,D^\dagger}(X) > t_{\beta/2,D^\dagger}(\theta,X)\} \cap 
\{\mathcal{B}\cap\mathcal{M}\}\big) \nonumber\\
&\leq&
\mathbb{P}_{\theta,b} \left(\left\{\|\theta-\theta_0\|^2 < 4C(\alpha,\beta) 
\varepsilon^2 \sqrt{\sum_{j=1}^{D^\dagger \wedge M}  X_j^{-4}} + (7+4\sqrt{x_{\alpha/2}})\left[\sigma^2 \ln^{3/2}(1/\sigma) \vee a^{-2}_{D^\dagger \wedge M}\right]  \right\}\cap \{\mathcal{B}\cap\mathcal{M}\} \right)\nonumber \\
&\leq&
\mathbb{P}_{\theta,b} \left(\|\theta-\theta_0\|^2 < \inf_{D \in \mathbb{N}}\left[16C(\alpha,\beta) 
\varepsilon^2 \sqrt{\sum_{j=1}^{D \wedge M_1}  b_j^{-4}} + (7+4\sqrt{x_{\alpha/2}}) \left[\sigma^2 \ln^{3/2}(1/\sigma) \vee a^{-2}_{D \wedge M_0}\right] \right]\right) \\
&=0,&
\label{eq:tab-control-2}
\end{eqnarray*}
as soon as
$$
\|\theta-\theta_0\|^2 \geq \inf_{D \in \mathbb{N}}\left[\tilde C(\alpha,\beta) 
\varepsilon^2 \sqrt{\sum_{j=1}^{D \wedge M_1}  b_j^{-4}} + (7+4\sqrt{x_{\alpha/2}}) \left [\sigma^2 \ln^{3/2}(1/\sigma) \vee a^{-2}_{D \wedge M_0}\right] \right],
$$
where $\tilde C(\alpha,\beta)$ is defined in (\ref{eq:constant_cyprus}).  Therefore, we immediately get that (\ref{eq:re-ub-1}) holds true.Finally, the validity of (\ref{eq:min-sep-rate-cor}) follows immediately on noting that 
\begin{eqnarray*}
\tilde r^2_{\varepsilon,\sigma} &:=& \inf_{\tilde \Psi_\alpha:\, \Alpha_{\varepsilon,\sigma}(\tilde \Psi_\alpha)\leq \alpha} r^2_{\varepsilon,\sigma}(\mathcal{E}_a, \tilde\Psi_\alpha,\beta) \nonumber\\
&\leq&
r^2_{\varepsilon,\sigma}(\mathcal{E}_a, \Psi_{D^\dagger,M},\beta) \nonumber\\
&\leq&
\inf_{D \in \mathbb{N}}\left[\tilde C(\alpha,\beta) 
\varepsilon^2 \sqrt{\sum_{j=1}^{D \wedge M_1}  b_j^{-4}} + (7+4\sqrt{x_{\alpha/2}}) \left[\sigma^2 \ln^{3/2}(1/\sigma) \vee a^{-2}_{D \wedge M_0}\right] \right].
\end{eqnarray*}
This completes the proof of the theorem. \begin{flushright} $\Box$ \end{flushright}

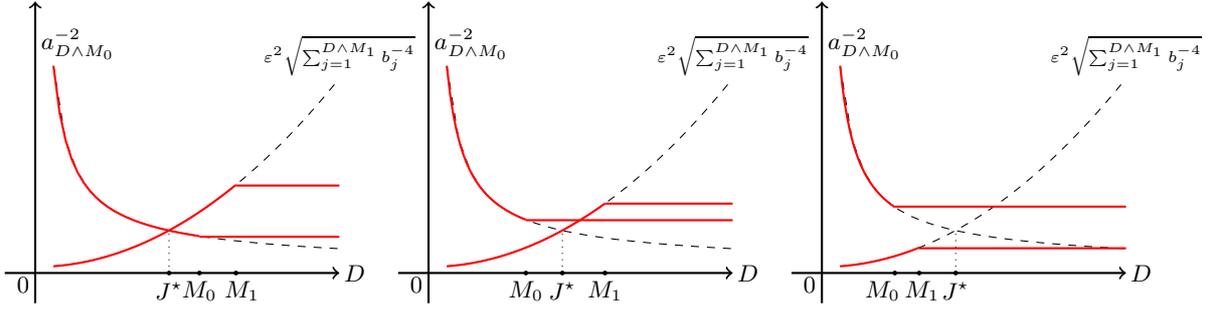
\begin{figure}
\begin{center}
\hspace{0cm}
\begin{tikzpicture}[scale=0.8]
\draw [thick,->] (-0.5,0) -- (5,0);
\draw (5.25,0) node {\footnotesize $D$};
\draw [thick, ->] (0,-0.5) -- (0,4.5);
\draw (-0.2,-0.2) node {\footnotesize $0$};
\draw [black,domain=0.3:5,thin,dashed] plot (\x, {(\x/ln(2.5+\x))^(-1)});
\draw [black,domain=0.3:5,thin,dashed] plot (\x, {0.1+\x^2/8});
\draw (5,3.6) node { \tiny $\varepsilon^2 \sqrt{\sum_{j=1}^{D \wedge M_1} b_j^{-4}}$ \normalsize};
\draw (0.7,3.8) node {\footnotesize $a_{D \wedge M_0}^{-2}$};
\draw [thin,dotted] (2.2,0) -- (2.2,0.73);
\draw[fill] (2.2,0) circle (0.03cm);
\draw (2.2,-0.3) node {\footnotesize $J^\star$};
\draw[fill] (2.7,0) circle (0.03cm);
\draw (2.7,-0.3) node {\footnotesize $M_0$};
\draw [red,thick] (2.7,0.6) -- (5,0.6);
\draw [red,domain=0.3:2.7,thick] plot (\x, {(\x/ln(2.5+\x))^(-1)});
\draw[fill] (3.3,0) circle (0.03cm);
\draw (3.4,-0.3) node {\footnotesize $M_1$};
\draw [red,thick] (3.3,1.45) -- (5,1.45);
\draw [red,domain=0.3:3.3,thick] plot (\x, {0.1+\x^2/8});
\end{tikzpicture} 
\hspace{-0.6cm} \hfill
\begin{tikzpicture}[scale=0.8]
\draw [thick,->] (-0.5,0) -- (5,0);
\draw (5.25,0) node {\footnotesize $D$};
\draw [thick, ->] (0,-0.5) -- (0,4.5);
\draw (-0.2,-0.2) node {\footnotesize $0$};
\draw [black,domain=0.3:5,thin,dashed] plot (\x, {(\x/ln(2.5+\x))^(-1)});
\draw [black,domain=0.3:5,thin,dashed] plot (\x, {0.1+\x^2/8});
\draw (5,3.6) node {\tiny $\varepsilon^2 \sqrt{\sum_{j=1}^{D \wedge M_1} b_j^{-4}}$};
\draw (0.7,3.8) node {\footnotesize $a_{D \wedge M_0}^{-2}$};
\draw [thin,dotted] (2.2,0) -- (2.2,0.73);
\draw[fill] (2.2,0) circle (0.03cm);
\draw (2.2,-0.3) node {\footnotesize $J^\star$};
\draw[fill] (1.6,0) circle (0.03cm);
\draw (1.6,-0.3) node {\footnotesize $M_0$};
\draw [red,thick] (1.6,0.875) -- (5,0.875);
\draw [red,domain=0.3:1.6,thick] plot (\x, {(\x/ln(2.5+\x))^(-1)});
\draw[fill] (2.9,0) circle (0.03cm);
\draw (2.9,-0.3) node {\footnotesize $M_1$};
\draw [red,thick] (2.9,1.15) -- (5,1.15);
\draw [red,domain=0.3:2.9,thick] plot (\x, {0.1+\x^2/8});
\end{tikzpicture}
\hspace{-0.6cm}\hfill
\vspace{0.5cm}
\begin{tikzpicture}[scale=0.8]
\draw [thick,->] (-0.5,0) -- (5,0);
\draw (5.25,0) node {\footnotesize $D$};
\draw [thick, ->] (0,-0.5) -- (0,4.5);
\draw (-0.2,-0.2) node {\footnotesize $0$};
\draw [black,domain=0.3:5,thin,dashed] plot (\x, {(\x/ln(2.5+\x))^(-1)});
\draw [black,domain=0.3:5,thin,dashed] plot (\x, {0.1+\x^2/8});
\draw (5,3.6) node {\tiny $\varepsilon^2 \sqrt{\sum_{j=1}^{D \wedge M_1} b_j^{-4}}$};
\draw (0.7,3.8) node {\footnotesize $a_{D \wedge M_0}^{-2}$};
\draw [thin,dotted] (2.2,0) -- (2.2,0.73);
\draw[fill] (2.2,0) circle (0.03cm);
\draw (2.2,-0.3) node {\footnotesize $J^\star$};
\draw[fill] (1.6,0) circle (0.03cm);
\draw (1.65,-0.3) node {\footnotesize $M_1$};
\draw [red,thick] (1.2,1.1) -- (5,1.1);
\draw [red,domain=0.3:1.2,thick] plot (\x, {(\x/ln(2.5+\x))^(-1)});
\draw[fill] (1.2,0) circle (0.03cm);
\draw (1.0,-0.3) node {\footnotesize $M_0$};
\draw [red,thick] (1.6,0.41) -- (5,0.41);
\draw [red,domain=0.3:1.6,thick] plot (\x, {0.1+\x^2/8});
\end{tikzpicture}
\end{center}
\caption{\textit{{\rm [}Case I: $a_{D\wedge M_0}^{-2} \gtrsim \sigma^2 \ln^{3/2}(1/\sigma)${\rm]} An illustration of the two resulting two terms (red color), namely $\varepsilon^2 \sqrt{\sum_{j=1}^{D \wedge M_1} b_j^{-4}}$ and $a_{D\wedge M_0}^{-2}$, for each $D \in \mathbb{N}$, involved in the (upper bound of the) minimax separation radius $\tilde r^2_{\varepsilon, \sigma}$ (see (\ref{eq:min-sep-rate-cor}), where the bandwidths $M_0$ and $M_1$ are defined in (\ref{eqn:M0-M1-M}). The bandwidth $J^\star$ corresponds to the value $J \in \mathbb{N}$ where the two dashed lines cross, i.e., $J \in \mathbb{N}:\,\varepsilon^2 \sqrt{\sum_{j=1}^{J} b_j^{-4}} = a_{J}^{-2}$. The computation of the separation radius $\tilde{r}_{\varepsilon,\sigma}$, for $0<\varepsilon \leq \varepsilon_0$, $\varepsilon_0 \in ]0,1[$ and $0<\sigma \leq \sigma_0$, $\sigma_0 \in ]0,1[$,  leads to three different scenarios: $J^\star \lesssim M_0 \lesssim M_1$ (left figure), $M_0 \lesssim J^\star \lesssim M_1$ (center figure) and $M_0 \lesssim M_1 \lesssim J^\star$ (right figure).}}
\label{fig:Fig-2}
\end{figure}
\normalsize

\begin{figure}
\begin{center}
\hspace{0cm}
\begin{tikzpicture}[scale=0.8]
\draw [thick,->] (-0.5,0) -- (5,0);
\draw (5.25,0) node {\footnotesize $D$};
\draw [thick, ->] (0,-0.5) -- (0,4.5);
\draw (-0.2,-0.2) node {\footnotesize $0$};
\draw [black,domain=0.3:5,thin,dashed] plot (\x, {(\x/ln(2.5+\x))^(-1)});
\draw [black,domain=0.3:5,thin,dashed] plot (\x, {0.1+\x^2/8});
\draw (5,3.6) node { \tiny $\varepsilon^2 \sqrt{\sum_{j=1}^{D \wedge M_1} b_j^{-4}}$ \normalsize};
\draw (0.7,3.8) node {\footnotesize $a_{D\wedge M^\star}^{-2}$};
\draw [thin,dotted] (2.2,0) -- (2.2,0.73);
\draw[fill] (2.2,0) circle (0.03cm);
\draw (2.2,-0.3) node {\footnotesize $J^\star$};
\draw[fill] (2.7,0) circle (0.03cm);
\draw (2.7,-0.3) node {\footnotesize $M^\star$};
\draw [red,thick] (2.7,0.6) -- (5,0.6);
\draw [red,domain=0.3:2.7,thick] plot (\x, {(\x/ln(2.5+\x))^(-1)});
\draw[fill] (3.3,0) circle (0.03cm);
\draw (3.4,-0.3) node {\footnotesize $M_1$};
\draw [red,thick] (3.3,1.45) -- (5,1.45);
\draw [red,domain=0.3:3.3,thick] plot (\x, {0.1+\x^2/8});
\end{tikzpicture} 
\hspace{-0.6cm} \hfill
\begin{tikzpicture}[scale=0.8]
\draw [thick,->] (-0.5,0) -- (5,0);
\draw (5.25,0) node {\footnotesize $D$};
\draw [thick, ->] (0,-0.5) -- (0,4.5);
\draw (-0.2,-0.2) node {\footnotesize $0$};
\draw [black,domain=0.3:5,thin,dashed] plot (\x, {(\x/ln(2.5+\x))^(-1)});
\draw [black,domain=0.3:5,thin,dashed] plot (\x, {0.1+\x^2/8});
\draw (5,3.6) node {\tiny $\varepsilon^2 \sqrt{\sum_{j=1}^{D \wedge M_1} b_j^{-4}}$};
\draw (0.7,3.8) node {\footnotesize $a_{D\wedge M^\star}^{-2}$};
\draw [thin,dotted] (2.2,0) -- (2.2,0.73);
\draw[fill] (2.2,0) circle (0.03cm);
\draw (2.2,-0.3) node {\footnotesize $J^\star$};
\draw[fill] (1.6,0) circle (0.03cm);
\draw (1.6,-0.3) node {\footnotesize $M^\star$};
\draw [red,thick] (1.6,0.875) -- (5,0.875);
\draw [red,domain=0.3:1.6,thick] plot (\x, {(\x/ln(2.5+\x))^(-1)});
\draw[fill] (2.9,0) circle (0.03cm);
\draw (2.9,-0.3) node {\footnotesize $M_1$};
\draw [red,thick] (2.9,1.15) -- (5,1.15);
\draw [red,domain=0.3:2.9,thick] plot (\x, {0.1+\x^2/8});
\end{tikzpicture}
\hspace{-0.6cm}\hfill
\vspace{0.5cm}
\begin{tikzpicture}[scale=0.8]
\draw [thick,->] (-0.5,0) -- (5,0);
\draw (5.25,0) node {\footnotesize $D$};
\draw [thick, ->] (0,-0.5) -- (0,4.5);
\draw (-0.2,-0.2) node {\footnotesize $0$};
\draw [black,domain=0.3:5,thin,dashed] plot (\x, {(\x/ln(2.5+\x))^(-1)});
\draw [black,domain=0.3:5,thin,dashed] plot (\x, {0.1+\x^2/8});
\draw (5,3.6) node {\tiny $\varepsilon^2 \sqrt{\sum_{j=1}^{D \wedge M_1} b_j^{-4}}$};
\draw (0.7,3.8) node {\footnotesize $a_{D\wedge M^\star}^{-2}$};
\draw [thin,dotted] (2.2,0) -- (2.2,0.73);
\draw[fill] (2.2,0) circle (0.03cm);
\draw (2.2,-0.3) node {\footnotesize $J^\star$};
\draw[fill] (1.6,0) circle (0.03cm);
\draw (1.65,-0.3) node {\footnotesize $M_1$};
\draw [red,thick] (1.2,1.1) -- (5,1.1);
\draw [red,domain=0.3:1.2,thick] plot (\x, {(\x/ln(2.5+\x))^(-1)});
\draw[fill] (1.2,0) circle (0.03cm);
\draw (1.0,-0.3) node {\footnotesize $M^\star$};
\draw [red,thick] (1.6,0.41) -- (5,0.41);
\draw [red,domain=0.3:1.6,thick] plot (\x, {0.1+\x^2/8});
\end{tikzpicture}
\end{center}
\caption{\textit{{\rm [}Case II: $a_{D \wedge M_0}^{-2} \sim \sigma^2 \ln^{3/2}(1/\sigma)${\rm]} An illustration of the two resulting two terms (red color), namely $\varepsilon^2 \sqrt{\sum_{j=1}^{D \wedge M_1} b_j^{-4}}$ and $a_{D\wedge M^\star}^{-2}$, for each $D \in \mathbb{N}$, involved in the (upper bound of the) minimax separation radius $\tilde r^2_{\varepsilon, \sigma}$ (see (\ref{eq:min-sep-rate-cor}), where the bandwidth $M_1$ is defined in (\ref{eqn:M0-M1-M}) and 
the bandwidth $M^\star$ is the value of $M \in \mathbb{N}$ such that the two terms $a_{M}^{-2}$ and $\sigma^2 \ln^{3/2}(1/\sigma)$ are of the same order, i.e., $M \in \mathbb{N}:\,a_{M}^{-2} \sim\sigma^2 \ln^{3/2}(1/\sigma)$.
The bandwidth $J^\star$ corresponds to the value $J \in \mathbb{N}$ where the two dashed lines cross, i.e., $J \in \mathbb{N}:\,\varepsilon^2 \sqrt{\sum_{j=1}^{J} b_j^{-4}} = a_{J}^{-2}$. 
The computation of the separation radius $\tilde{r}_{\varepsilon,\sigma}$, for $0<\varepsilon \leq \varepsilon_0$, $\varepsilon_0 \in ]0,1[$ and $0<\sigma \leq \sigma_0$, $\sigma_0 \in ]0,1[$,  leads to three different scenarios: $J^\star \lesssim M^\star \lesssim M_1$ (left figure), $M^\star \lesssim J^\star \lesssim M_1$ (center figure) and $M^\star \lesssim M_1 \lesssim J^\star$ (right figure).}}
\label{fig:Fig-3}
\end{figure}
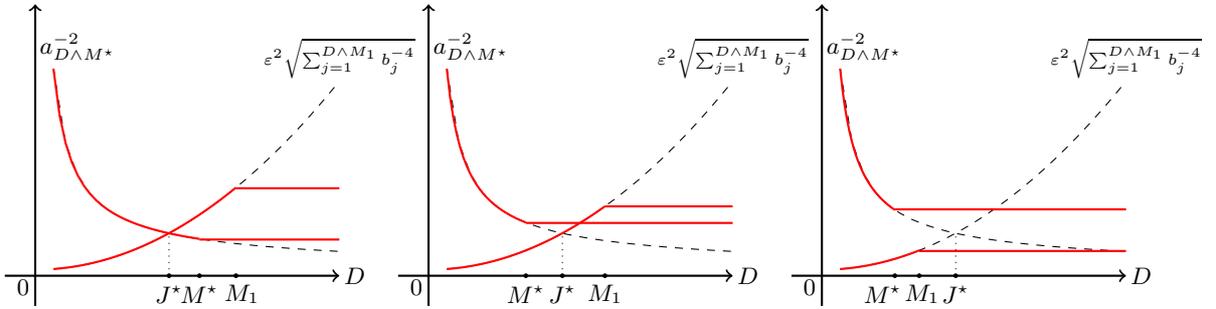
\normalsize

\subsection{Upper Bounds: Specific Cases}
\label{s:p_rates}

For the sake of convenience, we give the proof of each item (i)-(iv) in Theorem \ref{thm:asupper} in different sections. 

\subsubsection{Case (i): Mildly ill-posed problems with ordinary smooth functions}
\label{casei}
Recall that
\begin{equation}
\label{eq:a-b-mildly}
b_j \sim j^{-t}, \; t>0,  \quad \text{and} \quad a_j \sim j^s, \; s>0, \quad j \in \mathbb{N}.
\end{equation}

\begin{proposition}
\label{eq:prop-a-b-mildly}
Assume that the sequences $b=(b_j)_{j \in \mathbb{N}}$ and $a=(a_j)_{j \in \mathbb{N}}$ are given by (\ref{eq:a-b-mildly}). Then, there exists $\varepsilon_0, \sigma_0 \in ]0,1[$ such that, for all $0 < \varepsilon \leq \varepsilon_0$ and $0 <\sigma \leq \sigma_0$, the minimax separation radius $\tilde r_{\varepsilon, \sigma}$ satisfies
\begin{equation}
\label{eq:prop-mildly-upper}
\tilde r^2_{\varepsilon, \sigma} \lesssim \varepsilon^{\frac{4s}{2s+2t+1/2}}
\vee \left[\sigma \ln^{3/4}(1/\sigma) \right]^{2 \left(\frac{s}{t} \wedge 1\right)}.
\end{equation}
\end{proposition}

\noindent 
\textsc{Proof of Proposition \ref{eq:prop-a-b-mildly}} In a first time, we determine the order of the bandwidths $M_0$ and $M_1$. Setting 
$$ \bar M_1 :=  \left(\sigma \sqrt{ \frac{1}{t} \ln \left( 1/\sigma \right)}\right)^{-1/t} \quad \mathrm{and} \quad \tilde M_1 :=  \left(\sigma \sqrt{ \frac{1}{2t} \ln \left( 1/\sigma\right)}\right)^{-1/t} $$
we get 
$$ \sigma h_{1,\bar M_1} \sim \sigma \sqrt{ \ln (\bar M_1) } = \sigma \sqrt{ \frac{1}{t}\ln\left( 1/\sigma\right) - \frac{1}{2t} \ln \left( \frac{1}{t} \ln\left(1/\sigma\right) \right)} \leq \sigma \sqrt{ \frac{1}{t}\ln\left( 1/\sigma\right) } \sim b_{\bar M_1},   $$
which implies that $M_1 \gtrsim \bar M_1$. At the same time
\begin{eqnarray*} 
\sigma h_{1,\tilde M_1} \sim \sigma \sqrt{ \ln (\tilde M_1) } 
& = & \sigma \sqrt{ \frac{1}{t}\ln\left( 1/\sigma\right) - \frac{1}{2t} \ln \left( \frac{1}{2t} \ln\left(1/\sigma\right) \right)},\\
& = & \sigma \sqrt{ \frac{1}{2t}\ln\left( 1/\sigma\right) +\frac{1}{2t}\ln\left( 1/\sigma\right) - \frac{1}{2t} \ln \left( \frac{1}{2t} \ln\left(1/\sigma\right) \right)},\\
& \gtrsim & \sigma \sqrt{ \frac{1}{2t}\ln\left( 1/\sigma\right) } \sim b_{\bar M_1},   
\end{eqnarray*}
which implies that $M_1 \lesssim \tilde M_1$. Hence, we can conclude that
$$ M_1\sim \left(\sigma \sqrt{ \ln \left( 1/\sigma \right)}\right)^{-1/t} .$$
Similarly, we get that 
$$ M_0\sim \left(\sigma \sqrt{  \ln \left( 1/\sigma \right)}\right)^{-1/t} .$$

In order to control the terms involved in the upper bound on the minimax separation radius,  we consider the cases $s < t$ and $s\geq t$ separately.\\

Consider first the case $s < t$. In this case, for all $D \in \mathbb{N}$, 
\begin{eqnarray*}
a_{D \wedge M_0}^{-2} \gtrsim a_{M_0}^{-2} \sim M_0^{-2s} \sim
\left(\sigma \ln^{1/2}(1/\sigma)\right)^{2s/t} \gtrsim \sigma^{2} \ln^{3/2}(1/\sigma).
\end{eqnarray*}
Hence,
\begin{eqnarray*}
\tilde r^2_{\varepsilon, \sigma} \lesssim \inf_{D \in \mathbb{N}}\left[
\varepsilon^2 \sqrt{\sum_{j=1}^{D \wedge M_1}  b_j^{-4}} +  \left[\sigma^2 \ln^{3/2}(1/\sigma) \vee a^{-2}_{D \wedge M_0}\right] \right] 
\lesssim
\inf_{D \in \mathbb{N}}\left[
\varepsilon^2 \sqrt{\sum_{j=1}^{D \wedge M_1}  b_j^{-4}} +  a^{-2}_{D \wedge M_0}\right].  \nonumber\\
\end{eqnarray*}
Define now the value of $J \in \mathbb{N}$ that satisfies the following equation
\begin{eqnarray*}
\varepsilon^2 \sqrt{\sum_{j=1}^{J}  b_j^{-4}} \sim a^{-2}_{J} 
\; \Leftrightarrow \; \varepsilon^2 J^{2t+1/2} \sim J^{2s}
\; \Leftrightarrow \; J:=J^\star \sim \varepsilon^{\frac{-2}{2s+2t+1/2}}.
\end{eqnarray*}
We now consider the following situations (see Figure \ref{fig:Fig-2} for a graphical illustration):
\begin{itemize}
\item ($J^\star \lesssim M_0$) In this case, 
$$
\tilde r^2_{\varepsilon, \sigma} \lesssim a_{J^{\star}}^{-2} \lesssim 
\varepsilon^{\frac{4s}{2s+2t+1/2}}.
$$
\item ($J^\star \gtrsim M_1$) In this case, 
$$
\tilde r^2_{\varepsilon, \sigma} \lesssim 
\inf_{D \in \mathbb{N}}\left[a^{-2}_{D \wedge M_0}\right]
\lesssim a^{-2}_{M_0} \sim 
 \left[\sigma \ln^{1/2}(1/\sigma) \right]^{\frac{2 s}{t}}.
$$
\item ($M_0 \lesssim J^{\star} \lesssim M_1$) In this case,
\begin{eqnarray}
\tilde r^2_{\varepsilon, \sigma} &\lesssim&
\inf_{D \in \mathbb{N}}\left[
\varepsilon^2 \sqrt{\sum_{j=1}^{D \wedge M_1}  b_j^{-4}} +  a^{-2}_{D \wedge M_0}\right].  \nonumber\\
&=& \left\{\inf_{D \leq M_0} \left[
\varepsilon^2 \sqrt{\sum_{j=1}^{D \wedge M_1}  b_j^{-4}} +  a^{-2}_{D \wedge M_0}\right] \right\} \wedge \left\{\inf_{D > M_0} \left[
\varepsilon^2 \sqrt{\sum_{j=1}^{D \wedge M_1}  b_j^{-4}} +  a^{-2}_{D \wedge M_0}\right] \right\} \nonumber \\
&\lesssim& a_{M_0}^{-2} \wedge \left\{\inf_{D > M_0} \left[
\varepsilon^2 \sqrt{\sum_{j=1}^{D \wedge M_1}  b_j^{-4}} +  a^{-2}_{M_0}\right] \right\} \nonumber \\
&\lesssim& a_{M_0}^{-2} \wedge \left\{
\varepsilon^2 \sqrt{\sum_{j=1}^{M_0}  b_j^{-4}} +  a^{-2}_{M_0} \right\} \nonumber \\
&\lesssim& a_{M_0}^{-2} \sim \left[\sigma \ln^{1/2}(1/\sigma) \right]^{\frac{2 s}{t}}. \nonumber
\end{eqnarray}
\end{itemize}
Combining the above terms, we immediately get
\begin{equation}
\label{eq:M-F-01}
\tilde r^2_{\varepsilon, \sigma} \lesssim \varepsilon^{\frac{4s}{2s+2t+1/2}}
\vee \left[\sigma \ln^{1/2}(1/\sigma) \right]^{\frac{2s}{t}}.
\end{equation}

Consider now the case $s\geq t$. Define the value of $M \in \mathbb{N}$ that satisfies the following equation
$$
a_{M}^{-2} \sim \sigma^2 \ln^{3/2}(1/\sigma) \; \Leftrightarrow \; 
M:=M^{\star} \sim \left[ \sigma \ln^{3/4}(1/\sigma) \right]^{-\frac{1}{s}}.
$$
Hence,
\begin{eqnarray*}
\tilde r^2_{\varepsilon, \sigma} \lesssim \inf_{D \in \mathbb{N}}\left[
\varepsilon^2 \sqrt{\sum_{j=1}^{D \wedge M_1}  b_j^{-4}} +  \left[\sigma^2 \ln^{3/2}(1/\sigma) \vee a^{-2}_{D \wedge M_0}\right] \right] 
\lesssim
\inf_{D \in \mathbb{N}}\left[
\varepsilon^2 \sqrt{\sum_{j=1}^{D \wedge M_1}  b_j^{-4}} +  a^{-2}_{D\wedge M^{\star}}\right].  \nonumber\\
\end{eqnarray*}
Working along the lines of the case $s \leq t$, by replacing $M_0$ by $M^\star$ (see Figure \ref{fig:Fig-3}), we get
\begin{equation}
\label{eq:M-F-3}
\tilde r^2_{\varepsilon, \sigma} \lesssim \varepsilon^{\frac{4s}{2s+2t+1/2}}
\vee \left[\sigma \ln^{3/4}(1/\sigma) \right]^{2}.
\end{equation}

Hence, (\ref{eq:prop-mildly-upper}) follows thanks to (\ref{eq:M-F-01}) and (\ref{eq:M-F-3}). This completes the proof of the proposition. \begin{flushright} $\Box$ \end{flushright}

\subsubsection{ Case (ii): Mildly ill-posed problems with super smooth functions}
Recall that
\begin{equation}
\label{eq:a-b-super-mildly}
b_j \sim j^{-t}, \; t>0,  \quad \text{and} \quad a_j \sim \exp\{js\}, \; s>0, \quad j \in \mathbb{N}.
\end{equation}

\begin{proposition}
\label{eq:prop-a-b-super-mildly}
Assume that the sequences $b=(b_j)_{j \in \mathbb{N}}$ and $a=(a_j)_{j \in \mathbb{N}}$ are given by (\ref{eq:a-b-super-mildly}). Then, there exists $\varepsilon_0, \sigma_0 \in ]0,1[$ such that, for all $0 < \varepsilon \leq \varepsilon_0$ and $0 <\sigma \leq \sigma_0$, the minimax separation radius $\tilde r_{\varepsilon, \sigma}$ satisfies
\begin{equation}
\label{eq:prop-super-mildly-upper}
\tilde r^2_{\varepsilon, \sigma} \lesssim \varepsilon^2 \left[\ln\left(1/\varepsilon\right)\right]^{\left(2t+\frac{1}{2}\right)} \vee \sigma^{2} \ln^{\frac{3}{2}}\left(1/\sigma\right).
\end{equation}
\end{proposition}

\noindent 
\textsc{Proof of Proposition \ref{eq:prop-a-b-super-mildly}} According to Section \ref{casei}, we obtain again
$$ M_1\sim \left(\sigma \sqrt{ \ln \left( 1/\sigma \right)}\right)^{-1/t} \quad \mathrm{and} \quad M_0\sim \left(\sigma \sqrt{  \ln \left( 1/\sigma \right)}\right)^{-1/t} .$$
Then, for all $D \in \mathbb{N}$, 
\begin{eqnarray*}
a_{D \wedge M_0}^{-2} \gtrsim a_{M_0}^{-2} \sim \exp\{-2M_0 s\} \sim
\exp\left\{-2s \left(\sigma \ln^{1/2}(1/\sigma) \right)^{-1/t}\right\} \lesssim \sigma^{2} \ln^{3/2}(1/\sigma).
\end{eqnarray*}
Define as in the previous case the value $M \in \mathbb{N}$ that satisfies the following equation
$$
a_{M}^{-2} \sim \sigma^2 \ln^{3/2}(1/\sigma) \; \Leftrightarrow \; 
M=:M^\star \sim \frac{1}{s}  \ln \left[\frac{1}{\sigma \ln^{3/4}(1/\sigma)} \right].
$$
Hence,
\begin{eqnarray*}
\tilde r^2_{\varepsilon, \sigma} \lesssim \inf_{D \in \mathbb{N}}\left[
\varepsilon^2 \sqrt{\sum_{j=1}^{D \wedge M_1}  b_j^{-4}} +  \left[\sigma^2 \ln^{3/2}(1/\sigma) \vee a^{-2}_{D \wedge M_0}\right] \right] 
\lesssim
\inf_{D \in \mathbb{N}}\left[
\varepsilon^2 \sqrt{\sum_{j=1}^{D \wedge M_1}  b_j^{-4}} +  a^{-2}_{D\wedge M^{\star}}\right].  \nonumber\\
\end{eqnarray*}
Define now the value of $J \in \mathbb{N}$ that satisfies the following equation
\begin{eqnarray*}
\varepsilon^2 \sqrt{\sum_{j=1}^{J}  b_j^{-4}} \sim a^{-2}_{J} 
\; \Leftrightarrow \; \varepsilon^2 J^{2t+1/2} \sim \exp\{-2Js\}
\; \Leftrightarrow \; J:=J^\star \sim \frac{1}{\sigma}\ln \left(1/\varepsilon\right) - \ln \left[ \left(\frac{1}{\sigma}\ln \left(1/\varepsilon\right)\right)^{2t+\frac{1}{2}}\right].
\end{eqnarray*}

We now consider the following situations:
\begin{itemize}
\item ($J^\star \lesssim M^\star$) In this case, 
$$
\tilde r^2_{\varepsilon, \sigma} \lesssim a_{J^{\star}}^{-2} \lesssim 
\varepsilon^2 \left[\ln\left(1/\varepsilon\right)\right]^{\left(2t+\frac{1}{2}\right)}.
$$
\item ($J^\star \gtrsim M_1$) In this case, 
$$
\tilde r^2_{\varepsilon, \sigma} \lesssim 
\inf_{D \in \mathbb{N}}\left[a^{-2}_{D \wedge M^\star}\right]
\lesssim a^{-2}_{M^\star} \sim 
\sigma^{2} \ln^{\frac{3}{2}}\left(1/\sigma\right).
$$
\item ($M^\star \lesssim J^{\star} \lesssim M_1$) In this case,
\begin{eqnarray*}
\tilde r^2_{\varepsilon, \sigma} &\lesssim&
\inf_{D \in \mathbb{N}}\left[
\varepsilon^2 \sqrt{\sum_{j=1}^{D \wedge M_1}  b_j^{-4}} +  a^{-2}_{D \wedge M^\star}\right].  \nonumber\\
&=& \left\{\inf_{D \leq M^\star} \left[
\varepsilon^2 \sqrt{\sum_{j=1}^{D \wedge M_1}  b_j^{-4}} +  a^{-2}_{D \wedge M^\star}\right] \right\} \wedge \left\{\inf_{D > M^\star} \left[
\varepsilon^2 \sqrt{\sum_{j=1}^{D \wedge M_1}  b_j^{-4}} +  a^{-2}_{D \wedge M^\star}\right] \right\} \nonumber \\
&\lesssim& a_{M^\star}^{-2} \wedge \left\{\inf_{D > M^\star} \left[
\varepsilon^2 \sqrt{\sum_{j=1}^{D \wedge M_1}  b_j^{-4}} +  a^{-2}_{M^\star}\right] \right\} \nonumber \\
&\lesssim& a_{M^\star}^{-2} \wedge \left\{
\varepsilon^2 \sqrt{\sum_{j=1}^{M^\star}  b_j^{-4}} +  a^{-2}_{M^\star} \right\} \nonumber \\
&\lesssim& a_{M^\star}^{-2} \sim \sigma^{2} \ln^{\frac{3}{2}}\left(1/\sigma\right). \nonumber
\end{eqnarray*}
\end{itemize}
Combining the above terms, we immediately get (\ref{eq:prop-super-mildly-upper}). This completes the proof of the proposition. \begin{flushright} $\Box$ \end{flushright}

\subsubsection{Case (iii): Severely ill-posed problems with ordinary smooth functions}
\label{s:caseiii}
Recall that
\begin{equation}
\label{eq:a-b-mildly-super}
b_j \sim \exp\{-jt\}, \; t>0,  \quad \text{and} \quad a_j \sim j^s, \; s>0, \quad j \in \mathbb{N}.
\end{equation}

\begin{proposition}
\label{eq:prop-a-b-mildly-super}
Assume that the sequences $b=(b_j)_{j \in \mathbb{N}}$ and $a=(a_j)_{j \in \mathbb{N}}$ are given by (\ref{eq:a-b-mildly-super}). Then, there exists $\varepsilon_0, \sigma_0 \in ]0,1[$ such that, for all $0 < \varepsilon \leq \varepsilon_0$ and $0 <\sigma \leq \sigma_0$, the minimax separation radius $\tilde r_{\varepsilon, \sigma}$ satisfies
\begin{equation}
\label{eq:prop-mildly-super-upper}
\tilde r^2_{\varepsilon, \sigma} \lesssim \left[\ln\left(1/\varepsilon\right)\right]^{-2s} \vee \left[\ln\left(\frac{1}{\sigma \ln^{1/2}(1/\sigma)}\right)\right]^{-2s}.
\end{equation}
\end{proposition}

\noindent 
\textsc{Proof of Proposition \ref{eq:prop-a-b-mildly-super}} In a first time, we determine the order of the bandwidths $M_0$ and $M_1$. Setting 
$$ \bar M_1 :=  \frac{1}{t} \ln\left(\frac{1}{\sigma \ln^{1/2}(1/\sigma)}\right)  \quad \mathrm{and} \quad \tilde M_1 :=  \frac{1}{t} \ln\left(1/\sigma\right), $$
we get 
$$ \sigma h_{1,\bar M_1} \sim \sigma \sqrt{ \ln (\bar M_1) } = \sigma \sqrt{ \ln\left( \frac{1}{t} \ln\left(\frac{1}{\sigma \ln^{1/2}(1/\sigma)}\right)\right)} \lesssim  e^{-\bar M_1 t} = \sigma \sqrt{\ln(1/\sigma)} \sim b_{\bar M_1},   $$
which implies that $M_1 \geq \bar M_1$ for $\sigma$ small enough. At the same time
\begin{eqnarray*} 
\sigma h_{1,\tilde M_1} \sim \sigma \sqrt{ \ln (\tilde M_1) } 
& = &\sigma \sqrt{ \ln \left( \frac{1}{t} \ln\left(1/\sigma\right)\right) }  \\
& \gtrsim &b_{\tilde M_1} \sim e^{-\tilde M_1 t} =\sigma,   
\end{eqnarray*}
which implies that $M_1 \leq \tilde M_1$ for $\sigma$ small enough. Hence, we can conclude that
$$   \frac{1}{t} \ln\left(\frac{1}{\sigma \ln^{1/2}(1/\sigma)}\right)   \leq M_1\leq \frac{1}{t} \ln\left(1/\sigma\right) ,$$
for $\sigma$ small enough. Similarly, we get that 
$$   \frac{1}{t} \ln\left(\frac{1}{\sigma \ln^{1/2}(1/\sigma)}\right)   \leq M_0 \leq \frac{1}{t} \ln\left(1/\sigma\right) ,$$
for $\sigma$ small enough. \\

Now, we turn our attention to the proof of (\ref{eq:prop-mildly-super-upper}). For all $D \in \mathbb{N}$, 
\begin{eqnarray*}
a_{D \wedge M_0}^{-2} \gtrsim a_{M_0}^{-2} \gtrsim M_0^{-2 s} \sim
\left[\ln\left(\frac{1}{\sigma \ln^{1/2}(1/\sigma)}\right)\right]^{-2s} \gtrsim \sigma^{2} \ln^{3/2}(1/\sigma).
\end{eqnarray*}
Hence,
\begin{eqnarray*}
\tilde r^2_{\varepsilon, \sigma} \lesssim \inf_{D \in \mathbb{N}}\left[
\varepsilon^2 \sqrt{\sum_{j=1}^{D \wedge M_1}  b_j^{-4}} +  \left[\sigma^2 \ln^{3/2}(1/\sigma) \vee a^{-2}_{D \wedge M_0}\right] \right] 
\lesssim
\inf_{D \in \mathbb{N}}\left[
\varepsilon^2 \sqrt{\sum_{j=1}^{D \wedge M_1}  b_j^{-4}} +  a^{-2}_{D \wedge M_0}\right].  \nonumber\\
\end{eqnarray*}
Define now the value of $J \in \mathbb{N}$ that satisfies the following equation
\begin{eqnarray*}
\varepsilon^2 \sqrt{\sum_{j=1}^{J}  b_j^{-4}} \sim a^{-2}_{J} 
\; \Leftrightarrow \; \varepsilon^2 \exp\{2tJ\} \sim J^{-2s}
\; \Leftrightarrow \; J:=J^\star \sim \frac{1}{t}\ln \left(1/\varepsilon\right) - \ln \left[ \left(\frac{1}{t}\ln \left(1/\varepsilon\right)\right)^{2s}\right].
\end{eqnarray*}

We now consider the following situations:
\begin{itemize}
\item ($J^\star \lesssim M_0$) In this case, 
$$
\tilde r^2_{\varepsilon, \sigma} \lesssim a_{J^{\star}}^{-2} \lesssim 
\left[\ln\left(1/\varepsilon\right)\right]^{-2s}.
$$
\item ($J^\star \gtrsim M_1$) In this case, 
$$
\tilde r^2_{\varepsilon, \sigma} \lesssim 
\inf_{D \in \mathbb{N}}\left[a^{-2}_{D \wedge M_0}\right]
\lesssim a^{-2}_{M_0} \sim 
\left[\ln\left(\frac{1}{\sigma \ln^{1/2}(1/\sigma)}\right)\right]^{-2s}.
$$
\item ($M_0 \lesssim J^{\star} \lesssim M_1$) In this case,
\begin{eqnarray*}
\tilde r^2_{\varepsilon, \sigma} &\lesssim&
\inf_{D \in \mathbb{N}}\left[
\varepsilon^2 \sqrt{\sum_{j=1}^{D \wedge M_1}  b_j^{-4}} +  a^{-2}_{D \wedge M_0}\right].  \nonumber\\
&=& \left\{\inf_{D \leq M_0} \left[
\varepsilon^2 \sqrt{\sum_{j=1}^{D \wedge M_1}  b_j^{-4}} +  a^{-2}_{D \wedge M_0}\right] \right\} \wedge \left\{\inf_{D > M_0} \left[
\varepsilon^2 \sqrt{\sum_{j=1}^{D \wedge M_1}  b_j^{-4}} +  a^{-2}_{D \wedge M_0}\right] \right\} \nonumber \\
&\lesssim& a_{M_0}^{-2} \wedge \left\{\inf_{D > M_0} \left[
\varepsilon^2 \sqrt{\sum_{j=1}^{D \wedge M_1}  b_j^{-4}} +  a^{-2}_{M_0}\right] \right\} \nonumber \\
&\lesssim& a_{M_0}^{-2} \wedge \left\{
\varepsilon^2 \sqrt{\sum_{j=1}^{M_0}  b_j^{-4}} +  a^{-2}_{M_0} \right\} \nonumber \\
&\lesssim& a_{M_0}^{-2} \lesssim\left[\ln\left(\frac{1}{\sigma \ln^{1/2}(1/\sigma)}\right)\right]^{-2s}. \nonumber
\end{eqnarray*}
\end{itemize}
Combining the above terms, we immediately get (\ref{eq:prop-mildly-super-upper}). This completes the proof of the proposition. \begin{flushright} $\Box$ \end{flushright}

\subsubsection{Case (iv): Severely ill-posed problems with super smooth functions}
Recall that
\begin{equation}
\label{eq:a-b-super}
b_j \sim\exp\{-jt\}, \; t>0,  \quad \text{and} \quad a_j \sim \exp\{js\}, \; s>0, \quad j \in \mathbb{N}.
\end{equation}

\begin{proposition}
\label{eq:prop-a-b-super}
Assume that the sequences $b=(b_j)_{j \in \mathbb{N}}$ and $a=(a_j)_{j \in \mathbb{N}}$ are given by (\ref{eq:a-b-super}). Then, there exists $\varepsilon_0, \sigma_0 \in ]0,1[$ such that, for all $0 < \varepsilon \leq \varepsilon_0$ and $0 <\sigma \leq \sigma_0$, the minimax separation radius $\tilde r_{\varepsilon, \sigma}$ satisfies
\begin{equation}
\label{eq:prop-super-upper}
\tilde r^2_{\varepsilon, \sigma} \lesssim \varepsilon^{\frac{2s}{s+t}}
\vee [\sigma\ln^{1/2}(1/\sigma)]^{2 \left(\frac{s}{t} \wedge 1\right)}.
\end{equation}
\end{proposition}

\noindent 
\textsc{Proof of Proposition \ref{eq:prop-a-b-super}} According to Section \ref{s:caseiii}, we obtain again that
$$   \frac{1}{t} \ln\left(\frac{1}{\sigma \ln^{1/2}(1/\sigma)}\right)   \leq M_1\leq \frac{1}{t} \ln\left(1/\sigma\right) $$
and
$$   \frac{1}{t} \ln\left(\frac{1}{\sigma \ln^{1/2}(1/\sigma)}\right)   \leq M_0 \leq \frac{1}{t} \ln\left(1/\sigma\right) ,$$
for $\sigma$ small enough. Now, we consider the cases $s < t$ and $s\geq t$ separately.\\

Consider first the case $s < t$. In this case, for all $D \in \mathbb{N}$, 
\begin{eqnarray*}
a_{D \wedge M_0}^{-2} \gtrsim a_{M_0}^{-2} \sim \exp\{-2sM_0\} \gtrsim
\left(\sigma \ln^{1/2}(1/\sigma)\right)^{2s/t} \gtrsim \sigma^{2} \ln^{3/2}(1/\sigma).
\end{eqnarray*}
Hence,
\begin{eqnarray*}
\tilde r^2_{\varepsilon, \sigma} &\lesssim& \inf_{D \in \mathbb{N}}\left[
\varepsilon^2 \sqrt{\sum_{j=1}^{D \wedge M_1}  b_j^{-4}} +  \left[\sigma^2 \ln^{3/2}(1/\sigma) \vee a^{-2}_{D \wedge M_0}\right] \right] \lesssim
\inf_{D \in \mathbb{N}}\left[
\varepsilon^2 \sqrt{\sum_{j=1}^{D \wedge M_1}  b_j^{-4}} +  a^{-2}_{D \wedge M_0}\right].  \nonumber\\
\end{eqnarray*}
Define now the value of $J \in \mathbb{N}$ that satisfies the following equation
\begin{eqnarray*}
\varepsilon^2 \sqrt{\sum_{j=1}^{J}  b_j^{-4}} \sim a^{-2}_{J} 
\; \Leftrightarrow \; \varepsilon^2 \exp\{2tJ\} \sim \exp\{2sJ\}
\; \Leftrightarrow \; J:=J^\star \sim \frac{1}{s+t}\ln \left(1/\varepsilon\right).
\end{eqnarray*}
We now consider the following situations:
\begin{itemize}
\item ($J^\star \lesssim M_0$) In this case, 
$$
\tilde r^2_{\varepsilon, \sigma} \lesssim a_{J^{\star}}^{-2} \lesssim 
\varepsilon^{\frac{2s}{s+t}}.
$$
\item ($J^\star \gtrsim M_1$) In this case, 
$$
\tilde r^2_{\varepsilon, \sigma} \lesssim 
\inf_{D \in \mathbb{N}}\left[a^{-2}_{D \wedge M_0}\right]
\lesssim a^{-2}_{M_0} \lesssim 
 \left[\sigma \ln^{1/2}(1/\sigma) \right]^{\frac{2 s}{t}}.
$$
\item ($M_0 \lesssim J^{\star} \lesssim M_1$) In this case,
\begin{eqnarray}
\tilde r^2_{\varepsilon, \sigma} &\lesssim&
\inf_{D \in \mathbb{N}}\left[
\varepsilon^2 \sqrt{\sum_{j=1}^{D \wedge M_1}  b_j^{-4}} +  a^{-2}_{D \wedge M_0}\right].  \nonumber\\
&=& \left\{\inf_{D \leq M_0} \left[
\varepsilon^2 \sqrt{\sum_{j=1}^{D \wedge M_1}  b_j^{-4}} +  a^{-2}_{D \wedge M_0}\right] \right\} \wedge \left\{\inf_{D > M_0} \left[
\varepsilon^2 \sqrt{\sum_{j=1}^{D \wedge M_1}  b_j^{-4}} +  a^{-2}_{D \wedge M_0}\right] \right\} \nonumber \\
&\lesssim& a_{M_0}^{-2} \wedge \left\{\inf_{D > M_0} \left[
\varepsilon^2 \sqrt{\sum_{j=1}^{D \wedge M_1}  b_j^{-4}} +  a^{-2}_{M_0}\right] \right\} \nonumber \\
&\lesssim& a_{M_0}^{-2} \wedge \left\{
\varepsilon^2 \sqrt{\sum_{j=1}^{M_0}  b_j^{-4}} +  a^{-2}_{M_0} \right\} \nonumber \\
&\lesssim& a_{M_0}^{-2} \lesssim \left[\sigma \ln^{1/2}(1/\sigma) \right]^{\frac{2 s}{t}}. \nonumber
\end{eqnarray}
\end{itemize}
Combining the above terms, we immediately get
\begin{equation}
\label{eq:M-F-01-FF3}
\tilde r^2_{\varepsilon, \sigma} \lesssim \varepsilon^{\frac{2s}{s+t}}
\vee \left[\sigma \ln^{1/2}(1/\sigma) \right]^{\frac{2s}{t}}.
\end{equation}

Consider now the case $s\geq t$. Define the value $M \in \mathbb{N}$ that satisfies the following equation
$$
a_{M}^{-2} \sim \sigma^2 \ln^{3/2}(1/\sigma) \; \Leftrightarrow \; 
M:=M^\star \sim \frac{1}{s}  \ln \left[\frac{1}{\sigma \ln^{3/4}(1/\sigma)} \right].
$$
Hence,
\begin{eqnarray*}
\tilde r^2_{\varepsilon, \sigma} \lesssim \inf_{D \in \mathbb{N}}\left[
\varepsilon^2 \sqrt{\sum_{j=1}^{D \wedge M_1}  b_j^{-4}} +  \left[\sigma^2 \ln^{3/2}(1/\sigma) \vee a^{-2}_{D \wedge M_0}\right] \right] 
\lesssim
\inf_{D \in \mathbb{N}}\left[
\varepsilon^2 \sqrt{\sum_{j=1}^{D \wedge M_1}  b_j^{-4}} +  a^{-2}_{M^{\star}}\right].  \nonumber\\
\end{eqnarray*}
Working along the lines of the case $s < t$ by replacing $M_0$ by $M^\star$, we get
\begin{equation}
\label{eq:M-F-3-FF}
\tilde r^2_{\varepsilon, \sigma} \lesssim \varepsilon^{\frac{2s}{s+t}}
\vee \left[\sigma \ln^{1/2}(1/\sigma) \right]^{2}.
\end{equation}

Hence, (\ref{eq:prop-super-upper}) follows thanks to (\ref{eq:M-F-01-FF3}) and (\ref{eq:M-F-3-FF}). This completes the proof of the proposition. \begin{flushright} $\Box$ \end{flushright}

\subsection{Non-Asymptotic Lower Bounds}
\label{s:lower_bounds}

\subsubsection{Proof of Proposition \ref{prop:1-e=0}} 
\label{s:prlowerb}
Let $\theta_0 \in \mathcal{E}_a$ be given sequence (to be made precise below). Given a (prior) probability measure $\pi$ on the set associated with $H_1$, i.e., a probability measure $\pi$ on $
\tilde \Theta_{a,\theta_0}(r_\sigma,b):=
\Theta_{a,\theta_0}(r_\sigma) \times \mathcal{B}(b)$, where $\Theta_{a,\theta_0}(r_\sigma) =\theta_0 + \Theta_{a}(r_\sigma)$, by standard Bayesian arguments (see, e.g., Section 3.1 of \cite{MS_2014}), we arrive at
\begin{eqnarray}
\label{eqn:lower-bound-2nd-kep}
\beta_{0,\sigma,b}(\Theta(r_\sigma),\mathcal{B}(b))
&=& \inf_{\tilde \Psi_{\alpha}:\,\Alpha_{0,\sigma}(\tilde \Psi_\alpha)\leq \alpha}\; \sup_{\substack{\theta_0 \in \mathcal{E}_a,\; \theta-\theta_0\in \Theta_a(r_\sigma) \\ \bar b\in \mathcal{B}(b)}} \mathbb{P}_{\theta, \bar b}(\tilde \Psi_\alpha=0) \nonumber\\
&\geq&
1-\alpha -\frac{1}{2}(\mathbb{E}_0[L^2_\pi(Y,X)]-1)^{1/2},
\end{eqnarray} 
where $L_\pi(Y,X)$ denotes the likelihood ratio between the two measures $\mathbb{P}_\pi$ and $\mathbb{P}_0$, $\mathbb{E}_0$ denotes the expectation with respect to $\mathbb{P}_0$, with $\mathbb{P}_0=\int_{\tilde \Theta_{a,\theta_0}(r_\sigma,b)} \mathbb{P}_{\theta_0,\bar b} \;d\pi(\theta,\bar b)$ and $\mathbb{P}_\pi = \int_{\tilde \Theta_{a,\theta_0}(r_\sigma,b)} \mathbb{P}_{\theta,\bar b}\; d\pi(\theta,\bar b)$, and the last inequality is obtained by standard calculations (see, e.g., Section 3.1 of \cite{MS_2014}). \\

The probability measure $\pi$ on $\tilde \Theta_{a,\theta_0}(r_\sigma,b)$ is selected as product probability measure, i.e.,
$$
\pi=\prod_{j \in \mathbb{N}} \pi_j, \quad \pi_j=\pi_{j,1}\times\pi_{j,2}, \quad j \in \mathbb{N}.
$$
Then, given the sequence $\theta$ and the bandwidth $D \in \mathbb{N}$ (to be made precise below), we set
$$
\pi_{j,1}=\delta_{\theta_{j,0}} \quad \text{and} \quad \pi_{j,2}=\delta_{b_j}, \quad j \neq D,
$$
and 
$$
\pi_{D,1}=G_D^{-1}(C_0,C_1)  \delta_{\theta_{D}} \quad \text{and} \quad d\pi_{D,2}(t)=\frac{1}{\sigma \sqrt{2\pi}} \exp\left\{-\frac{1}{2\sigma^2}(t-b_D)^2\right\}dt,
$$
where $G_D(C_0,C_1)=(1/\sigma \sqrt{2\pi}) \int_{C_0b_D}^{C_1b_D}\exp\left\{-(t-b_D)^2/(2\sigma^2) \right\}dt$. In some sense, using the above product probability measure $\pi$, we deal with observations $(Y,X)=(Y_j,X_j)_{j\in\mathbb{N}}$ from the following Bayesian sequence model
$$
Y_j=b_j\theta_{j,0}, \quad X_j= b_j+\sigma \eta_j,\quad j \in \mathbb{N}\setminus\{D\},
$$
and
\begin{equation}
\label{eq:prior-model-lower}
Y_D=B_D\theta_{D}, \quad X_D=B_D+\sigma \eta_D,
\end{equation}
where $B_D$ is Gaussian random variable with mean $b_D$ and variance $\sigma^2$, that is independent of the standard Gaussian sequence $\{\eta_j\}_{j \in \mathbb{N}}$. Note that 
\begin{eqnarray}
\pi(\tilde \Theta_{a,\theta_0}(r_\sigma,b))&:=& \pi (\Theta_{a,\theta_0}(r_\sigma) \times \mathcal{B}(b)) \nonumber \\
&=&\pi_1(\Theta_{a,\theta_0}(r_\sigma)) \times \pi_2(\mathcal{B}(b)) \nonumber\\
&=&G_D^{-1}(C_0,C_1)  \frac{1}{\sigma \sqrt{2\pi}} \int_{C_0b_D}^{C_1b_D}\exp\left\{-\frac{1}{2\sigma^2}(t-b_D)^2 \right\}dt\\
&=& 1.
\end{eqnarray}

In view of the above, it is immediately seen that 
$$
L_\pi(Y,X)=\prod_{j \in \mathbb{N}}L_{\pi_j}(Y_j,X_j)=L_{\pi_D}(Y_D,X_D).
$$
Hence, as before, we arrive at 
\begin{equation}
\label{eqn:lower-bound-2nd-kep}
\beta_{0,\sigma,b}(\Theta(r_\sigma),\mathcal{B}(b)) \geq 1-\alpha -\frac{1}{2}(\mathbb{E}_0[L^2_{\pi_D}(Y_D,X_D)]-1)^{1/2}.
\end{equation}
Our task below is then to provide an upper bound on $\mathbb{E}_0[L^2_{\pi_D}(Y_D,X_D)]$. To this end, it is easily seen from model (\ref{eq:prior-model-lower}) that $Z_D=(X_D,Y_D)$, $D \in \mathbb{N}$, is Gaussian random vector with mean $U_{\theta,D}$ and covariance matrix $\sigma^2 \Sigma_{\theta,D}$, where
$$
U_{\theta,D}=\left(\begin{array}{c}
b_D\\
b_D \theta_D\\
\end{array}
\right),
\quad 
\Sigma_{\theta,D}=\left(\begin{array}{cc}
2&\theta_D\\
\theta_D&\theta^2_D\\
\end{array}
\right).
$$
Note that
$$
\Sigma^{-1}_{\theta,D}=\frac{1}{\theta_D^2}\left(\begin{array}{cc}
\theta_D^2&-\theta_D\\
-\theta_D&2\\
\end{array}
\right), 
$$
and
\begin{eqnarray}
(Z_D-U_{\theta,D})' \Sigma^{-1}_{\theta,D} (Z_D-U_{\theta,D})&=&\frac{1}{\theta_D^2}
(X_D-b_D, Y_D-b_D\theta_D) 
\left(\begin{array}{cc}
\theta_D^2&-\theta_D\\
-\theta_D&2\\
\end{array}
\right)
\left(\begin{array}{c}
X_D-b_D\\
Y_D-b_D \theta_D\\
\end{array}
\right)\nonumber \\
&=& 
\frac{1}{\theta_D^2}
(X_D-b_D, Y_D-b_D\theta_D) 
\left(\begin{array}{c}
\theta_D^2(X_D-b_D)-\theta_D(Y_D-b_D\theta_D)\\
-\theta_D(X_D-b_D) +2(Y_D-b_D\theta_D)\\
\end{array}
\right) \nonumber \\
&=& 
\frac{1}{\theta_D^2}
(X_D-b_D, Y_D-b_D\theta_D) 
\left(\begin{array}{c}
X_D\theta_D^2-Y_D\theta_D\\
2Y_D-b_D\theta_D -X_D\theta_D\\
\end{array}
\right) \nonumber \\
&=&
\frac{1}{\theta_D^2}((X_D-b_D)(X_D\theta_D^2-Y_D\theta_D)+(Y_D-b_D\theta_D)(2Y_D-b_D\theta_D -X_D\theta_D))\nonumber \\
&=&
\frac{1}{\theta_D^2}[(Y_D-X_D\theta_D)^2+(Y_D-b_D\theta_D)^2]. \nonumber
\end{eqnarray}
Hence,
\begin{eqnarray}
L_{\pi_D}(Z_D)&=&\exp\left\{(Z_D-U_{\theta_0,D})' \sigma^{-2}\Sigma^{-1}_{\theta_0,D} (Z_D-U_{\theta_0,D})-(Z_D-U_{\theta,D})' \sigma^{-2}\Sigma^{-1}_{\theta,D} (Z_D-U_{\theta,D})\right\}\nonumber\\
&=& \exp\left\{\frac{1}{\sigma^2}\left[\frac{(Y_D-X_D\theta_{D,0})^2}{\theta_{D,0}^2}+\frac{(Y_D-b_D\theta_{D,0})^2}{\theta_{D,0}^2}  
-\frac{(Y_D-X_D\theta_{D})^2}{\theta_{D}^2}-\frac{(Y_D-b_D\theta_{D})^2}{\theta_{D}^2}
\right]\right\} \nonumber\\
&=&
\exp\left\{\frac{1}{\sigma^2}\left[2Y_D^2\left(\frac{1}{\theta_{D,0}^2}-\frac{1}{\theta_{D}^2}\right)-2X_DY_D \left(\frac{1}{\theta_{D,0}}-\frac{1}{\theta_{D}}\right) -2b_DY_D \left(\frac{1}{\theta_{D,0}}-\frac{1}{\theta_{D}}\right) \right]\right\} \nonumber\\
&=&
\exp\left\{\frac{1}{\sigma^2}\left[2Y_D^2\left(\frac{1}{\theta_{D,0}^2}-\frac{1}{\theta_{D}^2}\right)-2Y_D(X_D+b_D)\left(\frac{1}{\theta_{D,0}}-\frac{1}{\theta_{D}}\right) \right]\right\}, \nonumber
\end{eqnarray}
Under $H_0$, $Y_D=B_D\theta_{D,0}$. Therefore, conditionally on $B_D$,
\begin{eqnarray}
\mathbb{E}_{0} \left[L^2_\pi(Z_D)\right]&=&G_D^{-1}(C_0,C_1)\,\mathbb{E} \left[ 
\exp\left\{\frac{2}{\sigma^2}\left(2B_D^2\left(1-\frac{\theta^2_{D,0}}{\theta_{D}^2}\right)-2B_D(X_D+b_D)\left(1-\frac{\theta_{D,0}}{\theta_{D}}\right) \right)\right\}
\right] \nonumber \\
&=&
G_D^{-1}(C_0,C_1) \mathbb{E} \left(\exp\left\{\frac{4B_D^2}{\sigma^2}\left(1-\frac{\theta^2_{D,0}}{\theta_{D}^2}\right)-\frac{4B_Db_D}{\sigma^2}\left(1-\frac{\theta_{D,0}}{\theta_{D}}\right) \right\} \right.\nonumber \\
&& \left. \times\; \mathbb{E}\left[ \exp\left\{ -\frac{4B_DX_D}{\sigma^2}\left(1-\frac{\theta_{D,0}}{\theta_{D}}\right) \right\}  \mid B_D \right]\right) \nonumber\\
&:=& G_D^{-1}(C_0,C1) \, \mathbb{E}\left(\mathbb{E}_{0} \left[L^2_\pi(Z_D) \mid B_D\right]\right). \nonumber
\end{eqnarray}
Using the formula
\begin{equation}
\label{eqn:laplace-G}
\mathbb{E} \left[ \exp(-(\lambda_1+\lambda_2 V) \right]=\exp(-\lambda_1+\lambda_2^2/2), \quad \lambda_1, \lambda_2 \in \mathbb{R}, 
\end{equation}
for any standard Gaussian random variable $V$, with 
$$
\lambda_1=\frac{4B^2_D}{\sigma^2}\left(1-\frac{\theta_{D,0}}{\theta_{D}} \right), \quad 
\lambda_2=\frac{4B_D}{\sigma}\left(1-\frac{\theta_{D,0}}{\theta_{D}} \right),
$$
we arrive at
\begin{eqnarray}
\mathbb{E}_{0} \left[L^2_\pi(Z_D) \mid B_D\right] &=&
\exp\left\{\frac{4B_D^2}{\sigma^2}\left(1-\frac{\theta^2_{D,0}}{\theta_{D}^2}\right)-\frac{4B_D}{\sigma^2}\left(1-\frac{\theta_{D,0}}{\theta_{D}}\right) (B_D+b_D)\right\}\nonumber \\
&& \times\; \exp\left\{\frac{8B^2_D}{\sigma^2}\left(1-\frac{\theta_{D,0}}{\theta_{D}}\right)^2 \right\} \nonumber \\
&=&\exp\{\sigma^{-2}B_D^2 \left[4(1-\rho_D^2)-4(1-\rho)+8(1-\rho_D)^2  \right]\} \nonumber \\
&&\times\;
\exp\{-4\sigma^{-2}B_Db_D(1-\rho_D)\}, \nonumber
\end{eqnarray}
where
$$
\rho_D =\left(1-\frac{\theta_{D,0}}{\theta_D}\right).
$$
Using simple algebra, we get
\begin{eqnarray}
\mathbb{E}_{0} \left[L^2_\pi(Z_D) \mid B_D\right] &=&
\exp\left\{\frac{4B_D^2}{\sigma^2}(\rho_D-1)(\rho_D-2)-\frac{4B_Db_D}{\sigma^2}(1-\rho_D)\right\} \nonumber \\
&=& 
\exp\left\{\frac{4}{\sigma^2}(\rho_D-1)[B_D^2(\rho_D-2)+B_Db_D]\right\}. \nonumber
\end{eqnarray}
It is easily seen that
\begin{eqnarray}
B_D^2(\rho_D-2)+B_Db_D&=&(b_D+\sigma \tilde \eta_D)^2(\rho_D-2)+(b_D+\sigma \tilde\eta_D)b_D \nonumber\\
&=&
b_D^2(\rho_D-1)+2\sigma \tilde \eta_D b_D(\rho_D-3/2)+\sigma^2\tilde \eta^2_D(\rho_D-2), \nonumber
\end{eqnarray}
where $\{\tilde \eta_D\}_{D \in \mathbb{N}}$ is a sequence of independent standard Gaussian random variables. Therefore,
\begin{eqnarray}
\mathbb{E}_{0} \left[L^2_\pi(Z_D) \mid B_D\right] &=&
\exp\left\{\frac{4b_D^2}{\sigma^2}(1-\rho_D)^2\right\}\nonumber \\
&& \times\; 
\exp\left\{\frac{8}{\sigma}b_D\tilde\eta_D(\rho_D-1)\left(\rho_D-\frac{3}{2}\right)+ 4\tilde \eta_D^2(\rho_D-2)\right\}. \nonumber
\end{eqnarray}
Since $\rho_D \in ]1,2[$, then $4 \tilde\eta_D^2(\rho_D-2)<0$ and, hence,
\begin{eqnarray}
\mathbb{E}_{0} \left[L^2_\pi(Z_D) \mid B_D\right] &\leq&
\exp\left\{\frac{4b_D^2}{\sigma^2}(1-\rho_D)^2+ 
\frac{8}{\sigma}b_D\tilde\eta_D(\rho_D-1)\left(\rho_D-\frac{3}{2}\right)\right\}. \nonumber
\end{eqnarray}
Using (\ref{eqn:laplace-G}) with
$$
\lambda_1=0, \quad 
\lambda_2=\frac{8}{\sigma}b_D(\rho_D-1)\left(\rho_D-\frac{3}{2}\right),
$$
we get
\begin{eqnarray}
\mathbb{E}_{0} \left[L^2_\pi(Z_D) \right] &=& G_D^{-1}(C_0,C_1)\,
\mathbb{E} \left(\mathbb{E}_{0} \left[L^2_\pi(Z_D) \mid B_D \right]\right)
\nonumber\\
&\leq& G_D^{-1}(C_0,C_1)\, \mathbb{E}\left( \exp\left\{\frac{4b_D^2}{\sigma^2}(1-\rho_D)^2+ 
\frac{8}{\sigma}b_D\tilde\eta_D(\rho_D-1)\left(\rho_D-\frac{3}{2}\right)\right\} \right)\nonumber\\
&=&
G_D^{-1}(C_0,C_1)\,
\exp\left\{\frac{4b_D^2}{\sigma^2}(1-\rho_D)^2 \left[1+\left(\rho_D-\frac{3}{2}\right)^2 \right]
\right\}\nonumber \\
&\leq& 
G_D^{-1}(C_0,C_1)\,
\exp\left\{\frac{5b_D^2}{\sigma^2}(1-\rho_D)^2 \right\} \nonumber\\
&\leq& 1+4(1-\alpha-\beta)^2, \nonumber
\end{eqnarray}
as soon as
$$
\frac{5b_D^2}{\sigma^2}(1-\rho_D)^2 \leq \ln (1+4(1-\alpha-\beta)^2) + \ln (G(C_0,C_1)),
$$
or, equivalently, as soon as 
$$
|\theta_D-\theta_{D,0}| \leq C_{\alpha,\beta,D}\; \sigma |\theta_D| b_D^{-1},
$$
where 
\begin{eqnarray}
C_{\alpha,\beta,D}&=&\ln (1+4(1-\alpha-\beta)^2) + \ln (G_D(C_0,C_1)) \nonumber \\
&:=&C_{\alpha,\beta}+\ln (G_D(C_0,C_1)).
\end{eqnarray}
(Note that, according to (\ref{eq:M2-def}), for all $D \leq M_2$, $C_{\alpha,\beta,D} \geq C_{\alpha,\beta}/2$.)\\

\noindent
{\bf Choice of $\theta$:} The sequence $\theta=(\theta_j)_{j \in \mathbb{N}}$ is chosen as follows
\begin{equation*}
\theta_j=
\begin{cases}
0 & \text{if $j\neq D$,}\\
a_D^{-1}/2 & \text{if $j= D$.}
\end{cases}
\end{equation*}
It can be easily seen that $\theta \in \mathcal{E}_a$.\\

\noindent
{\bf Choice of $\theta_0$:} The sequence $\theta_0=(\theta_{j,0})_{j \in \mathbb{N}}$ is chosen as follows
\begin{equation*}
\theta_{j,0}=
\begin{cases}
0 & \text{if $j\neq D$,}\\
a_D^{-1}/2 +C_{\alpha,\beta,D}\; \sigma a_D^{-1} b_D^{-1}/2& \text{if $j= D$.}
\end{cases}
\end{equation*}
Note that $\theta_0 \in  \mathcal{E}_a$ as soon as
\begin{equation}
\label{theta-0-cond}
C_{\alpha,\beta,D}\; \sigma b_D^{-1} \leq 1.
\end{equation}
Indeed, using the standard inequality $(x+y)^2 \leq 2(x^2+y^2)$, for $x, y \in \mathbb{R}$, we immediately get
\begin{eqnarray}
\sum_{j \in \mathbb{N}} a_j^2\theta_{j,0}^2=a_D^2\theta_{D,0}^2 \leq
a_D^2 \left(2\frac{a_D^{-2}}{4}+2 C^2_{\alpha,\beta,D}\frac{\sigma^2}{b_D^2}\frac{a_D^{-2}}{4}\right) \leq \frac{1}{2}+C^2_{\alpha,\beta,D}\frac{\sigma^2}{2b_D^2} \leq 1,
\end{eqnarray}
as soon as (\ref{theta-0-cond}) is satisfied. Furthermore, as soon as (\ref{theta-0-cond}) is satisfied, it is easily seen that $\theta-\theta_0 \in \mathcal{E}_a$.\\

Moreover, for the specific choices of $\theta$ and $\theta_0$ given above, it is immediately seen that
$$
|\theta_D-\theta_{D,0}| = C_{\alpha,\beta,D}\; \sigma b_D^{-1} |\theta_D| \Leftrightarrow 
\|\theta-\theta_{0}\| = \frac{C_{\alpha,\beta,D}}{2}\; \sigma b_D^{-1} a_D^{-1}.
$$
In other words, we have proved that for all $D \in \mathbb{N}$ satisfying (\ref{theta-0-cond}) then 
$$
\beta_{0,\sigma,b}(\Theta(r_{\sigma,D}),\mathcal{B}(b)) > \beta, \quad \text{where} \quad r_{\sigma,D}=\frac{C_{\alpha,\beta,D}}{2}\; \sigma b_D^{-1} a_D^{-1},
$$
for any given $\beta \in ]0,1-\alpha[$. This implies that, for every $\rho >0$,
$\beta_{0,\sigma,b}(\Theta(\rho),\mathcal{B}(b)) > \beta$ as soon as
$$
\rho \leq \frac{C_{\alpha,\beta,D}}{2}\; \sigma b_D^{-1} a_D^{-1}
\quad \text{for some} \quad D \in \mathbb{N}:\; C_{\alpha,\beta,D}\; \sigma b_D^{-1} \leq 1,
$$
which holds, as soon as
$$
\rho \leq \frac{C_{\alpha,\beta}}{4}\; \sigma b_D^{-1} a_D^{-1}
\quad \text{for some} \quad 1 \leq D \leq M_2,
$$
on noting that
$$
M_2:=\sup \left\{ D \in \mathbb{N}:\; C_{\alpha,\beta}\, \sigma |b_D^{-1}| \leq 2 \quad \text{and} \quad G_D(C_0,C1) \geq \frac{1}{\sqrt{1+4(1-\alpha-\beta)^2}}\right\},
$$
and that 
$$
C_{\alpha,\beta,D} \geq \frac{C_{\alpha,\beta}}{2}, \quad 1 \leq D \leq M_2.
$$

In particular, 
$$
\beta_{0,\sigma,b}(\Theta(\rho),\mathcal{B}(b)) > \beta \quad \text{for all} \quad \rho \leq\frac{C_{\alpha,\beta}}{4}\; \sigma \max_{1 \leq D \leq M_2} [ b_D^{-1} a_D^{-1} ].
$$
Hence, 
$$
\tilde r_{0,\sigma} \geq \frac{C_{\alpha,\beta}}{4}\; \sigma\, 
\max_{1 \leq D \leq M_2} [b_{D}^{-1} a_{D}^{-1}].
$$
This completes the proof of the proposition.
\begin{flushright} $\Box$ \end{flushright}

\subsubsection{Proof of Theorem \ref{thm-comb-lb}} 
\label{s:prlower_m}
The proof is splitted in two parts. We first show that $\tilde r_{\varepsilon,\sigma} \geq \tilde r_{\varepsilon,0}$ and then show that $\tilde r_{\varepsilon,\sigma} \geq \tilde r_{0,\sigma}$.\\

Consider observations $Y=(Y_j)_{j \in \mathbb{N}}$ from the GSM (\ref{eq:model-s=0}). Introduce the following goodness-of-fit testing algorithm:
\begin{itemize}
\item Generate a sequence $\tilde X=(\tilde X_j)_{j \in \mathbb{N}}$ according to
the GSM 
\begin{equation}
\label{eq:ind-copy-X}
\tilde X_j = b_j + \sigma \tilde \eta_j, \quad j \in \mathbb{N},
\end{equation}
where $\tilde \eta=(\tilde \eta_j)_{j \in \mathbb{N}}$ is a sequence of independent standard Gaussian random variables (that is also independent of the sequence $\xi=(\xi_j)_{j \in \mathbb{N}}$). (Note that the GSM (\ref{eq:ind-copy-X}) is an independent copy of the second equation in the GSM (\ref{1.0.0}).)

\item Let $\tilde \Psi_{\alpha}:=\tilde \Psi_{\alpha}(Y,\tilde X)$ be a given (non-randomized) $\alpha$-level goodness-of-fit testing procedure based on observations $(Y,\tilde X)=(Y_j, \tilde X_j)_{j \in \mathbb{N}}$ from the GSMs (\ref{eq:model-s=0}) and (\ref{eq:ind-copy-X}).
\item Define the randomized test $\Psi_{\alpha}:=\Psi_{\alpha}(Y)$\footnote{a measurable function 
of the observation $Y=(Y_j)_{j\in\mathbb{N}}$ from the GSM (\ref{eq:model-s=0})
with values in the interval $[0,1]$: for any given radius $\rho >0$, the null hypothesis is rejected with probability $\Psi_{\alpha}(Y)$ and it is not rejected with probability $1-\Psi_{\alpha}(Y)$. In this case, $\Alpha_\varepsilon(\Psi_{\alpha}):= \mathbb{E}_{\theta_0,b}( \Psi_{\alpha}(Y))$ and  $\Beta_\varepsilon(\Theta_a(\rho),\Psi_{\alpha}) := \sup_{\theta_0 \in \mathcal{E}_a\, \theta-\theta_0\in \Theta_a(\rho)} \mathbb{E}_{\theta,b}(1-\Psi_{\alpha}(Y)))$.} as
$$
\Psi_{\alpha}(Y):=\mathbb{E}[\tilde \Psi_{\alpha} \mid Y ],
$$
where $\mathbb{E}[\cdot]$ refers to expectation with respect to the independent standard Gaussian sequence $\tilde \eta$. 
\end{itemize}

In particular, for every $\varepsilon >0$ and $\sigma >0$, the randomized test $\Psi_{\alpha}$ is an $\alpha$-level test. Indeed,
\begin{eqnarray}
\Alpha_\varepsilon(\Psi_{\alpha}) &=& \mathbb{E}_{\theta_0,b}[\Psi_\alpha] \nonumber\\
&=&
\mathbb{E}_{\theta_0,b}[ \mathbb{E}[\tilde \Psi_{\alpha} \mid Y ] ] \nonumber \\
&=&  \mathbb{E}_{\theta_0,b}[\tilde \Psi_{\alpha}]  \nonumber \\
&=& \mathbb{P}_{\theta_0,b}(\tilde \Psi_{\alpha}=1) = \alpha,
\end{eqnarray}
since $\tilde \Psi_{\alpha}$ is an $\alpha$-level test.\\

Let $\theta \in l^2(\mathbb{N})$ and $\theta -\theta_0 \in \mathcal{E}_a$ be fixed. Then, the associated second kind error probability satisfies
\begin{eqnarray}
\mathbb{E}_{\theta,b}(1-\Psi_{\alpha}(Y)) &=& \mathbb{E}_{\theta,b}(1-\mathbb{E}[\tilde \Psi_{\alpha} \mid Y ]) \nonumber\\
&=& \mathbb{E}_{\theta,b}(1-\tilde \Psi_{\alpha}) \nonumber\\
&=& 
\mathbb{P}_{\theta,b}(\tilde \Psi_{\alpha}=0) \leq \beta,
\end{eqnarray}
as soon as
$$
\|\theta-\theta_0 \| \geq r_{\epsilon,\sigma}(\mathcal{E}_a,\tilde \Psi_\alpha, \beta).
$$

This implies that for any $\alpha$-level goodness-of-fit testing procedure $\tilde \Psi_\alpha$, based on observations $(Y,\tilde X)$ from the GSMs (\ref{eq:model-s=0})-(\ref{eq:ind-copy-X}), we can associate an $\alpha$-level goodness-of-fit testing procedure $\Psi_\alpha$, based on observations $Y$ from the GSM (\ref{eq:model-s=0}), such that the separation radius of $\Psi_\alpha$ is smaller than the separation radius of $\tilde \Psi_\alpha$, i.e.,  
$$
r_{\varepsilon,0}(\mathcal{E}_a,\Psi_\alpha, \beta) \leq r_{\varepsilon,\sigma}(\mathcal{E}_a, \tilde \Psi_\alpha, \beta).
$$
Hence, it is immediately seen that, for any $\alpha$-level goodness-of-fit testing procedure $\tilde \Psi_\alpha$, based on observations $(Y,\tilde X)$ from the GSMs (\ref{eq:model-s=0}) and (\ref{eq:ind-copy-X}),
\begin{eqnarray}
\tilde r_{\varepsilon,0}&:= &\inf_{\Psi_\alpha: \,\Alpha_{\varepsilon,0}( \bar \Psi_\alpha)\leq \alpha} r_{\varepsilon,0}(\mathcal{E}_a, \bar \Psi_\alpha,\beta) \nonumber\\
&\leq& r_{\varepsilon,0}(\mathcal{E}_a,\Psi_\alpha, \beta) \nonumber\\
&\leq& r_{\varepsilon,\sigma}(\mathcal{E}_a, \tilde \Psi_\alpha, \beta),
\end{eqnarray}
implying that
$$
\tilde r_{\varepsilon,0} \leq \tilde r_{\varepsilon,\sigma}.
$$

The proof of the assertion
$$
\tilde r_{0,\sigma} \leq \tilde r_{\varepsilon,\sigma}.
$$
follows similarly, along the lines of the proof of the previous assertion, and it is therefore omitted. This completes the proof of (\ref{eq:com-lower-bound-1})\\

Finally, (\ref{eq:com-lower-bound-2}) follows immediately form (\ref{eq:com-lower-bound-1}), taking into account (\ref{eq:prop-e=0}) and (\ref{eq:prop-s=0}).
This completes the proof of the theorem.
\begin{flushright} $\Box$ \end{flushright}

\subsection{Lower Bounds: Specific Cases}
\label{s:p_rateslow}

For the sake of convenience, we give the proof of each item (i)-(iv) in Theorem \ref{thm:alower} in different sections.

\subsubsection{Case (i): Mildly ill-posed problems with ordinary smooth functions}
We assume that (\ref{eq:a-b-mildly}) holds true, i.e.,
$$
b_j \sim j^{-t}, \; t>0,  \quad \text{and} \quad a_j \sim j^s, \; s>0, \quad j \in \mathbb{N}.
$$

\begin{proposition}
\label{eq:prop-a-b-mildly-lower}
Assume that the sequences $b=(b_j)_{j \in \mathbb{N}}$ and $a=(a_j)_{j \in \mathbb{N}}$ are given by (\ref{eq:a-b-mildly}). Then, there exists $\varepsilon_0, \sigma_0 \in ]0,1[$ such that, for all $0 < \varepsilon \leq \varepsilon_0$ and $0 <\sigma \leq \sigma_0$, the minimax separation radius $\tilde r_{\varepsilon, \sigma}$ satisfies
\begin{equation}
\label{eq:prop-mildly-lower}
\tilde r^2_{\varepsilon, \sigma} \gtrsim \varepsilon^{\frac{4s}{2s+2t+1/2}}
\vee \sigma^{2 \left(\frac{s}{t} \wedge 1\right)}.
\end{equation}
\end{proposition}

\noindent 
\textsc{Proof of Proposition \ref{eq:prop-a-b-mildly-lower}} For the second term in (\ref{eq:com-lower-bound-2}), it is known that (see \cite{LLM_2012}, \cite{ISS_2012}),
$$
\sup_{D \in \mathbb{N}} \left[c_{\alpha,\beta}\,
\varepsilon^2 \sqrt{\sum_{j=1}^D b_j^{-4}} \wedge a_D^{-2} \right] \sim 
\varepsilon^{\frac{4s}{2s+2t+1/2}}.
$$

Consider now the first term in (\ref{eq:com-lower-bound-2}). If $s>t$, then the sequence $\{b_j^{-1} a_j^{-1}\}_{j \in \mathbb{N}}$ is non-increasing and, thus,
$$
\frac{C^2_{\alpha,\beta}}{16}\, \sigma^2\, 
\max_{1 \leq D \leq M_2} [ b_D^{-2} a_D^{-2} ] \sim \sigma^2.
$$
On the other hand, if $s \leq t$, then the sequence $\{b_j^{-1} a_j^{-1}\}_{j \in \mathbb{N}}$ is non-decreasing. Hence, thanks to (\ref{eq:M2-def}), 
$$
\sigma^2 \sim b_{M_2}^2 \Leftrightarrow M_2 \sim \sigma^{-1/t},
$$ 
and, thus,
$$
\frac{C^2_{\alpha,\beta}}{16}\, \sigma^2\, \max_{1 \leq D \leq M_2} [ b_D^{-2} a_D^{-2} ]   \sim \sigma^2 b_{M_2}^{-2} a_{M_2}^{-2} \sim a_{M_2}^{-2} \sim \sigma^{\frac{2s}{t}}.
$$
Combining the above terms, we arrive at (\ref{eq:prop-mildly-lower}). This completes the proof of the proposition.
\begin{flushright} $\Box$ \end{flushright}

\subsubsection{Case (ii): Mildly ill-posed problems with super smooth functions}
We assume that (\ref{eq:a-b-super-mildly}) holds true, i.e.,
$$
b_j \sim j^{-t}, \; t>0,  \quad \text{and} \quad a_j \sim \exp\{ js\}, \; s>0, \quad j \in \mathbb{N}.
$$

\begin{proposition}
\label{eq:prop-a-b-super-mildly-lower}
Assume that the sequences $b=(b_j)_{j \in \mathbb{N}}$ and $a=(a_j)_{j \in \mathbb{N}}$ are given by (\ref{eq:a-b-super-mildly}). Then, there exists $\varepsilon_0, \sigma_0 \in ]0,1[$ such that, for all $0 < \varepsilon \leq \varepsilon_0$ and $0 <\sigma \leq \sigma_0$, the minimax separation radius $\tilde r_{\varepsilon, \sigma}$ satisfies
\begin{equation}
\label{eq:prop-super-mildly-lower}
\tilde r^2_{\varepsilon, \sigma} \gtrsim \varepsilon^2 \left[\ln\left(1/\varepsilon\right)\right]^{\left(2t+\frac{1}{2}\right)} \vee \sigma^{2}.
\end{equation}
\end{proposition}

\noindent 
\textsc{Proof of Proposition \ref{eq:prop-a-b-super-mildly-lower}} For the second term in (\ref{eq:com-lower-bound-2}), it is known that (see \cite{LLM_2012}, \cite{ISS_2012}),
$$
\sup_{D \in \mathbb{N}} \left[c_{\alpha,\beta}\,
\varepsilon^2 \sqrt{\sum_{j=1}^D b_j^{-4}} \wedge a_D^{-2} \right] \sim 
\varepsilon^2 \left[\ln\left(1/\varepsilon\right)\right]^{\left(2t+\frac{1}{2}\right)}.
$$

Consider now the first term in (\ref{eq:com-lower-bound-2}). Then, the sequence $\{b_j^{-1} a_j^{-1}\}_{j \in \mathbb{N}}$ is non-increasing for each $s,t>0$, and, thus,
$$
\frac{C^2_{\alpha,\beta}}{16}\, \sigma^2\, 
\max_{1 \leq D \leq M_2} [ b_D^{-2} a_D^{-2} ] \sim \sigma^2.
$$
Combining the above terms, we arrive at (\ref{eq:prop-super-mildly-lower}). This completes the proof of the proposition.
\begin{flushright} $\Box$ \end{flushright}

\subsubsection{Case (iii): Severely ill-posed problems with ordinary smooth functions}
We assume that (\ref{eq:a-b-mildly-super}) holds true, i.e.,
$$
b_j \sim \exp\{ -jt\}, \; t>0,  \quad \text{and} \quad a_j \sim j^s, \; s>0, \quad j \in \mathbb{N}.
$$

\begin{proposition}
\label{eq:prop-a-b-mildly-super-lower}
Assume that the sequences $b=(b_j)_{j \in \mathbb{N}}$ and $a=(a_j)_{j \in \mathbb{N}}$ are given by (\ref{eq:a-b-mildly-super}). Then, there exists $\varepsilon_0, \sigma_0 \in ]0,1[$ such that, for all $0 < \varepsilon \leq \varepsilon_0$ and $0 <\sigma \leq \sigma_0$, the minimax separation radius $\tilde r_{\varepsilon, \sigma}$ satisfies
\begin{equation}
\label{eq:prop-mildly-super-lower}
\tilde r^2_{\varepsilon, \sigma} \gtrsim \left[\ln\left(1/\varepsilon\right)\right]^{-2s} \vee \left[\ln\left(1/\sigma\right)\right]^{-2s}.
\end{equation}
\end{proposition}

\noindent 
\textsc{Proof of Proposition \ref{eq:prop-a-b-mildly-super-lower}} For the second term in (\ref{eq:com-lower-bound-2}), it is known that (see \cite{LLM_2012}, \cite{ISS_2012}),
$$
\sup_{D \in \mathbb{N}} \left[c_{\alpha,\beta}\,
\varepsilon^2 \sqrt{\sum_{j=1}^D b_j^{-4}} \wedge a_D^{-2} \right] \sim 
\left[\ln\left(1/\varepsilon\right)\right]^{-2s}.
$$

Consider now the first term in (\ref{eq:com-lower-bound-2}). Then, the sequence $\{b_j^{-1} a_j^{-1}\}_{j \in \mathbb{N}}$ is non-decreasing. Hence, thanks to (\ref{eq:M2-def}), 
$$
\sigma^2 \sim b_{M_2}^2 \Leftrightarrow M_2 \sim \frac{1}{t} \ln\left(1/\sigma\right)
$$ 
and, thus,
$$
\frac{C^2_{\alpha,\beta}}{16}\, \sigma^2\, \max_{1 \leq D \leq M_2} [ b_D^{-2} a_D^{-2} ]   \sim \sigma^2 b_{M_2}^{-2} a_{M_2}^{-2} \sim a_{M_2}^{-2} \sim \left[\ln\left(1/\sigma\right)\right]^{-2s}.
$$
Combining the above terms, we arrive at (\ref{eq:prop-mildly-super-lower}). This completes the proof of the proposition.
\begin{flushright} $\Box$ \end{flushright}

\subsubsection{Case (iv): Severely ill-posed problems with super smooth functions}
We assume that (\ref{eq:a-b-super}) holds true, i.e.,
$$
b_j \sim \exp\{-jt\}, \; t>0,  \quad \text{and} \quad a_j \sim \exp\{js\}, \; s>0, \quad j \in \mathbb{N}.
$$

\begin{proposition}
\label{eq:prop-a-b-super-lower}
Assume that the sequences $b=(b_j)_{j \in \mathbb{N}}$ and $a=(a_j)_{j \in \mathbb{N}}$ are given by (\ref{eq:a-b-super}). Then, there exists $\varepsilon_0, \sigma_0 \in ]0,1[$ such that, for all $0 < \varepsilon \leq \varepsilon_0$ and $0 <\sigma \leq \sigma_0$, the minimax separation radius $\tilde r_{\varepsilon, \sigma}$ satisfies
\begin{equation}
\label{eq:prop-super-lower}
\tilde r^2_{\varepsilon, \sigma} \gtrsim \varepsilon^{\frac{2s}{s+t}}
\vee \sigma^{2 \left(\frac{s}{t} \wedge 1\right)}.
\end{equation}
\end{proposition}

\noindent 
\textsc{Proof of Proposition \ref{eq:prop-a-b-super-lower}} For the second term in (\ref{eq:com-lower-bound-2}), it is known that (see \cite{LLM_2012}, \cite{ISS_2012}),
$$
\sup_{D \in \mathbb{N}} \left[c_{\alpha,\beta}\,
\varepsilon^2 \sqrt{\sum_{j=1}^D b_j^{-4}} \wedge a_D^{-2} \right] \sim 
\varepsilon^{\frac{2s}{s+t}}.
$$

Consider now the first term in (\ref{eq:com-lower-bound-2}). If $s>t$, then the sequence $\{b_j^{-1} a_j^{-1}\}_{j \in \mathbb{N}}$ is non-increasing and, thus,
$$
\frac{C^2_{\alpha,\beta}}{16}\, \sigma^2\, 
\max_{1 \leq D \leq M_2} [ b_D^{-2} a_D^{-2} ] \sim \sigma^2.
$$
On the other hand, if $s \leq t$, then the sequence $\{b_j^{-1} a_j^{-1}\}_{j \in \mathbb{N}}$ is non-decreasing. Hence, thanks to (\ref{eq:M2-def}), 
$$
\sigma^2 \sim b_{M_2}^2 \Leftrightarrow M_2 \sim \frac{1}{t} \ln\left(1/\sigma\right),
$$ 
and, thus,
$$
\frac{C^2_{\alpha,\beta}}{16}\, \sigma^2\, \max_{1 \leq D \leq M_2} [ b_D^{-2} a_D^{-2} ]   \sim \sigma^2 b_{M_2}^{-2} a_{M_2}^{-2} \sim a_{M_2}^{-2} \sim \sigma^{\frac{2s}{t}}.
$$
Combining the above terms, we arrive at (\ref{eq:prop-mildly-lower}). This completes the proof of the proposition.
\begin{flushright} $\Box$ \end{flushright}

\bigskip
\bigskip
\bibliography{Survey}
\bibliographystyle{plain}

\end{document}